\input amssym.def
\input amssym
\font\goth=eufm10

\font\Big=cmbx10 scaled\magstep2
\def\g{{\hbox{\goth g}}}

\newdimen\ybox\ybox=5pt 
\setbox249=\vbox{\hbox{\vrule{\vbox to\ybox{}}\kern\ybox}\hrule}
\setbox248=\hbox{\copy249}

\newcount\n
\def\R#1,{%
\n=#1%
\nointerlineskip
\hbox{%
\loop\ifnum\n>0\copy248\advance\n by-1\repeat%
\vrule}%
}

\def\YD#1.{\,\vcenter{\hrule#1}\,}
\def\half{{1\over 2}}
\def\mapright#1{\smash{\mathop{\longrightarrow}\limits^{#1}}}

\Big \centerline{OPERATOR ALGEBRAS AND CONFORMAL FIELD THEORY III}
\vskip .05in
\bf \centerline{FUSION OF POSITIVE ENERGY REPRESENTATIONS
OF LSU(N)} 
\bf \centerline{USING BOUNDED OPERATORS}
\vskip .05in
\bf \centerline{Antony Wassermann}
\bf \centerline{University of Cambridge}
\vskip .1in

\def\ue{{\bf\"{\goth u}}}
\def\oe{{\bf\"{\goth o}}}
\def\s{{\rm\ss}}
\goth\centerline{Schm\ue cke dich, o liebe Seele,}

 \centerline{La\s e die dunkle S\ue ndenh\oe hle,}

 \centerline{Komm ins helle Licht gegangen,}

  \centerline{Fange herrlich an zu prangen!}
\vskip .02in

\centerline{J~Franck~1649}

\vskip .2in
\bf \centerline{CONTENTS}
\vskip .2in
\rm
{\baselineskip = 10pt
\noindent 1. Introduction.
\vskip .05in
\noindent CHAPTER I.~POSITIVE ENERGY REPRESENTATIONS OF LSU(N).
\vskip .05in
\noindent 2. Irreducible representations of SU(N).

\noindent 3. Fermions and quantisation.

\noindent 4. The fundamental representation.

\noindent 5. The central extension ${\cal L}G$.

\noindent 6. Positive energy representations.

\noindent 7. Infinitesimal action of $L^0\g$ on finite energy vectors.

\noindent 8. The exponentiation theorem.

\noindent 9. Classification of positive energy representations of
level $\ell$.

\vskip .05in
\noindent CHAPTER II.~LOCAL LOOP GROUPS AND THEIR VON NEUMANN ALGEBRAS.
\vskip .05in
\noindent 10. von Neumann algebras. 

\noindent 11. Abstract modular theory. 

\noindent 12. Modular operators and Takesaki devissage for von
Neumann algebras. 

\noindent 13. Araki duality and modular theory for Clifford algebras.

\noindent 14. Prequantised geometric modular theory. 

\noindent 15. Haag--Araki duality and geometric modular theory
for fermions on the circle.

\noindent 16. Ergodicity of the modular group.

\noindent 17. Consequences of modular theory for local loop groups.

\vskip .05in
\noindent CHAPTER III.~THE BASIC ORDINARY DIFFERENTIAL EQUATION.
\vskip .05in
\noindent 18. The basic ODE and the transport problem.

\noindent 19. Analytic reduction of the ODE.

\noindent 20. Algebraic reduction of the ODE.

\noindent 21. Symmetry and analyticity properties of transport matrices.

\noindent 22. Projected power series solutions.

\noindent 23. Euler--Thomae integral representation of projected
solutions. 

\noindent 24. Computation of transport matrices.
\vskip .05in

\noindent CHAPTER IV.~VECTOR AND DUAL VECTOR PRIMARY FIELDS.
\vskip .05in
\noindent 25. Existence and uniqueness of vector and dual vector
primary fields. 

\noindent 26. Transport equations for four--point functions and
braiding of primary fields.

\noindent 27. Sugawara's formula.

\noindent 28. The Knizhnik--Zamolodchikov ODE.

\noindent 29. Braiding relations between vector and dual vector
primary fields.
\vskip .05in
\noindent CHAPTER V.~CONNES FUSION OF POSITIVE ENERGY REPRESENTATIONS.
\vskip .05in
\noindent 30. Definition and elementary properties of Connes fusion
for positive energy representations.

\noindent 31. Connes fusion with the vector
representation.

\noindent 32. Connes fusion with exterior powers of the vector
representation.

\noindent 33. The fusion ring. 

\noindent 34. General fusion rules (Verlinde formulas).
\vskip .05in
\noindent References.} 
\vfill\eject
\noindent \bf 1. Introduction. \rm This is one of a series
of papers devoted to the study of
conformal field theory from the point of view of operator algebras
(see [39] and [40] for an overview of the whole series). 
In order to make the paper accessible and self--contained, we have
not assumed a detailed knowledge of either operator algebras or
conformal field theory, inlcuding short--cuts and direct proofs
wherever possible. This research
programme was originally motivated by V.~Jones' suggestion that there
might be 
a deeper `operator algebraic' explanation of the coincidence between
certain unitary representations of the infinite braid group that had
turned up independently in the theory of
subfactors, exactly solvable models in statistical mechanics and
conformal field theory (CFT). To understand why there should be any link
between these subjects, recall that, amongst other things,
the classical `additive' theory of von Neumann algebras [25] was
developed to provide a framework for studying unitary 
representations of Lie groups. In concrete examples,
for example the Plancherel theorem for semisimple groups, this
abstract framework had to be complemented by a considerably harder
analysis of intertwining operators and associated differential equations.
The link between CFT and operator algebras comes from the recently
developed `multiplicative' (quantum?) theory of von Neumann algebras. 
This theory has three
basic sources: firstly the algebraic approach to quantum field theory (QFT)
of Doplicher, Haag and Roberts [9]; then in Connes' theory of
bimodules and their tensor products or fusion [8]; and lastly in Jones' theory
of subfactors [17]. Our work reconciles these ideas with the theory of
primary fields, one of the fundamental concepts in CFT.
Our work has the following consequences, some of which will be taken
up in subsequent papers:
\vskip .05in
\item{(1)} Several new constructions of subfactors.
\item{(2)} Non--trivial algebraic QFT's
in $1+1$ dimensions with finitely many sectors and non--integer
statistical (or quantum) dimension (``algebraic CFT''). 
\item{(3)} A definition of quantum invariant theory without using quantum
groups at roots of unity.
\item{(4)} A computable and manifestly unitary definition of
fusion for positive energy representations (``Connes fusion'') making
them into a tensor category. 
\item{(5)} Analytic properties of primary fields (``constructive CFT'').
\vskip .05in
\noindent To our knowledge, no previous work has suceeded in
integrating the theory of primary fields with the ideas of algebraic QFT
nor in revealing the very simple analytic structure of primary fields.
As we explain below, the main thrust of our work is the explicit
computation of Connes
fusion of positive energy representations. Finiteness of statistical
dimension (or Jones index) is a natural consequence, not a technical
mathematical inconvenience. It is perhaps worth
emphasising that the theory of operator algebras only provides a framework for
studying CFT. As in the case of group representations, it must be
complemented by a detailed analysis of certain interwining operators,
the primary fields, and their associated differential equations. As we
discuss later, however, the operator algebraic point of view can be
used to reveal basic positivity and unitarity properties in
CFT that previously seem to have been overlooked.

Novel features of our treatment are the construction of
representations and primary fields from fermions. This makes unitarity
of the representations and boundedness properties of smeared vector
primary fields obvious. The only formal ``vertex algebra'' aspects 
of primary fields borrowed from [37] are the trivial proof of
uniqueness and the statement of the Knizhnik--Zamolodchikov equation;
our short derivation of the KZ equation circumvents the well--known
contour integral proof [22], implicit 
but not given in [37]. The proof that the axioms of algebraic QFT are
satisfied in the non--vacuum sectors is new and relies heavily on our
fermionic construction; the easier properties in the vacuum sector have been
known for some time [14]. The treatment of braiding relations for
smeared primary fields is new but inspired by the
Bargmann--Hall--Wightman theorem ([19],[34]). To our knowledge, the
application of Connes 
fusion to a non--trivial model in QFT is quite new. Our definition is
a slightly simplified version of Connes' original
definition, tailor--made for CFT because of the ``four--point function
formula''; no general theory is required.

The finite--dimensional 
irreducible unitary representations of $SU(N)$ and their tensor
product rules are well known to mathematicians and physicists. The
representations $V_f$ are classified by signatures or Young diagrams
$f_1\ge f_2\ge \cdots \ge f_N$ and, if $V_{[k]}=\lambda^k {\Bbb C}^N$,
we have the tensor product rule
$V_f\otimes V_{[k]}=\bigoplus_{g>_k f} V_g$, where $g$ ranges over
all diagrams that can be obtained by adding $k$ boxes to $f$ with no
two in the same row. For the infinite--dimensional loop group
$LSU(N)= C^\infty(S^1,SU(N))$, the appropriate unitary representations
to consider in place of finite--dimensional representations are the
projective unitary
representations of positive energy. Positive energy representations
form one of the most important foundation stones
of conformal field theory ([5], [11], [22]). The classification of
positive energy representations is straightforward and has been known
for some time now. A positive energy representation $H_f$ is
classified by its level $\ell$, a positive integer, and its signature
$f$, which must satisfy the permissibility condition $f_1-f_N\le
\ell$.  Extending 
the tensor product rules to representations of a fixed level, however,
presents a problem. It is already 
extremely difficult just giving a coherent definition of the 
tensor product, since the naive one fails hopelessly because it does
not preserve the level. On the other hand physicists have known for
years how to `fuse' 
representations in terms of short range expansions of products of
associated quantum fields (primary fields).
We provide one solution to this `problem of fusion' in conformal field
theory by giving a
mathematically sound definition of the tensor product
that ties up with the intuitive picture of physicists. Our solution
relates positive energy representations of loop groups to bimodules
over von Neumann algebras. Connes defined a tensor product operation
on such bimodules -- ``Connes fusion'' -- which translates directly
into a definition of fusion for positive energy representations. 
The general fusion rules follow from the particular rules $H_f\boxtimes
H_{[k]} =\bigoplus_{g>_k f} H_g$, where $g$ must now also be
permissible. In this way the level $\ell$ representations of $LSU(N)$
exhibit a structure similar to that of the irreducible representations
of a finite group. There are several other
approaches to fusion of positive energy 
representations, notably those of Graeme Segal [33] and Kazhdan \&
Lusztig [21]. Our picture seems to be a unitary boundary value
of Segal's holomorphic proposal for fusion, based on a disc
with two smaller discs removed. When the discs shrink to points on the
Riemann sphere, Segal's definition should degenerate to the algebraic
geometric fusion of Kazhdan \& Lusztig. We now give an informal
summary of the paper. 
\vskip .05in

\noindent \bf Fermions. \rm Let ${\rm
Cliff}(H)$  be the Clifford algebra of a complex Hilbert space $H$
is generated by a linear map $f\mapsto a(f)$ ($f\in H$) satisfying 
$a(f)a(g)+a(g)a(f)=0$ and $a(f)a(g)^*+a(g)^*a(f)=(f,g)$. It acts
irreducibly on Fock space $\Lambda H$ via $a(f) \omega=f\wedge 
\omega$. Other representations of ${\rm Cliff}(H)$ arise by
considering the real linear map $c(f)=a(f)+a(f)^*$ which satisfies 
$c(f)c(g)+c(g)c(f)=2{\rm Re}\,(f,g)$; note that $a(f)=\half
(c(f)-ic(if))$. Since $c$ relies only on the underlying real Hilbert
space $H_{\Bbb R}$, complex 
structures on $H_{\Bbb R}$ commuting with $i$ give new irreducible
representations of ${\rm Cliff}(H)$. The structures correspond to
projections $P$ with multiplication by $i$ given by $i$ on $PH$ and
$-i$ on $(PH)^\perp$. The corresponding representation $\pi_P$ is
given by $\pi_P(a(f))=\half(c(f) -i c(i(2P-I)f))$. Using ideas that go
back to Dirac and von Neumann, we give our own short proof of I.~Segal's
equivalence criterion: if $P-Q$ is a Hilbert--Schmidt operator, then
$\pi_P$ and $\pi_Q$ are unitarily equivalent. On the other hand if
$u\in U(H)$, then $a(uf)$ and $a(ug)$ also 
satisfy the complex Clifford algebra relations. Thus $a(f)\rightarrow
a(uf)$ gives an automorphism of ${\rm
Cliff}(H)$. We say that this ``Bogoliubov'' automorphism is implemented
in $\pi_P$ iff $\pi_P(a(uf))=U\pi_P(a(f))U^*$ for some unitary $U$.
This gives a projective representation of the subgroup of
implementable unitaries $U_P(H)$. 
Segal's equivalence criterion leads immediately to a
quantisation criterion: if $[u,P]$ is a Hilbert--Schmidt operator,
then $u\in U_P(H)$. 
\vskip .05in
\noindent \bf Positive energy representations. \rm Let $G=SU(N)$ and let
$LG=C^\infty(S^1,G)$ be the loop group, with the rotation group ${\rm
Rot}\, S^1$ acting as automorphisms. If $H=L^2(S^1,{\Bbb C}^N)$ and
$P$ is the projection onto Hardy 
space $H^2(S^1,{\Bbb C}^N)$,  $LSU(N)\rtimes {\rm Rot}\, S^1\subset
U_P(H)$ so we get a projective representations $\pi_P^{\otimes\ell}
:LU(N)\rtimes {\rm Rot}\,S^1 \rightarrow PU({\cal F}^{\otimes\ell})$
where ${\cal F}$ denotes Fock space $\Lambda H_P$. Now ${\rm Rot}\,
S^1$ acts with positive energy, where an action $U_\theta$ on $H$ is
said to have positive energy if $H=\bigoplus_{n\ge 0} H(n)$ with
$U_\theta\xi=e^{in\theta}\xi$ for $\xi\in H(n)$,
$H(n)$ is finite--dimensional and $H(0)\ne (0)$. This implies that
${\cal F}^{\otimes \ell}$ splits as a direct sum of irreducibles
$H_i$, called the level $\ell$ positive energy representations. 
The $H_i$'s are classified by their lowest energy subspaces $H_i(0)$,
which are irreducible modules for the constant loops $SU(N)$. Their
signatures $f_1\ge f_2\ge \cdots \ge f_N$ must satisfy $f_1-f_N\le \ell$, so
${\cal F}_V^{\otimes \ell}$ has only finitely many inequivalent
irreducible summands. This classification is achieved by defining
an infinitesimal action of the algebraic Lie
algebra $L^0\g$ on the 
finite energy subspace $H^0=\sum H(n)$ using bilinear terms
$a(f) a(g)^*$. Our main contribution here is to match up these
operators with the skew--adjoint operators predicted by analysis.
The quantisation criterion also implies that the M\"obius
transformations of determinant $1$ act projectively on each positive
energy representation compatibly with $LG$. The vacuum representation
$H_0$ corresponds to the trivial representation of $G$; the
M\"obius transformations of determinant $-1$ also act on $H_0$, but
this time by conjugate--linear isometries. This presentation of the
theory of positive energy representations is adequate for the needs of
this paper; in [40] we show from scratch that any irreducible positive
energy representation of $LSU(N)\rtimes {\rm Rot}\, S^1$ arises as a
subrepresentation of some ${\cal F}_V^{\otimes \ell}$.
\vskip .05in

\noindent \bf von Neumann algebras. \rm We briefly summarise those
parts of the general theory of operator algebras that are background
for this paper. (They will serve only as motivation, since all the
advanced results we need will be proved directly for local fermion or
loop group algebras.)  A {\it von Neumann algebra} is simply the commutant  
${\cal S}^\prime = \{T\in B(H): \hbox{$Tx=xT$ for all $x\in {\cal S}$}\}$ of
a subset ${\cal S}$ of $B(H)$ with ${\cal S}^*={\cal S}$. Typically
${\cal S}$ will be a *--subalgebra of $B(H)$ or a subgroup of $U(H)$;
the von Neumann algebra generated by ${\cal S}$ is then just ${\cal
S}^{\prime\prime}$. 
A von Neumann algebra $M$ is called a {\it factor} if its centre contains only
scalar operators. {\it Modules} over a factor were classified by Murray and
von Neumann [25] using a {\it dimension function}, the range of values
giving an invariant of the factor: the non--negative integers (type
I), the non--negative reals (type II) and $\{0,\infty\}$ (type
III). Further structure comes from the modular operators $\Delta^{it}$
and $J$ of Tomita--Takesaki [7]: if $\Omega$ is a cyclic vector for $M$ and
$M^\prime$ and $S=J\Delta^{1/2}$ is the polar decomposition of the map
$S:M\Omega\rightarrow M\Omega, a\Omega\mapsto a^*\Omega$, then
$JMJ=M^\prime$ and $\Delta^{it}M\Delta^{-it}=M$. 
On the one hand the operators $\Delta^{it}$ provide a further
invariant of type III factors, the Connes spectrum 
$\bigcap_\Omega {\rm Spec}\,\Delta_\Omega^{it}$, a closed subgroup
of ${\Bbb R}$ ([8]; see also [40]); while on the other hand $J$ makes
the underlying Hilbert space $H_0$ into a bimodule 
over $M$, the {\it vacuum bimodule}, with the action of the opposite algebra
$M^{{\rm op}}$ given by $a\mapsto Ja^*J$. Bimodules are
closely related to subfactors and endomorphisms: a bimodule 
defines a subfactor by the inclusion $M^{\rm op}\subset M^\prime$; and
an endomorphism $\rho:M\rightarrow M$ can be used to define a new
bimodule structure on $H_0$. {\it Connes fusion} [8] gives an
associative tensor 
product operation on bimodules that generalises composition of
endomorphisms: given bimodules $X$ and $Y$, their fusion $X\boxtimes Y$
is the completion of ${\rm Hom}_{M^{\rm op}}(H_0,X)\otimes 
{\rm Hom}_M(H_0,Y)$ with respect to the pre--inner product 
$\langle x_1\otimes y_1,x_2\otimes y_2\rangle
=(x_2^*x_1y_2^*y_1\Omega,\Omega)$. Roughly speaking Jones,
Ocneanu and Popa ([17], [18], [28], [40]) proved that an {\it irreducible
bimodule} is classified by the {\it tensor category} it generates under
fusion, provided the category contains only finitely many isomorphism classes
of irreducible bimodules.

\vskip .05in

\noindent \bf Modular theory for fermions. \rm For fermions and
bosons, modular theory provides the most convenient framework for
proving the much older result in algebraic quantum field theory known
as ``Haag--Araki duality''. This deals with the symmetry between
observables in a region and its (space--like) complement. As in [23],
we consider more generally a {\it modular
subspace} $K$ of a complex Hilbert space $H$, i.e.~a closed real subspace
such that $K\cap iK=(0)$ and $K+iK$ is dense in $H$. (Thus
$K=\overline{M_{\rm sa} \Omega}$ in Tomita--Takesaki theory.) If
$S=J\Delta^{1/2}$ is the polar decomposition of the map 
$S:K+iK\rightarrow K+iK, \xi+i\eta \mapsto \xi-i\eta$, then
$JK=iK^{\perp}$ and $\Delta^{it}K=K$; in the text following [31] we
avoid unbounded operators by taking the equivalent definitions $J={\rm
phase}(E-F)$ and $\Delta^{it}=(2I-E-F)^{it}(E+F)^{it}$, where $E$ and
$F$ are the projections onto $K$ and $iK$. The modular operators $J$
and $\Delta^{it}$ are uniquely characterised by the Kubo--Martin--Schwinger
(KMS) condition: commuting operators $J$ and $\Delta^{it}$ give the modular
operators if $\Delta^{it}K=K$ and, for each $\xi\in K$,
$f(t)=\Delta^{it}\xi$ extends to 
a bounded continuous function on the strip $-\half \le {\rm Im}\, z\le
0$, holomorphic in the interior, with $f(t-i/2) =Jf(t)$. 

This theory can be used to prove an abstract result, implicit in the
work of Araki ([1], [2]). Let $K$ be a modular subspace of $H$
and let $M(K)$ be the von Neumann algebra on $\Lambda H$
generated by the operators $c(\xi)$ for $\xi\in H$.
Then $M(K^\perp)$ is the graded commutant of $M(K)$ (``Araki
duality'') and the modular
operators for $M(K)$ on $\Lambda H$ come from the quantisations of the
corresponding operators for $K$. This reduces computations to
``one--particle states'', i.e.~the prequantised Hilbert space. We then
perform the prequantised computation explicitly when $H=L^2(S^1,V)$ and
$K=L^2(I,V)$, with $I$ a proper subinterval of $S^1$ with complement
$I^c$. We deduce that 
if $M(I)$ is the von Neumann algebra on $\Lambda H_P$ by $a(f)$'s with
$f\in L^2(I,V)$, then $M(I^c)$ is the graded commutant of $M(I)$
(Haag--Araki duality) $\Delta^{it}$ and $J$ come from the M\"obius
flow and flip fixing the end points of $I$.

\vskip .05in
\noindent \bf Local loop groups. \rm Let $L_ISU(N)$ be
the subgroup of $LSU(N)$ consisting of loops equal to 
$1$ off $I$. The von Neumann algebra $N(I)$ generated by $L_I G$ is a
subalgebra of the local fermion algebra $M(I)$ invariant under
conjugation by the modular group $\Delta^{it}$, since it is geometric.
The modular operators of $N(I)$ can therefore be read off from those
of $M(I)$ by a result in [35] (``Takesaki devissage''); we give our own
short proof of a slightly modified version of Takesaki's result.
We deduce the following properties of the local subgroups, predicted
by the Doplicher--Haag--Roberts axioms [9]. The use of devissage,
relating different models, is new and seems unavoidable in proving
factoriality and local equivalence.

\vskip .05in
\noindent 1. {\it Locality} In any positive energy representation
$L_ISU(N)$ and $L_{I^c}SU(N)$ commute. 
\vskip .02in

\noindent  2. {\it Factoriality.}
$\pi_i(L_ISU(N))^{\prime\prime}$ is a factor if $(\pi_i,H_i)$ is an
irreducible positive energy representation. 

\vskip .02in

\noindent  3. {\it Local equivalence.} There is a unique
*--isomorphism $\pi_i:\pi_0(L_I 
G)^{\prime\prime} \rightarrow \pi_i(L_IG)^{\prime\prime}$ sending 
such that $Ta=\pi_i(a) T$ for all $T\in {\rm Hom}_{L_IG}(H_0,H_i)$.
\vskip .02in
\noindent  4. {\it Haag duality.} If $\pi_0$ is the vacuum representation
at level $\ell$, then
$\pi_0(L_ISU(N))^{\prime\prime} =\pi_0(L_{I^c}SU(N))^\prime$.

\vskip .02in
\noindent  5. {\it Irreducibility.} Let $A$ be a finite subset of $S^1$
and let $L^ASU(N)$ be the subgroup of $LSU(N)$ of loops trivial to all
orders at points of $A$. If $\pi$ is positive energy, then
$\pi(L^ASU(N))^{\prime} =\pi(LSU(N))^\prime$, so the irreducible
positive energy representations of $LSU(N)$ stay irreducible and
inequivalent when restricted to $L^ASU(N)$. 

\vskip .05in
\noindent \bf Vector primary fields. \rm Let $P_i$ and $P_j$ be
projections onto the 
irreducible summands $H_i$ and $H_j$ of $\pi_P^{\otimes \ell}$ and
fix an $SU(N)$--equivariant embedding of ${\Bbb C}^N$ in 
${\Bbb C}^N\otimes {\Bbb C}^\ell$. If $f\in L^2(S^1,{\Bbb C}^N)\subset
L^2(S^1,{\Bbb C}^N\otimes{\Bbb C}^\ell)$, 
we may ``compress'' the smeared fermion field $a(f)$ to get an operator
$\phi_{ij}(f)=P_i a(f) P_j\in {\rm Hom}(H_j,H_i)$. By construction
$\phi_{ij}(f)$ satisfies a group covariance relation 
$g\phi(f)g^{-1}=\phi(g\cdot f)$
for $g\in LSU(N)\rtimes {\rm Rot}\, S^1$ as well as the $L^2$ bound
$\|\phi(f)\|\le \|f\|_2$.
If $f$ is supported in $I^c$, then $\phi(f)$ gives a concrete 
element in ${\rm Hom}_{L_ISU(N)}(H_j,H_i)$; this space of
intertwiners is known to be non--zero by local 
equivalence. Clearly $\phi$ defines a map $L^2(S^1,{\Bbb C}^N)\otimes
H_j\rightarrow H_i$ which intertwines the action of $LSU(N)\rtimes
{\rm Rot}\, S^1$. The modes $\phi(v,n)= \phi(z^{-n} v)$ satisfy
Lie algebra covariance relations
$[D,\phi(v,n)]=-\phi(v,n), \quad [X(m),\phi(v,n)]=\phi(Xv,n+m)$.
Exactly as in [37], the field $\phi$ is uniquely determined by these
relations and its initial term $\phi(v,0)$ in ${\rm Hom}_G(V_j\otimes
V,V_i)$. Our main new result is that all vector primary fields arise
by compressing fermions and therefore satisfy the $L^2$ bound above.

\bf \noindent Braiding relations. \rm If $f$ and $g$ have disjoint
supports in $S^1$, then $a(f)a(g)=-a(g)a(f)$ and 
$a(f)a(g)^*=-a(g)^*a(f)$. Similar but more complex ``braiding
relations'' hold for vector primary fields and their adjoints.
These may be summarised as follows. Let $a,b \in L^2(S^1,{\Bbb C}^N)$ 
be supported in intervals $I$ and $J$ in
$S^{1}\backslash\{1\}$ with $J$ anticlockwise after $I$. Define
$a_{gf} =\phi_{gf}^\square(e_{-\alpha} a)$ and
$b_{gf}=\phi^\square_{gf}(e_{-\alpha} b)$, with
$\alpha=(\Delta_g-\Delta_f-\Delta_\square)/2(N+\ell)$.  
Then $$b_{gf}a_{fh}=\sum \mu_{f_1} a_{gf_1}b_{f_1h}, \quad b_{gf}a_{g_1
f}^*=\sum \nu_h a_{hg}^* b_{hg_1},$$
with all coefficients non--zero. The proof of these relations is similar
to that of the Bargmann--Hall--Wightman theorem ([10],[19],[34]). 
To prove the first for example let  $F_k(z) =\sum
(\phi_{ik}(u,n)\phi_{kj}(v,-n)v_j,v_i)z^n$, a
power series convergent for $|z|<1$ with values in $W={\rm
Hom}_{SU(N)}(V_j\otimes U\otimes V,V_i)$. To prove the braiding
relation, it suffices to show that $F_k$ extends continuously to $S^1\backslash
\{1\}$ and $F_k(e^{i\theta}) =\sum c_{kh}
e^{i\mu_{kh}\theta} F_h(e^{-i\theta})$ there. 
Using Sugawara's formula for $D$, we show directly that the
$F_k$'s satisfy the Knizhnik--Zamoldchikov ODE [22]
$${dF\over dz} ={P F\over z} + {Q F\over 1-z},$$ 
where $P, Q\in {\rm
End}(W)$ (the original proof in [22], referred to in [37], is
different and less elementary).  In all cases we need, the matrix $P$
has distinct eigenvalues, none of 
which differ by positive integers, and $Q$ is a non--zero multiple of
a rank one idempotent in general position with respect to $P$. For two
vector primary fields this ODE reduces to the classical hypergeometric ODE
and the required relation on $S^{1}\backslash\{1\}$ follows from
Gauss' formula for transporting solutions at $0$ to $\infty$. In
general the ODE can be related to the
generalised hypergeometric ODE for which the corresponding transport
relations were first obtained by Thomae [36] in 1867 in terms of
products of gamma functions. Such a link exists because there is a
basis of $W$ for which $P$ and $P-Q$ are both in rational canonical
form. In this basis, the ODE is just the matrix form of the
generalised hypergeometric ODE. 
\vskip .05in
\noindent \bf Transport formulas. \rm The operator $a_{\square
0}^*a_{\square 0}$ on $H_0$ commutes with $L_{I^c}SU(N)$, so lies in
$\pi_0(L_ISU(N))^{\prime\prime}$. Therefore, by local equivalence, we
have the right to 
consider its image under $\pi_f$. We obtain the fundamental 
``transport formula'':
$\pi_f(a_{\square 0}^*a_{\square 0})=\sum \lambda_g
a_{gf}^*a_{gf}$,
with $\lambda_g>0$. Thus for $ T\in {\rm Hom}_{L_I G}(H_0,H_f)$, we have
$$T a_{\square 0}^*a_{\square 0})=\sum \lambda_g
a_{gf}^*a_{gf}T.$$ We will prove the transport formula
by induction using the braiding 
relations; the original proof in [41] used the transport relations between
$0$ and $1$ of the basic ODE above. 
\vskip .05in
\noindent \bf Definition of Connes fusion. \rm We develop the ideas of fusion
directly at the level of loop groups without appeal to the general theory of
bimodules over von Neumann algebras ([8], [32], [40], [41]) Let $X$,
$Y$ be positive energy representations of $LSU(N)$ at level $\ell$.
Let ${\cal X}={\rm Hom}_{L_{I^c}SU(N)} (H_0,X)$ and
${\cal Y}={\rm Hom}_{L_ISU(N)}(H_0,Y)$. These spaces of bounded
intertwiners or fields replace vectors or states in $X$ and $Y$. Thus
$x\in {\cal X}$ ``creates'' the state $x\Omega$ from the vacuum
$\Omega$. The fusion 
$X\boxtimes Y$ is defined to be the completion of ${\cal X}\otimes
{\cal Y}$ with respect to the pre--inner product 
$\langle x_1\otimes y_1,x_2\otimes y_2\rangle
=(x_2^*x_1y_2^*y_1\Omega,\Omega)$,
a four--point function. $X\boxtimes Y$ admits a natural action of
$L_ISU(N)\times L_{I^c}SU(N)$. The map 
${\cal X}\otimes {\cal Y}\rightarrow X\boxtimes Y$ extends to
continuous maps $X\otimes {\cal Y}\rightarrow X\boxtimes Y$ and ${\cal
X}\otimes Y\rightarrow X\boxtimes Y$. This implies that if ${\cal
X}_0\subset {\cal X}$ and $\overline{{\cal X}_0\Omega}=X$, then 
${\cal X}_0\otimes {\cal Y}$ has dense image in $X\boxtimes Y$. Fusion
is associative and $X\boxtimes H_0=X=H_0\boxtimes X$. 
\vskip .05in
\noindent \bf Explicit computation of fusion. \rm We use the transport
formula to prove the fusion
rule $H_\square \boxtimes H_f=\bigoplus H_g$ where 
$g$ ranges over permissible signatures obtained by adding a box to
$f$. The transport formula is still true if $a_{gf}$ is
replaced by linear combinations $x_{gf}$ of intertwiners $\pi_g(h)a_{gf}$ with
$h\in L_IG$. But then for $y\in {\rm Hom}_{L_ISU(N)}(H_0,H_f)$ we have
$(x^*xy^*y\Omega,\Omega)=(y^*\pi_f(x^*x)y\Omega,\Omega)=
\sum \lambda_g \|x_{gf}y\Omega\|^2$. Thus $U(x\otimes y)=\bigoplus
\lambda_g^{1/2} x_{gf}y\Omega$ gives the required unitary intertwiner from
$H_\square\boxtimes H_f$ onto $\bigoplus H_g$.
Similar reasoning can be used to prove that $H_f\boxtimes H_{[k]}\le
\bigoplus_{g>_k f} H_g$, where $g$ runs over all permissible
signatures that can be obtained by adding $k$ boxes to $f$ with no two
boxes in the same row. This time a transport formula must be proved
with $a_{\square 0}$ replaced by a path
$a_{k,k-1} a_{k-1,k-2} \cdots
a_{\square 0}$ indexed by exterior powers. This
device of considering products of vector primary fields means that 
we can avoid the use of smeared primary fields corresponding to the exterior
powers $\lambda^k {\Bbb C}^N$ which need not be bounded [41]. 
\vskip .05in
\noindent \bf The fusion ring. \rm It follows immediately from the
fusion rule with $H_\square$ that the $H_f$'s are
closed under fusion. Moreover, if
$R$ denotes the operator corresponding to rotation by $180^{\rm o}$, then
the formula $B(x\otimes y)=R^*[RyR^*\otimes RxR^*]$ gives a unitary
intertwining $X\boxtimes Y$ and $Y\boxtimes X$; this is a less refined
form of the braiding operation that makes the level $\ell$
representations into a braided tensor category [42]. Thus the representation
ring ${\cal R}$ of formal sums $\sum m_f H_f$ becomes a commutative
ring. For each permitted signature $h$, let
$z_h\in SU(N)$ be the diagonal matrix with entries
$\exp(2\pi i (h_k + N-k -H) /(N+\ell))$ where $H=(\sum h_k
+N-k)/N$; these give a subset ${\cal T}$. Let ${\cal S}\subset {\Bbb
C}^{\cal T}$ be the image of $R(SU(N))$ under the map of restriction
of characters. Our main result asserts that the natural ${\Bbb
Z}$--module isomorphism ${\rm 
ch}:{\cal R}\rightarrow {\cal S}$ defined by $[H_f]\mapsto [V_f]$ is a
ring isomorphism. This completely determines the fusion rules. They
agree with the well--known ``Verlinde formulas'' ([38],[20]), in which 
the usual tensor product rules for $SU(N)$ are modified by an action
of the affine Weyl group.

\vskip .05in

\noindent \bf Discussion. \rm Many of the early versions of the
results in Chapter~II were worked out in 
discussions with Jones in 1989--1990 (see [18] and [40]). We were mainly
interested in the inclusion $\pi_i(L_IG)^{\prime\prime}
\subseteq \pi_i(L_{I^c}G)^\prime$ defined by the ``failure of Haag
duality''. Algebraic quantum field theory [14] provided a series of
predictions about these local loop group algebras which we interpreted
(in the language of [28])
and verified. In particular two of the main theorems, Haag--Araki
duality and loop group irreducibility, were originally obtained with Jones. 
In the case of geometric modular theory for fermions on $S^1$,
each of us came up with different proofs which appear
in simplified form here (see also [40]). The original proofs
of irreducibility have been superseded by the simpler and more widely
applicable method described above. One of our original proofs followed
from the stronger result that $L^A G$ is dense in $LG$ in the natural
topology on $U_P(H)$, so that $\pi(L^A G)$ is strong
operator dense in $\pi(LG)$ for any positive energy representation;
the analogous result fails for ${\rm Diff}\, S^1$ and its discrete
series representations. The geometric method of
descent from local fermion algebras to local loop group algebras and
its application to 
Haag duality and local equivalence were first suggested by me, but it
was Jones who pointed out that this approach tacitly assumed
Takesaki's result [35] (``Takesaki devissage'').

The first paper of this series [40] is an expanded version of
expository lectures given in the Borel seminar in Bern in 1994. Since
it was intended as an
introduction to the general theory, we included a complete treatment
of the whole theory of fusion, braiding and subfactors for the
important special case of $LSU(2)$. In the
second paper of the 
series [41] we made a detailed study of primary fields from several points
of view. (See Jones' S\'eminaire Bourbaki talk [46] for a detailed summary.) 
We constructed all primary fields as compressions of tensor
products of fermionic operators, thus establishing their analytic
properties. To do so, we had to complete and extend the Lie algebraic
approach of Tsuchiya and Kanie [37] and in particular prove the conjectured
four--point property of physicists. Fusion of positive energy
representations was computed
using the braiding properties of primary fields. The braiding
coefficients appeared
as transport coefficients between different singular points of the
basic ODE studied here; these coefficients were
derived using Karamata's Tauberian theorem and a unitary trick. Since
the smeared primary 
fields could be unbounded, their action had to be controlled by
Sobolev norms; and a detailed argument had to be supplied for
extending the braiding relations to arbitrary bounded intertwiners.

In this paper we give a more elementary approach to fusion using only
vector primary fields and their adjoints. It is not possible to
overemphasise the central r\^ole 
(prophesied by Connes) played by the fermionic model in our work, nor
the importance of considering the relationships between
different models (stressed by P.~Goddard). The boundedness of the
corresponding smeared fields is very significant. Not only does it
simplify the analysis, 
but more importantly it can be seen to lie at the heart of the crucial
irreducibility result (due to the duality between smeared
primary fields and loop group observables). This is an
example of Goddard's philosophy that ``vertex operators tell you
what to do.''
With the important exception of the Lie algebra
operators (indispensable for proving the KZ equation), we have tried
to keep exclusively to bounded operators. This is in line with Rudolf Haag's 
philosophy that quantum field theory can and should be understood in
terms of (algebras of) bounded operators [14]. Here, because of the
boundedness of vector primary fields, there is no choice.

In the fourth paper of this series [42] we explain how the positive energy
representations at a fixed level become a braided tensor category. We
have already seen a simplified version of the braiding operation when
proving that Connes fusion is commutative. The key to understanding
this braiding structure lies in the ''monodromy'' action of the 
braid group on products of vector primary fields.
The important feature of braiding allows us to make contact with the
subfactors of the hyperfinite type II$_1$ factor defined by Jones
and Wenzl ([17],[18],[43]) using special traces on the infinite braid
group. In particular this explains the coincidence between the
monodromy representation of the braid group in [37] and the Hecke
algebra representations of Jones and Wenzl. Further developments
include understanding the ``modularity'' of the category, the property
which allows 3--manifold invariants to be defined. This involves 
studying the elliptic KZ equations as well as finding and versifying
precise versions of the axioms for a CFT; the ideas behind our
computation of fusion seem to give a general method for understanding
unitarity and positivity properties of quite general CFTs. 
Specifically one can enlarge the monodromy action on paths of
primary fields to include Jones projections; our transport formulas
can then be interpreted as defining an inner product on such paths
making the monodromy action unitary. In addition the
analytic properties of primary fields implied by our construction
(such as the fact that $q^{L_0} \phi(z)$ is a
Hilbert--Schmidt operator for $|q|<1$)  should allow primary fields to
be interpreted as morphisms corresponding to 3--holed spheres or
trinions in Graeme Segal's language. This should yield
a precise analytic version of the Segal's ``modular functor'',
using the ``operator 
formalism'' for trinion decompositions of Riemann surfaces.

The braiding properties of vector primary fields can also be developed
through a more systematic use of the {\it
conformal inclusion} $SU(N)\times SU(\ell)\subset SU(N\ell)$. The level one
representations and vector primary fields of $SU(N\ell)$, when restricted
to $SU(N)\times SU(\ell)$ and decomposed into tensor products, yield
all representations and vector primary fields of $SU(N)$ at levels
$\ell$ and $SU(\ell)$ at level
$N$. The level one representations of $LSU(N\ell)$ arise by
restricting the fermionic representation of $LU(N\ell)$ to
$LSU(N)\times LU(1)$ (here $U(1)$ is the centre of $U(N\ell)$). Our
fermionic construction of primary fields for $LSU(N)$ in this and the
previous paper have been a simplification of the more sophisticated 
picture provided by the above conformal inclusion, first considered
from this point of view by Tsuchiya \& Nakanishi [26]. Here we have
ignored the 
r\^ole of the group $SU(\ell)$. If it is brought into play, it is
possible to give a less elementary but more conceptual
non--computational proof that all 
the braiding coefficients are non--zero, based on the Abelian braiding
of fermions or vector primary fields at level one. This approach,
which will be taken up in detail when we consider subfactors defined
by conformal inclusions, has the advantage firstly that it makes the
non--vanishing of the coefficients manifest and secondly that it does
not require the explicit solutions of the KZ ODE and their monodromy
properties that we have used here and in the second paper. It therefore
extends to 
other groups where less information about the KZ ODE is available at
present.

\vskip .1in
\noindent \bf Acknowledgements. \rm I would like to thank various
people for their help: Jurg Fr\"ohlich for suggesting that positive
energy representations might be studied 
from the point of view of algebraic quantum field theory;
Alain Connes for signalling the importance of fermions and
tensor products of bimodules in such a study; Klaus Fredenhagen for
outlining the literature in algebraic quantum field theory; Sorin Popa
for his results on type~III subfactors; Peter Goddard for guiding me
through the literature in conformal field theory; Richard Borcherds
for logistic aid; Terence Loke for
simplifying the proof of associativity of Connes fusion; Hans
Wenzl for explaining his results on Hecke algebra fusion and
subfactors; and Vaughan Jones for many helpful discussions
and encouragement, particularly during the early collaborative stage of
this research.

\vfill\eject
\noindent \bf CHAPTER I.~POSITIVE ENERGY REPRESENTATIONS OF LSU(N).
\vskip .2in
\noindent \bf 2. Irreducible representations of SU(N). \rm We give
a brief account of the representation theory of $SU(N)$ from a point
of view relevant to this paper. This account closely parallels 
our development of the classification and fusion of positive energy
representations of $LSU(N)$, so provides a simple prototype. Let
$V={\Bbb C}^N$ be the vector representation. We shall
consider irreducible representations of $SU(N)$ appearing in tensor
powers $V^{\otimes m}$. Let $R(SU(N))$ denote the 
representation ring of $SU(N)$, the ring of formal integer
combinations of such irreducible representations. Let $\g$ be the Lie
algebra of $SU(N)$, the traceless skew--adjoint matrices. Thus $\g$
acts on $V^{\otimes m}$, hence each irreducible representation $W$, and 
${\rm End}_G(W) = {\rm End}_{\g}(W)$. This representation of $\g$
extends linearly to a *--representation of its
complexification $\g_{\Bbb C}$, the traceless matrices. $\g_{\Bbb C}$
is spanned by the elementary matrices $E_{ij}$ ($i\ne j$) and
traceless diagonal matrices.  Let $T$
denote the subgroup of diagonal matrices $z=(z_1,z_2,\dots,z_N)$ in
$SU(N)$.Given an irreducible representation
$SU(N)\rightarrow U(W)$, we can write $W=\bigoplus_{g\in
{\Bbb Z}^N} W_g$ with
$\pi(z) v=z^g v$ for $v\in W_g$, $z\in T$. We call $g$ a
{\it weight} and $W_g$ a {\it weight space}; $g$ is only determined up
to addition of a vector $(a,a,\dots,a)$ for $a\in {\Bbb Z}$. The 
monomial matrices in $SU(N)$ permute the weight spaces
by permuting the entries of $g=(g_1,\dots,g_N)$, so there is always a
weight with $g_1\ge g_2\ge \dots \ge g_N$. Such a weight is called a
{\it signature}. If the weights are ordered lexicographically, the
raising operators 
$\pi(E_{ij})$ ($i<j$) carry weight spaces into weight
spaces of higher weight; their adjoints $\pi(E_{ij})$  ($i>j$) are
called lowering operators and decrease weight.

  Clearly every irreducible projective representation $W$ 
contains a highest weight vector $v$. Now $W$ is irreducible for
$\g_{\Bbb C}$ and every monomial $A$ of operators in $\g_{\Bbb C}$ is a sum
of products $LDR$ where $L$ is a product of lowering operators, $D$ is
a product of diagonal operators and $R$ is a product of raising
operators. Since $LDRv$ is proportional to $v$ or has lower weight,
$v$ is unique up to a multiple. On the 
other hand $(A_1v,A_2 v)$ is 
uniquely determined by the weight of $v$ and the $A_i$'s, since 
$A_2^*A_1$ can be written as a sum of operators $LDR$ and
$(LDRv,v)=(DRv,L^*v)$ with $L^*$ a raising operator.  Thus if $W^\prime$ is
another irreducible representation with the same highest weight and
corresponding vector $v^\prime$, 
$Av\mapsto Av^\prime$ is a unitary $W\rightarrow W^\prime$
intertwining $\g$ and hence $G=\exp(\g)$. Thus irreducible
representations are classified by their signatures. Every signature
occurs: if $f_1\ge f_2\ge \cdots \ge f_N\ge 0$, the 
vector $e_f=e_1^{\otimes (f_1-f_2)} \otimes (e_1\wedge e_2)^{\otimes
(f_2-f_3)} \otimes \cdots \otimes 
(e_1\wedge e_2\wedge \cdots \wedge e_N)^{\otimes f_N}$ is the unique
highest weight vector in $\lambda^1 V^{\otimes(f_1-f_2)}\otimes
\lambda^2 V^{\otimes(f_2-f_3)} \otimes \cdots \otimes \lambda^N
V^{\otimes f_N} \subseteq V^{\otimes(\sum f_i)}$. By uniqueness, $e_f$
generates an irreducible submodule.

 A signature $f$ with $f_N\ge 0$ is represented by a Young diagram with 
at most $N$ rows and $f_i$ boxes in the $i$th row. Thus $V$
corresponds to the diagram $\square$ and $\lambda^k V$ to the diagram
$[k]$ with $k$ rows, with one box in each row. We write $g> f$ if $g$
can be obtained by adding one box to $f$. More generally we write
$g>_k f$ if $g$ can be obtained by adding $k$ boxes to $f$ with no two
in the same row.
\vskip .1in
\noindent \bf Lemma. \it ${\rm Hom}_G(V_f\otimes V_{[k]}, V_g)$
is at most one--dimensional and only non--zero if $g>_k f$. When
$k=1$, it is non--zero iff $g>f$. Hence $V_f\otimes V_\square
=\bigoplus_{g>f} 
V_g$ and $V_f\otimes \lambda^k V\le \bigoplus _{g>_k f} V_g$.
\vskip .05in
\noindent \bf Proof. \rm Let $v_f$ and $v_g$ be highest weight vectors
in $V_f$ and $V_g$. If $T\in {\rm Hom}_G(V_f \otimes V_{[k]},
V_g)$ with $T(v_f\otimes v)=0$ for all $v\in \lambda^k V$, then applying
lowering operators we see that $T=0$. If $T\ne 0$, we take
$w=e_{i_1}\wedge \dots 
\wedge e_{i_k}$ of highest weight such that $T(v_f\otimes w)\ne 0$.
Applying raising operators, we see that $T(v_f\otimes w)$ is highest
weight in $V_g$ so is proportional to $v_g$. So the weight of
$v_f\otimes w$ is a signature and $g>_k f$. If $S$ is another
non--zero intertwiner, we may choose $\alpha$ such that $R=S-\alpha T$
satisfies $R(v_f\otimes w)=0$. If $R\ne 0$, we may choose $w^\prime$
of highest weight such that $R(v_f\otimes w)\ne 0$. But this gives a
contradiction, since $R(v_f\otimes w)$ would be annihilated by all
raising operators and have weight lower than $v_g$. So ${\rm
Hom}_G(V_f\otimes V_{[k]}, V_g)$ is at most one--dimensional.

If $g$ is obtained by adding a box to the $i$th row of
$f$, then the map
$$T:\lambda^1 V^{\otimes(f_1-f_2)}\otimes
\lambda^2 V^{\otimes(f_2-f_3)} \otimes \cdots \otimes \lambda^N
V^{\otimes f_N} \, \bigotimes \, V\rightarrow \lambda^1
V^{\otimes(g_1-g_2)}\otimes 
\lambda^2 V^{\otimes(g_2-g_3)} \otimes \cdots \otimes \lambda^N
V^{\otimes g_N}$$
given by exterior multiplication by $V$ on the $(f_1-f_i)$th copy of
$\Lambda V$ commutes with $G$ and satisfies $T(e_f\otimes e_i)=e_g$.
Thus if $P$ and $Q$ denote the projections onto the submodules
generated by $e_f$ and $e_g$ respectively, $QT(P\otimes I)$ gives a
non--zero intertwiner $V_f\otimes V\rightarrow V_g$.

\vskip .1in
\rm For $z_i\in {\Bbb C}$ and a signature $f$, we define the symmetric
function $X_f(z)={{\rm det} (z_j^{f_i+n-i})/{\rm det}(z_j^{n-i})}$.
The denominator here is a Vandermonde determinant given by
$\prod_{i<j} (z_i-z_j)$. If $X_k(z)= \sum_{i_1<\cdots i_k}$ $ z_{i_1}
\dots z_{i_k}$, then it is elementary to show that $X_f X_k =
\sum_{g>_k f} X_g$ for $k=1,\dots,N$. In particular 
$X_k(z)$ coincides with $X_{[k]}(z)$; and it follows, by induction
on $f_1-f_N$ and the number of boxes in the $f_1$th column, that
each $X_f(z)$ is an integral polynomial in the $X_k(z)$'s.

\vskip .1in

\bf \noindent Theorem. \it (1) $V_f\otimes V_{[k]}=\bigoplus_{g>_k f}
V_g$. 

\noindent (2) $R(SU(N))$ is generated by the exterior powers and the
map ${\rm ch}:[V_f]\rightarrow X_f$ gives a ring 
isomorphism 
between $R(SU(N))$ and ${\cal S}_N$, the ring of symmetric integral
polynomials in $z$, where $\prod z_i =1$. 

\noindent (3) (Weyl's character formula [44]) $\chi_f(z) \equiv {\rm
Tr}(\pi_f(z)) =X_f(z)$ for all $f$.

\vskip .1in
\noindent \bf Proof. \rm (1) We know that $V_f\otimes \lambda^k V \le
\bigoplus_{g>_k f} V_g$. We prove
by induction on $|f|=\sum f_j$ that 
$V_f\otimes V_k = \bigoplus_{g_1>_k f} V_{g_1}$. It suffices to show
that if this holds for $f$ then it holds for all 
$g$ with $g>f$. Now, comparing the coefficients of $X_h$ in
$(X_fX_k) X_\square= (X_f X_\square)X_k$, we see that $|\{g_1:h>_k g_1
> f\}|=|\{g_2:h> g_2>_k f\}|$. Tensoring by $V_\square$,
we deduce that $\oplus_{g>f} V_g\otimes V_{[k]} =\oplus_{g_1>_k f}
\oplus_{h>g_1} V_h=\oplus_{g>f} \oplus_{h>_k g} V_h$.
Since $ V_g\otimes V_{[k]} \le \bigoplus_{h>_k g} V_h$, we must have
equality for all $g$, completing the induction.

\noindent (2) Let ${\rm ch}$ be the ${\Bbb Z}$--linear isomorphism ${\rm
ch}:R(SU(N))\rightarrow {\cal S}_N$ extending ${\rm ch}(V_f)=X_f$.
Then by (1), ${\rm ch}(V_{[k]}V_f) =X_k X_f$. This implies that ${\rm
ch}$ restricts to a ring homomorphism on the subring of $R(SU(N))$
generated by the exterior powers. On the other hand the $X_k$'s
generate ${\cal S}_N$, so the image of this subring is the whole of
${\cal S}_N$. Since ${\rm ch}$ is injective, the ring generated by the
exterior powers must be the whole of $R(SU(N))$ and ${\rm ch}$ is thus
a ring homomorphism, as required.

\noindent (3) The maps $[V_f]\rightarrow \chi_f(z)$ and $[V_f]\mapsto
X_f(z)$ define ring homomorphisms $R(SU(N))\rightarrow {\Bbb C}$.
These coincide on the exterior powers and therefore everywhere.

\vskip .1in
\noindent \bf 3. Fermions and quantisation. \rm Given a complex
Hilbert space $H$, the complex Clifford algebra ${\rm Cliff}(H)$ is
the unital *--algebra generated
by a complex linear map $f\mapsto a(f)$ ($f\in H$) satisfying the
anticommutation relations
$a(f)a(g)+a(g)a(f)=0$ and $a(f)a(g)^*+a(g)^*a(f)=(f,g)$ (complex
Clifford algebra relations). The Clifford algebra has a natural action
$\pi$ on $\Lambda H$ (fermionic Fock space) given by $\pi(a(f))\omega
=f\wedge \omega$, called the complex wave representation. The complex wave
representation is irreducible. For $\Omega$ is the unique vector such
that $a(f)^*\Omega=0$ for all $f$ (this condition is equivalent to
orthogonality to $\sum_{k\ge 1} \lambda^k H$) and $\Omega$ is cyclic
for the $a(f)$'s. Thus if $T\in {\rm End}(\Lambda H)$ commutes with
all $a(f)^*$'s, $T\Omega =\lambda \Omega$ for $\lambda\in
{\Bbb C}$; and if $T$ also commutes with all $a(f)$'s, $T=\lambda I$.  

\vskip .05in
To produce other irreducible representations of ${\rm
Cliff}(H)$, we introduce the operators $c(f)=a(f)+a(f)^*$.  Thus $c$
satisfies 
$c(f)=c(f)^*$, $f\mapsto c(f)$ is real--linear and
$c(f)c(g)+c(g)c(f)=2{\rm Re}(f,g)$ (real Clifford algebra relations).
The equations $c(f)=a(f)+a(f)^*$ and 
and $a(f)=(c(f)-ic(if))/2$ give a correspondence between complex and
real Clifford algebra relations. Since
$c$ relies only on the underlying real Hilbert space $H_{\Bbb R}$, complex
structures on $H_{\Bbb R}$ commuting with $i$ give new irreducible
representations of ${\rm Cliff}(H)$. These complex structures correspond to
projections $P$ in $H$: multiplication by $i$ is given by $i$ on $PH$
and 
$-i$ on $(PH)^\perp$. Unravelling this definition, we find that the
projection $P$ defines an irreducible representation $\pi_P$ of ${\rm
Cliff}(H)$ on fermionic Fock space ${\cal F}_P=
\Lambda\,PH\widehat{\otimes} \Lambda\, (P^\perp 
H)^*$ given by $\pi_P(a(f))= a(Pf) + a((P^\perp f)^*)^*$.
(Equivalently $\pi_p(a(f))=(c(f) -ic(i(2P-I)f))/2$ on $\Lambda H$.)
\vskip .1in
\noindent \bf Theorem (Segal's equivalence criterion [3]). \it Two
irreducible 
representations $\pi_P$ and $\pi_Q$ are unitarily 
equivalent if $P-Q$ is a Hilbert--Schmidt operator.
\vskip .05in
\noindent \bf Remark. \rm The converse also holds ([3],[40]), but will
not be needed.
\vskip .05in
\noindent \bf Proof. \rm If $PH$ (or $P^\perp H$) is
finite--dimensional, then so is $QH$ (or $Q^\perp H$) and the
representations are easily seen to be equivalent to the irreducible
representation on $\Lambda H$ (or $\Lambda H^*$). So we may assume
that ${\rm dim}\, PH ={\rm dim}\, P^\perp H=\infty$. 

The operator $T=(P-Q)^2$ is compact, so by the spectral theorem
$H=\bigoplus_{\lambda\ge 0} 
H_\lambda$ where $T\xi=\lambda \xi$ for $\xi\in H_\lambda$. 
Moreover ${\rm dim}\,H_\lambda<\infty $ for $\lambda>0$ while $P=Q$ on
$H_0$.  
Now $T$ commutes with $P$ and $Q$, so that each $H_\lambda$ is
invariant under $P$ and $Q$. Thus $H$ can be written as a direct sum
of finite--dimensional irreducible submodules $V_i$ for $P$ and $Q$,
with $(P-Q)^2$ a scalar $\lambda$ on each. Since the images of $P$
and $Q$ (and $I$) should generate ${\rm End}(V_i)$, the identity
$(P-Q)^2=\lambda I$ forces ${\rm dim}\, {\rm End}(V_i)\le 4$. Hence
${\rm dim}\, V_i=1$ or $2$. 

Pick an orthonormal basis $(e_i)_{i\ge -a}$ of $P^\perp H$ with each $e_i$
lying in 
some $V_j$. We may assume that $Q^\perp e_{-1}=Q^\perp e_{-2}= \cdots
=Q^\perp e_{-a}=0$ and
that $Q^\perp e_i\ne 0$ for $i\ge 0$. Complete $(e_i)$ to an orthonormal
basis $(e_i)_{i\in {\Bbb Z}}$ by adding remaining vectors from the
$V_j$'s. We can 
also choose an orthonormal basis $(f_j)_{j\ge -b}$ of $Q^\perp H$ with
$f_i$ lying in the same $V_j$ as $e_i$ if $i\ge 0$; we shall even pick
$f_i$ so that $(e_i,f_i)>0$ in this case. A simple computation shows
that if $(P-Q)^2=\lambda_i I$ on $V_j$, then $(e_i,f_i)=\sqrt{1-\lambda_i}$
(so $\lambda_i=0$ when ${\rm dim}\, V_j=1$). Note that, using these
bases, we get $\|P-Q\|_2^2={\rm
Tr}\, (P-Q)^2= a+b +2\sum \lambda_i$, so that $\sum\lambda_i<\infty$.

The ``Dirac sea'' model ${\cal H}$ for $\Lambda H_P$ is
the Hilbert space with orthonormal basis given by all symbols
$e_{i_1}\wedge e_{i_2} \wedge e_{i_3} \wedge \cdots $ where
$i_1<i_2<i_3 <\cdots $ and $i_{k+1} =i_k +1$ for $k$ sufficiently
large. If $A(e_i)$ denotes exterior multiplication by $e_i$, then
$A(e_i)A(e_j)+A(e_j)A(e_i)=0$ and
$A(e_i)A(e_j)^*+A(e_j)^*A(e_i)=\delta_{ij}I$. By linearity and
continuity, these extend to operators $A(f)$ ($f\in H$) satisfying the
complex Clifford algebra relations so give a representation $\pi$ of
${\rm Cliff}(H)$. Let $\xi=e_{-a} \wedge
e_{-a+1} \wedge \cdots $. Then the $A(f)$ and $A(f)^*$'s act
cyclically on $\xi$ and 
$(A(f_1) \dots A(f_m) \xi, A(g_1) \cdots
A(g_n)\xi)= \delta_{mn} \det (Pf_i,g_j)$.
On the other hand $(\pi_P(a(f_1))\dots
\pi_P(a(f_m))\Omega_P, \pi_P(a(g_1))\dots
\pi_P(a(g_n))\Omega_P)=\delta_{mn} \det (Pf_i,g_j)$, 
where $\Omega_P$ is the vacuum vector in $\Lambda H_P$. Thus
$(\pi(a)\xi,\xi)= (\pi_P(a)\Omega_P,\Omega_P)$ for $a\in {\rm
Cliff}(H)$. Replacing $a$ by $a^*a$ and recalling that $\xi$ and
$\Omega_P$ are cyclic, we see that 	
$U(\pi_P(a)\Omega_P) =\pi(a)\xi$ defines a unitary from $\Lambda H_P$
onto ${\cal H}$ such that $\pi(a) =U\pi_P(a) U^*$. The same
``Gelfand--Naimark--Segal'' argument shows that unitary equivalence of
$\pi_P$ and $\pi_Q$ will follow as soon as we
find $\eta\in {\cal H}$ such that
$(\pi(a)\eta,\eta)=(\pi_Q(a)\Omega_Q,\Omega_Q)$. (Note that $\eta$ is
automatically cyclic, since ${\cal H}\cong \Lambda H_P$ is
irreducible.)

Let $\eta_N= f_{-b} \wedge \cdots \wedge f_{-1} \wedge f_0\wedge
\cdots \wedge 
f_N\wedge e_{N+1} \wedge e_{N+2} \wedge \cdots $. Clearly if $a$ lies
in the *--algebra generated by the $a(e_i)$'s, then $(\pi(a)\eta_N,\eta_N)
=(\pi_Q(a)\Omega_Q,\Omega_Q)$ for $N$ sufficiently large. Thus it will
suffice to show that $\eta_N$ has a limit $\eta$, i.e.~$(\eta_N)$ is a
Cauchy sequence. Since $\|\eta_N\|=1$, this follows if ${\rm Re}\,
(\eta_M,\eta_N)\rightarrow 1$ as $M\le N\rightarrow \infty$. But
$(\eta_M,\eta_N)=\prod_{i=M+1}^N (e_i,f_i)=\prod_{i=M+1}^N
\sqrt{1-\lambda_i}$ and, as $\sum \lambda_i <\infty$, this tends to
$1$ if $M,N\rightarrow \infty$, as required. 

\vskip .05in
\bf \noindent Corollary of proof. \it If $\pi_P$ and $\pi_Q$ are
unitarily equivalent and $\Omega_Q$ is the image of the vacuum vector
in ${\cal F}_Q$ in ${\cal F}_Q$, then $|(\Omega_P,\Omega_Q)|^2= \prod
(1-\mu_i)$ where $\mu_i$ are the eigenvalues of $(P-Q)^2$. 
\vskip .05in
\noindent \bf Proof. \rm We have
$|(\Omega_P,\Omega_Q)|=|(\xi,\eta)|= \lim |(\xi,\eta_N)|=\prod
(1-\mu_i)^{1/2}$.
\vskip .1in
Any $u\in U(H)$ gives rise to a Bogoliubov
automorphism of ${\rm 
Cliff}(H)$ via $a(f)\mapsto a(uf)$. This automorphism is said to be
implemented in $\pi_P$ (or on ${\cal F}_P$) if
$\pi_P(a(uf))=U\pi_P(a(f))U^*$ for some unitary $U\in U({\cal F}_P)$ 
unique up to a phase. Since $\pi_P(a(uf))=\pi_Q(a(f))$ with $Q=u^*Pu$,
we immediately deduce:
\vskip .05in
\noindent \bf Corollary (Segal's quantisation criterion [3],
[28],[40]). \it $u$ is implemented in ${\cal F}_P$ if $[u,P]$ is a
Hilbert--Schmidt operator.   
\vskip .05in
\rm We define the restricted unitary group $U_P(H)=\{u\in
U(H):\,\hbox{$[u,P]$ Hilbert--Schmidt}\}$, a topological group under
the strong operator topology combined with the metric
$d(u,v)=\|[u-v,P]\|_2$. By the corollary, there is 
a homomorphism $\pi$ of $U_P(H)$  into $PU({\cal F}_P)$,
called the {\it basic} projective 
representation. 
\vskip .1in
\noindent \bf Lemma. \it The basic representation is continuous. 
\vskip .05in
\noindent \bf Proof. \rm It is
enough to show continuity at the identity. Thus if $u_n\mapright{s} I$ and
$\|[u_n,P]\|_2\rightarrow 0$, we must find a lift $U_n\in U({\cal
F}_P)$ of $\pi(u_n)$ such that $U_n\mapright{s} I$.
Now $\|[u_n, P]\|_2=\|P-Q_n\|_2$ where
$Q_n=u_n^*Pu_n$. So ${\rm Tr}(P-Q_n)^2\rightarrow 0$. On the other
hand $|(\Omega_P,\Omega_{Q_n})|^2=\prod (1-\mu_i)$ where $\mu_i$ are the
(non--zero) eigenvalues of $(P-Q_n)^2$. Since ${\rm Tr}(P-Q_n)^2=\sum
\mu_i$ and $\prod (1-\mu_i) \ge
\exp (-2\sum \mu_i)$ for $\sum \mu_i$ small, it follows that
$|(\Omega_P,\Omega_{Q_n})|\rightarrow 1$ as $n\rightarrow \infty$. If
$u_n$ is implemented by $U_n$ in ${\cal F}_P$, then $U_n\Omega_P$ and
$\Omega_{Q_n}$ are equal up to a phase. So
$|(U_n\Omega_P,\Omega_P)|\rightarrow 1$. Adjusting $U_n$ by a phase,
we may assume $(U_n \Omega_P,\Omega_P)>0$ eventually so that
$U_n\Omega_P \rightarrow \Omega_P$. Now, taking operator norms, $\|U_n
\pi(a(f)) U_n^* - \pi(a(f))\| =\|\pi(a(u_n f-f))\|\le \|u_n f-f\|$. 
It follows that $\|U_n aU_n^* -a\|\rightarrow 0$ for any $a\in
\pi_P({\rm Cliff}\, H)$. Thus $U_n a\Omega_P= (U_n a U_n^*)
(U_n\Omega_P)\rightarrow a\Omega_P$ as $n\rightarrow \infty$. Since
vectors $a\Omega_P$ are dense in ${\cal F}_P$, we get $U_n\mapright{s}
I$, as required. 
\vskip .1in
Note that if $[u,P]=0$, so that $u$ commutes with $P$,
then $u$ is canonically implemented in Fock space ${\cal F}_P$ and we
may refer to the {\it canonical quantisation} of $u$. If on the contrary
$uPu^*=I-P$, then $u$ is canonically implemented by a
conjugate--linear isometry in Fock space, also called the canonical
quantisation of $u$. Thus the canonical quantisations correspond to
unitaries that are complex--linear or conjugate--linear for the
complex structure defined by $P$.

\vskip .1in
\noindent \bf 4. The fundamental representation. \rm Let $G=SU(N)$ (or
$U(N)$) 
and define the loop group $LG=C^\infty(S^1,G)$, the smooth maps of the
circle into $G$. Let $H=
L^2(S^1)\otimes V$ ($V={\Bbb C}^N$) and let $P$ be the projection onto
the Hardy space 
$H^2(S^1)\otimes V$ of functions with vanishing negative Fourier
coefficients (or equivalently boundary values of functions holomorphic
in the unit disc). Now $LG$ acts unitarily by multiplication on $H$.
In fact if $f\in C^\infty(S^1,{\rm End}\, V)$ and $m(f)$ is the
corresponding multiplication operator, then it is easy to check,
using the Fourier coefficients of $f$, that
$\|[P,m(f)]\|_2\le \|f^\prime\|_2$. In particular $LG$ satisfies Segal's
quantisation criterion  
for $P$ and we therefore get a projective representation 
of $LU(N)$ on ${\cal F}_P$ ([28],[40]), continuous for the $C^\infty$
topology on $LU(N)\subset C^\infty(S^1,{\rm End}\, V)$. The rotation
group ${\rm 
Rot}\,S^1$ acts by 
automorphisms on $LG$ by $(r_\alpha f)(\theta)=f(\theta+\alpha)$.
The same formula defines a unitary action on $L^2(S^1)\otimes V$ which
leaves $H^2(S^1)\otimes V$ invariant. Therefore this action of ${\rm
Rot}\, S^1$ is canonically quantised and we thus get a projective
representation of $LG\rtimes {\rm Rot}\, S^1$ on ${\cal F}_P$ which 
restricts to an ordinary representation on ${\rm Rot}\, S^1$.

Let 
$$SU_\pm(1,1)=\{\pmatrix{\alpha &\beta\cr \overline{\beta} &
\overline{\alpha}\cr}:|\alpha|^2-|\beta|^2 =\pm 1\}$$
and let $SU_+(1,1)=SU(1,1)$ and $SU_-(1,1)$ denote the elements with
determinant $+1$ or $-1$. Thus $SU_-(1,1)$ is a coset of $SU_+(1,1)$
with representative $F=\pmatrix{0 & -1\cr -1& 0\cr}$, for example. The
matrices $g\in SU_\pm(1,1)$ act by 
fractional linear transformations on $S^1$, 
$g(z)=(\alpha z +\beta)/(\overline{\beta} z+\overline{\alpha})$. This
leads to a unitary 
action on $L^2(S^1,V)$ via $(V_g \cdot f)(z) =(\alpha
-\overline{\beta} z)^{-1}f(g^{-1}(z))$. Since $(\alpha
-\overline{\beta} z)^{-1}$ is 
holomorphic for $|z|<1$ and $|\alpha|>|\beta|$, it follows that $V_g$
commutes with the Hardy 
space projection $P$ for $g\in SU_+(1,1)$. The matrix $F$ acts on
$L^2(S^1,V)$ via $(F \cdot f)(z) =z^{-1}f(z^{-1})$ and clearly
satisfies $FPF=I-P$. It follows that $F$ is canonically implemented in
fermionic Fock space ${\cal F}_V$ by a conjugate--linear isometry
fixing the vacuum vector. Since $SU_-(1,1)= SU_+(1,1)F$, the same
holds for each $g\in SU_-(1,1)$. Thus we get
an orthogonal representation of $SU_\pm(1,1)$ for the underlying real
inner product on ${\cal F}_V$ with $SU_+(1,1)$ preserving the complex
structure and $SU_-(1,1)$ reversing it.  The same is true in ${\cal
F}_V^{\otimes \ell}$.  

\vskip .1in
Let $U_z$ denote the canonically quantised action of the gauge group
$U(1)$ on ${\cal F}_V$, corresponding to multiplication by $z$ on $H$.
The ${\Bbb Z}_2$--grading on ${\cal F}_V$ is given by the operator
$U=U_{-I}$. 

\vskip .1in
\noindent \bf Lemma. \it $\pi(g)U_z \pi(g)^*=U_z$ for all $g\in
LSU(N)$ and $z\in U(1)$. 
\vskip .05in
\noindent \bf Proof. \rm The group $SU(N)$ is simply connected, so the
group $LSU(N)$ is connected (any path can be smoothly contracted to
a constant path and $SU(N)$ is connected). The map $U(H)
\times U(H)\rightarrow U(H)$, $(u,v)\mapsto uvu^*v^*$ is continuous
and descends to $PU(H)\times PU(H)$. So $(u,v)\mapsto uvu^*v^*$
defines a continuous map $PU(H)\times PU(H)\rightarrow U(H)$. Since
$g$ and $z$ commute on the prequantised space $H$, $\pi(g)$ and $U_z$
commute in $PU(H)$. Hence $\pi(g)U_z\pi(g)^*U_z^*=\lambda(g,z)$ where
$\lambda(g,z)\in {\Bbb T}$ depends continuous on $g$ and $z$. Writing
this equation as $\pi(g)U_z\pi(g)^*=\lambda(g,z)U_z$, we see that
$\lambda(g,\cdot)$ defines a character $\lambda_g$ of $U(1)$. Clearly
$\lambda_g\lambda _h =\lambda_{gh}$, so we get a continuous homomorphism
of $LSU(N)$ into $\widehat{U(1)}$, the group of characters of $U(1)$.
Since $\widehat{U(1)}={\Bbb Z}$ and $LSU(N)$ is connected,
$\lambda_g=1$ for all $g$. So $\lambda(g,z)=1$ for all $g,z$ as
required.
\vskip .1in
\noindent \bf Corollary. \it Each operator $\pi(g)$ with $g\in LSU(N)$
is even (it commutes with $U=U_{-1}$).
\vskip .1in

\noindent \bf 5. The central extension ${\cal L}G$. \rm 
We introduce
the central extension of $LG$ 
$$1\rightarrow {\Bbb T} \rightarrow {\cal L}G \rightarrow LG
\rightarrow 1$$
obtained by pulling back the central extension
$1\rightarrow {\Bbb T} \rightarrow U({\cal F}_V) \rightarrow PU({\cal
F}_V) \rightarrow 1$
under the map $\pi:LG\rightarrow PU({\cal F}_V)$. In other words it is the
closed subgroup of $LG\times U({\cal F}_V)$ given by $\{(g,u):
\pi(g)=[u]\}$: it contains ${\Bbb T}=1\times {\Bbb T}$ as a central
subgroup and has quotient $LG$. By definition ${\cal L}G$ has a unique
unitary representation $\pi$ on ${\cal F}_V$ given by
$\pi(g,u)=u$. This extension is compatible with the action of
$SU_{\pm}(1,1)$ and ${\rm Rot}\, S^1$.
\vskip .1in

\noindent \bf Lemma. \it If $\pi(\gamma)$ denotes the canonical
quantisation of
$\gamma\in SU_\pm(1,1)$ on fermionic Fock space ${\cal F}_V$ and ${\cal
L}G=\{(g,u):\pi(g)=[u]\}$, then the operators $(\gamma,\pi(\gamma))$
normalise $\pi({\cal L}G)$ acting on the centre ${\Bbb T}$ as the
identity if $\gamma\in 
SU_+(1,1)$ and as complex conjugation if $\gamma\in SU_-(1,1)$. 
\vskip .05in
\noindent \bf Proof. \rm This follows because
$\pi(\gamma)\pi(g)\pi(\gamma)^{-1}$ has the same image as
$\pi(g\cdot\gamma^{-1})$ in $PU({\cal F}_V)$. 

\vskip .1in

\noindent \bf 6. Positive energy representations. \rm We may consider
the decomposition of ${\cal F}_P=\Lambda(PH)\otimes
\Lambda(P^\perp H)^*$ into weight spaces of ${\rm Rot}\, S^1={\Bbb T}$,
writing ${\cal F}_P=\bigoplus_{n\ge 0} 
{\cal F}_P(n)$, where $z\in {\Bbb T}$ acts on ${\cal F}_P(n)$ as
multiplication by $z^n$. Since ${\rm Rot}\, S^1$ acts with finite
multiplicity and only non--negative weight spaces on $PH$ and $(P^\perp
H)^*$, it is easy to see that ${\cal F}_P(n)$ is finite--dimensional
for $n\ge 0$ and ${\cal F}_P(n)=(0)$ for $n<0$. Moreover ${\cal
F}_P(0)=\Lambda(V)$. We define a representation of ${\Bbb T}$ on $H$
to have {\it positive energy} if in the decomposition $H=\bigoplus
H(n)$ we have $H(n)=0$ for $n<0$ and $H(n)$ finite--dimensional for
$n\ge 0$. (Usually we will also insist on the normalisation $H(0)\ne
(0)$, which can always be achieved through tensoring by a character of
${\Bbb T}$.) Thus ${\rm Rot}\, S^1$ acts on ${\cal F}_V$ with positive
energy. 
\vskip .1in
\bf \noindent Proposition. \it Suppose that $\Gamma$ is a subgroup of
$U(H)$
and that ${\Bbb T}$ acts on $H$ with positive energy normalising $\Gamma$.
Let $U_t$ be the action (with $t\in [0,2\pi]$). 
\vskip .05in
\item{(a)} If $H$ is irreducible as an $\Gamma\rtimes {\Bbb
T}$--module, then 
it is irreducible as a $\Gamma$--module.
\item{(b)} If $H_1$ and $H_2$ are irreducible $\Gamma\rtimes {\Bbb
T}$--modules 
which are isomorphic as $\Gamma$--modules, then one is obtained from the
other by tensoring with a character of ${\Bbb T}$.
\item{(c)} If $H$ is the cyclic $\Gamma$--module generated by a lowest
energy vector, it contains an irreducible $\Gamma\rtimes {\Bbb
T}$--module generated by some lowest energy vector. 
\item{(d)} Any positive energy representation is a direct sum of
irreducibles.
\vskip .1in
\noindent \bf Proof. \rm (a) Let $M=\Gamma^\prime$, the commutant of
$\Gamma$, so 
that $M=\{T:\hbox{$Tg=gT$ for all $g\in \Gamma$}\}$. By Schur's lemma,
$M\cap \langle U_t\rangle^\prime ={\Bbb C}I$ since $\Gamma$ and 
${\Bbb T}$ act irreducibly. Note that $U_t$ normalises $M$, since it
normalises $\Gamma$. Let $v$ be a  
lowest energy vector in $H$. $v$ is cyclic for $\Gamma$ and ${\Bbb
T}$ and hence $\Gamma$, so $av\ne 0$
for $a\ne 0$ in $M$. If $M\ne {\Bbb C}$, there is a non--scalar
self--adjoint element $T\in M$. Define $T_n\in B(H)$ by  $(T_n
\xi,\eta) =(2\pi)^{-1} \int_0^{2\pi} e^{-int}
(U_t T_n U_t^*\xi,\eta)\, dt$. Then $T_n\in M$, $U_t TU_t^*=e^{int}
T_n$, $T_n^*=T_{-n}$ and $Tv=\oplus T_nv$. By assumption $T_0$ must be
a scalar. Since $T\notin{\Bbb C}I$,
$Tv$ cannot be a multiple of $v$ and therefore $T_n\ne 0$ for some
$n\ne 0$. Since $T_n^*=T_{-n}$, we may assume $n<0$. But then $T_nv\ne
0$ gives a vector of lower energy than $v$. So $M={\Bbb C}$ and
$\Gamma$ acts irreducibly. 

\noindent (b) Let $T:H_1\rightarrow H_2$ be a unitary
intertwiner for $\Gamma$. Then $V_t^*TU_t$ is also a unitary intertwiner, 
so must be of the form $\lambda(t) T$ for $\lambda(t)\in {\Bbb T}$ by
Schur's lemma. Since $TU_tT^*=\lambda(t) V_t$, $\lambda(t)$ must be a
character of ${\Bbb T}$. 

\noindent (c) Let $V$ be the subspace of lowest energy. Let $K$ be any
$\Gamma\rtimes {\Bbb T}$--invariant subspace of $H$ with corresponding
projection $p\in \Gamma^\prime$. Since $H=\overline{{\rm lin}}(\Gamma
V)$, $K=pH = 
\overline {\rm lin}(\Gamma pV)$. But $pV\subseteq V$, since $p$
commutes with ${\Bbb T}$. Choosing $pV$ in $V$ of smallest dimension,
we see that $K=\overline {\rm lin}(\Gamma pV)$ must be irreducible as a
$\Gamma\rtimes {\Bbb T}$--module and hence as a $\Gamma$--module. Thus
$H$ contains an irreducible submodule $K$ generated by any non--zero
$pv$ with $v\in V$.

\noindent (d) Take the cyclic module generated by
a vector of lowest energy. This contains an irreducible submodule
generated by another vector of lowest energy $H_1$ say. Now repeat
this process for $H_1^\perp$, to get $H_2$, $H_3$, etc. The positive
energy assumption shows that $H=\bigoplus H_i$.
\vskip .1in
\noindent \bf Corollary. \it If $\pi:LG\rtimes {\rm Rot}\,
S^1\rightarrow PU(H)$ is a projective representation which restricts
to an ordinary positive energy representation of ${\rm Rot}\, S^1$,
then $H$ decomposes as a direct sum $\bigoplus H_i\otimes K_i$ where
the $H_i$'s are representations of $LG\rtimes {\rm Rot}\, S^1$
irreducible on $LG$ with $H_i(0)\ne (0)$ and the multiplicity spaces
are positive energy representations of ${\rm Rot}\, S^1$.
\vskip .1in

\rm We apply this result to the positive energy representation ${\cal
F}_P^{\otimes \ell}$ of $LG\rtimes {\rm Rot}\, S^1$. The irreducible
summands of ${\cal F}_P^{\otimes \ell}$ are called the level $\ell$
irreducible representations of $LG$. By definition any positive energy
representation extends to $LG\rtimes {\rm Rot}\, S^1$. More generally
the vacuum representation at level $\ell$ extends (canonically) to
$LG\rtimes SU_\pm(1,1)$. 
In fact, since $SU_{\pm}(1,1)$ fixes the vacuum vector and this
generates the vacuum 
representation at level $\ell$ as an $LG$--module, it follows that
the vacuum representation at level $\ell$ admits a compatible orthogonal
representation of $SU_{\pm}(1,1)$, unitary on $SU_+(1,1)$ and
antiunitary on $SU_-(1,1)$. We also need the less obvious fact
that $SU(1,1)$ is implemented by a  
projective unitary representation in any level $\ell$ representation;
this follows from a global form of the Goddard--Kent--Olive construction
[11].
\vskip .1in
\noindent \bf Lemma (coset construction). \it Let $H=\bigoplus
H_i\otimes K_i$ and let 
$M=\bigoplus B(H_i)\otimes I$. Let $\pi:{\cal G}\rightarrow PU(H)$ be
a projective unitary representation of the connected topological group
${\cal G}$ such 
that $\pi(g)M\pi(g)^*=M$ for all $g\in {\cal G}$. Then there exist
projective unitary representations $\pi_i$ and $\sigma_i$ of ${\cal
G}$ on $H_i$ and $K_i$ such that $\pi(g)=\oplus \pi_i(g)\otimes
\sigma_i(g)$. 
\vskip .05in
\noindent \bf Proof. \rm ${\cal G}$ acts by automorphisms on $M$
through conjugation. It therefore preserves the centre and hence
the minimal central projections. Since ${\cal G}$ is connected and the action
strong operator continuous, it must fix the central projections. Thus
it fixes each block $H_i\otimes K_i$. It also normalises $B(H_i)$.
If $W_i$ denotes the restriction of $\pi(g)$ to $H_i\otimes K_i$, then
${\rm Ad}\,W_i$ restricts an automorphism $\alpha_i$ of $B(H_i)$. But,
if $K$ is a Hilbert space, any
automorphism $\alpha$ of $B(K)$ is inner: indeed if $\xi$ is a fixed
unit vector in $K$ and $P_\xi$ is the rank one projection onto ${\Bbb C}\xi$,
then $\alpha(P_\xi)=P_\eta$ for some unit vector $\eta$ and $U(T\xi)
=\alpha(T)\eta$ ($T\in B(K)$) defines a unitary with $\alpha={\rm
Ad}\, U$. Hence $\alpha_i={\rm Ad}\,U_i$ for $U_i\in U(H_i)$. But then
$(U_i^*\otimes I)W_i$ commutes with $B(H_i)\otimes I$ and hence lies
in $I\otimes  B(K_i)$. Hence $(U_i^*\otimes I)W_i =I\otimes V_i$, so
that $W_i=U_i\otimes V_i$. Thus we get the required homomorphism
${\cal G}\rightarrow  \prod PU(H_i)\times PU(V_i)$, which is clearly
continuous. 

\vskip .1in
\noindent \bf Corollary. \it There is a (unique) projective
representation $\pi_i$ of 
$SU(1,1)$ on $H_i$ satisfying
$\pi_i(\gamma)\pi_i(g)\pi_i(\gamma)^*=\pi_i(g\cdot\gamma^{-1})$ 
for $g\in {\cal L}G$ and $\gamma\in SU(1,1)$.
\vskip .05in
\noindent \bf Proof. \rm If $H={\cal F}_V^{\otimes\ell}$, we may write
$H=\bigoplus H_i\otimes K_i$ where the $H_i$'s are the distinct
level~$\ell$ irreducible representations of ${\cal L}G$ and the
$K_i$'s are multiplicity spaces. Then $\pi({\cal
L}G)^{\prime\prime}=\bigoplus B(H_i)\otimes I$ and the unitary
representation of $SU(1,1)$ normalises this algebra.
By the coset construction, each
$\gamma\in SU(1,1)$ has a decomposition $\pi(\gamma)=\bigoplus
\pi_i(\gamma) \otimes \sigma_i(\gamma)$, where
$\tau_i(\gamma)=\pi_i(\gamma)\otimes \sigma_i(\gamma)$ is an ordinary
representation of $SU(1,1)$ on $H_i\otimes K_i$. But 
$\pi_i(g\cdot\gamma^{-1})\otimes I=\tau_i(\gamma)(\pi_i(g)\otimes I)
\tau_i(\gamma)^* =\pi_i(\gamma)\pi_i(g)\pi_i(\gamma)^*$.
Hence 
$\pi_i(\gamma)\pi_i(g)\pi_i(\gamma)^*=\pi_i(g\cdot\gamma^{-1})$.
So, as before, the representation of $SU(1,1)$, now projective, is
compatible with the central extension ${\cal L}G$. 
\vskip .1in

\noindent \bf 7. Infinitesimal action of $L^0\g$ on finite energy
vectors. 
\rm If $\g={\rm
Lie}(G)$, then ${\rm Lie}(LG)=L\g=C^\infty(S^1,\g)$. Let $L^0\g$ be
the algebraic Lie algebra consisting of trigonometric polynomials with
values in $\g$.  Its
complexification is spanned by the functions
$X_n(\theta)=e^{-in\theta}X$ with  
$X\in \g$. ${\rm Rot}\, S^1$ and its Lie algebra act on $L^0\g$. The
Lie algebra of ${\rm Rot}\, S^1$ is generated by $id$ where
$[d,f](\theta) =-if^\prime(\theta)$ for $f\in L^0\g$. Thus $d$ may be
identified with the operator $-id/d\theta$. We obtain the Lie algebra
relations  $[X_n,Y_m]=[X,Y]_{n+m}$ and $[d,X_n]=-nX_n$. For $v\in
V$, let $v(n)=a(v_n)$ where $v_n\in L^2(S^1,V)$ is 
given by $v_n(\theta)=e^{-in\theta}v$. In particular, if $(e_i)$ is an
orthonormal basis of $V$, then we have fermions $e_i(n)$ for all $n$.
If $\Omega$ denotes the vacuum vector in ${\cal F}_V$, then it is easy
to see from its description as an exterior algebra that an orthonormal
basis of ${\cal F}_V$ is given by
$$e_{i_1}(n_1)e_{i_2}(n_2) \cdots e_{i_p}(n_p) e_{j_1}(m_1)^*
e_{j_2}(m_2)^* \dots e_{j_q}(m_q)^*\Omega$$
where $n_i\le 0$ and $m_j> 0$. Moreover $e_i(n)\Omega=0$ for $n\ge
0$ and $e_i(n)^*\Omega=0$ for $n<0$. Since ${\rm Rot}\, S^1$ commutes
with the Hardy space projection on $L^2(S^1,V)$, it is canonically
quantised. Let $R_\theta$ be the corresponding representation on
${\cal F}_V$. Then $R_\theta=e^{iD\theta}$ where $D$ is self--adjoint.
If $r_\theta$ is the action of ${\rm Rot}\, S^1$ on $L^2(S^1,V)$ given
by $(r_\theta f)(z)=f(e^{i\theta}z)$, then $r_\theta=e^{id}$ where
$d=-i{d\over d\theta}$ (we always regard functions on $S^1$ as
functions either of $z\in {\Bbb T}$ or of $\theta\in [0,2\pi]$). Now
$R_\theta a(f) R_\theta^*=a(r_\theta f)$. Hence $R_\theta v(m)
R_\theta^* =e^{-im\theta}v(m)$, so that $R_\theta$ acts on the basis vector
$e_{i_1}(n_1)e_{i_2}(n_2) \cdots e_{i_p}(n_p) e_{j_1}(m_1)^*
e_{j_2}(m_2)^* \cdots e_{j_q}(m_q)^*\Omega$ as multiplication by
$e^{iM\theta}$ where $M=\sum m_j -\sum n_i$. Since $R_\theta
=e^{iD\theta}$, it follows that $D$ acts on this basis vector as
multiplication by $M$, i.e.~this vector has energy $M=\sum m_j -\sum
n_i$. In particular $D\Omega=0$ and we can check that
$[D,v(n)]=-nv(n)$. Thus if $f$ is a trigonometric power series with
values in $V$, we have $[D,a(f)]=a(df)$. Note that if $T$ is a linear
operator on ${\cal F}^0_V$ commuting with the $e_i(a)$'s and
$e_i(a)^*$'s, then $T=\lambda I$ for $\lambda\in {\Bbb C}$: for, as in
section 3, $\Omega$ is the unique vector such that $e_i(n)^*\Omega=0$
($n\ge 0$), $e_i(n)\Omega=0$ ($n>0$) and $\Omega$ is cyclic.

\vskip .1in
\noindent \bf Theorem. \it Let $E_{ij}(n) =\sum_{m>0} e_i(n-m) e_j(-m)^* -
\sum_{m\ge 0} e_j(m)^*e_i(m+n)$,
and define $X(n)=\sum a_{ij} E_{ij}(n)$ for $X=\sum a_{ij}E_{ij}\in
{\rm Lie}\, U(V)\subset {\rm End}(V)$.
Then, as operators on $H^0$, we have
\vskip .05in
\item{(a)} $[X(m),a(f)]=a(X_m\cdot f)$ if $f$ is a trigonometric
polynomial with values in $V$; equivalently $[X(n),v(m)]=(Xv)(n+m)$.

\item{(b)} $[D,X(m)]=-mX(m)$.

\item{(c)} $[X(n),Y(m)]=[X,Y](n+m) + n (X,Y) \delta_{n+m,0} I$ where
$(X,Y)=-{\rm Tr}(XY)={\rm Tr}(XY^*)$ for $X,Y\in {\rm Lie}\, U(V)$. 

\vskip .05in
\noindent \bf Proof. \rm (a) Observe that
$[e_i(a)^*e_j(b),e_k(c)]=-\delta_{ac}\delta_{ik}e_j(b)$ and
$[e_j(b)e_i(a)^*,e_k(c)]=\delta_{ac} \delta_{ik}
e_j(b)$. Moreover if $i\ne j$, then $e_i(a)$ anticommutes with both
$e_j(b)$ and $e_j(b)^*$. Using these identities, it is easy to check
that $E_{ij}(n)$ satisfies the commutation relations (a)
with respect to the $e_i(n)$'s. Note that $X(n)\Omega=0$ for $n\ge 0$
since $e_i(n)\Omega=0$ for $n\ge 0$,  $e_i(n)^*\Omega=0$ for $n<0$ and
(formally) $X(n)^*=-X(-n)$ for $X\in {\rm Lie}\, U(V)$. 

\noindent (b) Since $[D,e_i(m)]=-m e_i(m)$
and $[D,e_i(m)^*]=me_i(m)^*$, it follows that $[D,X(m)]=-m X(m)$.

\noindent (c) From (a) we find that $T=[X(m),Y(n)]-[X,Y](m+n)$ commutes
with all 
$e_i(a)$'s and hence also all $e_i(a)^*$'s by the adjointness property.
Hence $[X(m),Y(n)]=[X,Y](m+n)
+\lambda(X,Y)(m,n)I$,
where $\lambda(m,n)$ is a scalar, bilinear in $X$
and $Y$. Now from (b), $[X(m),Y(n)]-[X,Y](m+n)$ lowers the energy by
$-m-n$, so that $\lambda(X,Y)(m,n)=0$ unless $m+n=0$. To compute the
value of 
$\lambda$ when $m=-n$, we note that we may assume that
$m\ge 0$, since $\lambda(X,Y)(m,n)^*=\lambda(Y,X)(-n,-m)$ by
the adjoint relations. Taking vacuum expectations, we get
$$\lambda(X,Y)(-m,m)= ([X(-m),Y(m)]\Omega,\Omega)=
(X(-m)\Omega,Y(-m)\Omega)= -m {\rm
Tr}(XY)=m(X,Y).$$ 
In fact if $X=\sum a_{ij}E_{ij}$ and $Y=\sum b_{ij} E_{ij}$, we have
$$(X(-m)\Omega,Y(-m)\Omega)=\sum_{ijpq} \sum_{r,s=0}^{m-1} (a_{ij}
e_j(r)^*e_i(r-m)\Omega, b_{pq}
e_q(s)^*e_p(s-m)\Omega) =m \sum a_{ij} \overline{b_{ij}}=m (X,Y),$$
since the terms $e_i(a)^*e_j(b)\Omega$ with
$a\ge 0$ and $b<0$ are orthonormal.

\vskip .2in
\bf \noindent 8. The exponentiation theorem. \rm 
 We wish to show that
the Lie algebra action just defined on ${\cal F}_V$ exponentiates to
give the fundamental representation of $LSU(N)\rtimes {\rm Rot}\,
S^1$. We have already discussed the action of ${\rm Rot}\, S^1$, which
is canonically quantised. So we now must show that if $x$ is an element of
$L^0\g$ and $X$ is the corresponding operator constructed above, then
$\pi \exp x$ and $\exp X$ have the same image in $PU({\cal F})$. To
see that this completely determines $\pi$ on $LG$, we need the
following result on products of exponentials.
\vskip .1in
\noindent \bf Exponential lemma. \it Every element of $LG$ is a
product of exponentials in $L\g=C^{\infty}(S^1,\g)$. Products of
exponentials in $L^0\g$ are dense in $LG$.
\vskip .05in
\noindent \bf Proof. \rm If $g\in LG\subset C(S^1,M_N({\Bbb C}))$
satisfies $\|g-I\|_\infty <1$, then $\log g= \log(I-(I-g))$ lies in
$C^\infty(S^1, \g)=L\g$. Thus $\exp L\g$ contains an open
neighbourhood of $I$ in $LG$. Since $LG$ is connected, $\exp L\g$ must
generate $LG$, as required. 
\vskip .1in
The bilinear formulas for the Lie algebra operators $X$ immediately imply
Sobolev type estimates for the infinitesimal action of $L^0\g$ on
finite energy vectors. We define the Sobolev 
norms by  $\|\xi\|_s = \|(I+D)^s\xi\|$
for $s\in {\Bbb R}$, usually a half--integer. Recall that
if $A$ is a skew--adjoint operator, the smooth vectors for $A$ are the
subspace $C^\infty(A)=\bigcap {\cal D}(A^n)$ and for any $\xi\in
C^\infty(A)$ we have $e^{At}\xi=\sum_{i=0}^n {t^k\over
k!} A^k \xi   + O(t^{n+1})$.
\vskip .1in

\bf \noindent Exponentiation Theorem. \it Let $H= {\cal F}_V$ be the
level one fermionic representation of $LSU(V)$ and let $H^0$ be the
subspace of finite energy vectors. 
\vskip .02in
\item{(1)} For $x\in L^0\g$, there is a constant $K$
depending on $s$ and $x$ such that
$\|X\cdot\xi\|_s \le K \|\xi\|_{s+1}$
for $\xi\in H^0$, $X=\pi(x)$. 
\item{(2)} For each $x\in L^0\g$, the corresponding operator $X$ is
essentially skew--adjoint on $H^0$ and leaves $H^0$ invariant.
\item{(3)} Each vector in $H^0$ is a $C^\infty$ vector for any
$x\in L^0\g$.
\item{(4)} For $x\in L^0\g$, the unitary $\exp(X)$ agrees up to a
scalar with $\pi(\exp(x))$.
\vskip .05in

\noindent \bf Proof. \rm (1) It clearly suffices to prove the
estimates in the lemma for $X=E_{ij}(n)$ and $\xi$ of fixed
energy, say $D \xi =\mu \xi$. Then
$E_{ij}(n)\xi=\sum_{m>0} e_i(n-m) e_j(-m)^*\xi -
\sum_{m\ge 0} e_j(m)^*e_i(m+n)\xi$. 
So $\|E_{ij}(n)\xi\|\le 2(|n|+\mu)\|\xi\|$, since at most $2(|n|+\mu)$ of
the terms in the sums can be non--zero and each has norm bounded by
$\|\xi\|$. Hence for $s\ge 0$, 
$$\eqalign{\|E_{ij}(n)\xi\|_s
&\le (1+|n|+\mu)^s\|E_{ij}(n)\xi\|\le 2
(1+|n|+\mu)^s(|n|+\mu)\cr
&\le 2(1+|n|)^{s+1} (1+\mu)^{s+1}\|\xi\|\le
2(1+|n|)^{s+1}\|\xi\|_{s+1}.\cr}$$
\vskip .05in
\noindent (2) Clearly any $X\in L^0\g$ acts on $H^0$.
We need the Glimm--Jaffe--Nelson commutator theorem (see [10], [29] or
[40]): if $D$ be is the energy operator on $H^0$ and 
$X:H^0\rightarrow 
H^0$ is formally skew--adjoint with $X(D+I)^{-1}$, $(D+I)^{-1}X$ and
$(D+I)^{-1/2} [X,D](D+I)^{-1/2}$ bounded, then the closure of $X$ is
skew--adjoint. The Sobolev estimates show that these conditions hold
for $D$ and $X$, since $[D,X]$ is actually in $L^0\g$. 
\vskip .05in
\noindent (3) Since
$XH^0 \subset H^0$ and the $C^\infty$ vectors for $X$
are just $\cap {\cal D}(X^n)$, it follows that the
vectors in $H^0$ are $C^\infty$ vectors for $X$. 
\vskip .05in
\noindent (4) We prove the commutation relation
$e^{tX} a(f) e^{-tX} = a(e^{tx}f)$ for $f\in L^2(S^1)\otimes V$. We
start by noting that 
$$a(Xf)\xi= X a(f)\xi - a(f)X\xi$$
for $f$ a trigonometric polynomial with values in $V$, $X\in L^0\g $
and $\xi\in 
H^0$. We fix $X$ and $f$ and denote by $C^\infty(X)$ the space
of $C^\infty$ vectors for $X$, i.e.~$\cap {\cal D}(X^n)$. Now say $\xi\in
{\cal D}(X)$ and $f\in L^2(S^1,V)$. Take
$\xi_n\in H^0$, such that $\xi_n\rightarrow \xi$ and
$X\xi_n\rightarrow X\xi$, and $f_n$ trigonometric polynomials with
values in $V$ such that $f_n\rightarrow f$. Then $a(f_n)\xi_n\rightarrow
a(f)\xi$ and 
$Xa(f_n)\xi_n= a(Xf_n)\xi_n + a(f_n)X\xi_n\rightarrow a(Xf)\xi + a(f)
X\xi$. Since $X$ is closed, we deduce that $a(f)\xi$ lies in ${\cal D}(X)$
and $a(Xf)\xi=Xa(f)\xi -a(f)X\xi$. Successive applications of this
identity imply that 
$a(f)\xi$ lies in ${\cal D}(X^n)$ for all $n$ if $\xi$ lies in
$C^\infty(X)$, so that $a(f)C^\infty(X)\subset C^\infty(X)$.

Now take $\xi,\eta\in C^\infty(X)$ and consider
$F(t)=(e^{-Xt} a(e^{xt}f) e^{Xt}\xi,\eta)=(a(e^{xt}f)
e^{Xt}\xi,e^{Xt}\eta)$.
Since $\xi,\eta$ lie in $C^\infty(X)$, we have $e^{X(t+s)}\xi =
e^{Xt}\xi + s 
Xe^{Xt}\xi + O(s^2)$ and $e^{X(t+s)}\eta=e^{Xt}\eta + sXe^{Xt}\eta +
O(s^2)$. For any $f$, 
we have $e^{x(t+s)}f= e^{xt}f + sxe^{xt}f + O(s^2)$ in
$L^2(S^1)\otimes V$. Since 
$\|a(g)\|=\|g\|$, it follows that $a(e^{x(t+s)}f) = a(e^{xt}f)
+sa(xe^{xt}f) + O(t^2)$ 
in the operator norm. Hence we get
$$\eqalign{F(t+s)&
=(a(e^{xt}f)e^{Xt}\xi,e^{Xt}\eta) \cr
&\quad +s
[(a(e^{xt}f)Xe^{Xt}\xi,e^{Xt} \eta) + (a(xe^{xt}f)e^{Xt}\xi,e^{Xt}\eta) +
(a(e^{xt}f)e^{Xt}\xi,Xe^{Xt}\eta)] + O(s^2)\cr
&=(a(e^{xt}f) e^{Xt}\xi,e^{Xt}\eta) +O(s^2).\cr}$$
since $[X,a(g)]=a(xg)$. Thus $F(t)$ is differentiable with
$F^\prime(t)\equiv 0$. Hence $F(t)$ is constant and therefore equal to
$F(0)$. This proves that $e^{-tX}a(e^{tx}f)e^{tX}\xi = a(f)\xi$ for
$\xi\in H^0\subset C^\infty(X)$. Hence $a(e^{tx}f)=e^{tX}a(f)e^{-tX}$,
as required. Thus $e^{tX}$ implements the Bogoliubov automorphism
corresponding to $e^{tx}$. 
\vskip .1in
\noindent \bf Corollary. \it Let $H$ be a level $\ell$
positive energy representation of $LSU(N)$ and let $H^0$ be the
subspace of finite energy vectors. 
\vskip .05in
\item{(1)} There is a
projective representation of $L^0\g\rtimes {\Bbb R}$ on $H^0$ such
that $[D,X(n)]=-nX(n)$, $D^*=D$, $X(n)^*=-X(-n)$ and 
$[X(m),Y(n)]=[X,Y](n+m) +m\ell\delta_{m+n,0}\,(X,Y)$.
\item{(2)} For each $x\in L^0\g$, the corresponding operator $X$ is
essentially skew--adjoint on $H^0$ and leaves $H^0$ invariant.
\item{(3)} For $x\in L^0\g$, the unitary $\exp(X)$ agrees up to a
scalar with the corresponding group element in $LG$.
\item{(4)} Each vector in $H^0$ is a $C^\infty$ vector for any
$X$. 

\vskip .05in
\noindent \bf Proof. \rm We observe that the embedding $LSU(N)\subset
LU(N\ell)$ gives all 
representations of $LSU(N)$ at level $\ell$. The continuity properties
of the action of the larger group and its Lie algebra are immediately
inherited by $LSU(N)$. Note that it is clear from the functoriality of
the fermionic construction that the restriction of the fermionic
representation of $LU(N\ell)$ to $LU(N)$ can be identified with ${\cal
F}^{\otimes \ell}$ where ${\cal F}$ is the (level $1$) fermionic
representation of $LU(N)$. The other properties follow immediately
from the following result, applied to irreducible summands $K$ of $H={\cal
F}^{\otimes \ell}$.
\vskip .1in
\noindent \bf Lemma. \it Let $X$ be a skew--adjoint operator on $H$
with core $H^0$ such that $X(H^0)\subseteq H^0$. Let $K$ be a closed
subspace such that $P(H^0)\subseteq H^0$, where $P$ is the projection
onto $K$. Let $K^0=K\cap H^0$. Then $X(K^0)\subseteq K^0$ iff
$\exp(Xt) K=K$ for all $t$. In this case $K^0$ is a core for $X|_K$. 
\vskip .05in
\noindent \bf Proof. \rm Suppose that $K$ is invariant under $\exp
(Xt)$. Then $\exp(Xt)\xi =\xi + tX\xi +\cdots$ for $\xi\in K^0$ and
hence $XK^0\subseteq K\cap H^0 =K^0$.
Conversely, if $X(K^0)\subseteq K^0$, take $\xi\in {\cal D}(X)$
and let $P$ be the orthogonal projection onto $K$. It will suffice to
show that $P\xi\in {\cal D}(X)$ and $XP\xi=PX\xi$, for then $X$
commutes with $P$ in the sense 
of the spectral theorem. Since $P(H^0)\subseteq H^0$, we have $H^0=
H^0\cap K \bigoplus H^0\cap K^\perp$. Since $X$ is skew--adjoint and
$X(K^0)\subseteq K^0$, it
follows that $X$ leaves $H^0\cap K^\perp$ invariant. Thus $PX=XP$ on
$H^0$. Take $\xi_n\in H^0$ such that $\xi_n\rightarrow \xi$ and
$X\xi_n\rightarrow X\xi$. Then $XP\xi_n=PX\xi_n\rightarrow PX\xi$ and
$P\xi_n\rightarrow \xi$. Since $X$ is closed, $XP\xi=PX\xi$ as
required. Finally since $P\xi_n\rightarrow P\xi$ and
$XP\xi_n\rightarrow XP\xi$, it follows that $K^0$ is a core for
$X|_K$.

\vskip .1in

\noindent \bf 9. Classification of positive energy representations of
level $\ell$. \rm 
\vskip .1in
\noindent \bf Proposition. \it Let $(\pi,H)$ be an irreducible
positive energy 
projective representation of $LG\rtimes {\rm Rot}\, S^1$ of level
$\ell$.
\vskip .05in
\item{(1)} The action of $L^0\g\rtimes {\Bbb R}$ on $H^0$ is algebraically
irreducible.

\item{(2)} $H(0)$ is irreducible as an $SU(N)$--module.

\item{(3)} If $H(0)=V_f$, then $f_1-f_N\le \ell$.

\item{(4)} (Existence) If $f_1-f_N\le \ell$, there is a an irreducible
positive energy representation of $LG$ of level $\ell$ of the above
form with $H(0)\cong V_f$ as $SU(N)$--modules.

\item{(5)} (Uniqueness) If $H$ and $H^\prime$ are irreducible positive
energy representations of level $\ell$ of the above form with
$H(0)\cong H^\prime(0)$ as $SU(N)$--modules, then $H$ and $H^\prime$ are
unitarily equivalent as projective representations of $LG\rtimes {\rm
Rot}\, S^1$. 
\vskip .05in

\noindent \bf Proof. \rm (1) Recall that $H$ is irreducible as
an $LG\rtimes {\Bbb T}$--module iff it is irreducible as an
$LG$--module by the proposition in section~6. Any subspace $K$ of
$H^0$ invariant under $L^0\g\rtimes {\Bbb R}$ 
is clearly invariant under ${\rm Rot}\, S^1$. It therefore coincides
with the space of finite energy vectors of its closure. By the lemma
in section~8, its
closure is invariant under all operators $\exp(X)$ for $x\in L^0\g$.
But $\exp(L^0\g)$ generates a dense subgroup of $LG$, so the closure
must be invariant under $LG$ and therefore coincide with the whole of
$H$ by irreducibility. Hence $K=H^0$ as required.

\noindent (2) Let $V$ be an irreducible $SU(N)$--submodule of $H(0)$.
From (1), the $L^0\g\rtimes {\Bbb R}$--module generated by $V$ is the
whole of $H^0$. Since $D$ fixes $V$, it follows that the $L^0\g$--module
generated by $V$ is the whole of $H^0$. The commutation rules 
show that any monomial in the $X(n)$'s can be written as a sum of
monomials of the form $P_- P_0 P_+$, where $P_-$ is a monomial in the
$X(n)$'s for $n<0$ (energy raising operators), $P_0$ is a monomial in
the $X(0)$'s (constant energy operators) and $P_+$ is 
a monomial in the $X(n)$'s with $n>0$ (energy lowering operators).
Hence $H^0$ is spanned by 
products $P_- v$ ($v\in V$). Since the lowest energy subspace of this
$L^0\g$--module is $V$, 
$H(0)=V$, so that $H(0)$ is irreducible as a $G$--module.

\noindent (3) Suppose that $H(0)\cong V_f$ and let $v\in H(0)$ be a
highest weight vector, so that
$(E_{ii}(0)-E_{jj}(0))v =(f_i-f_j)v$ and $E_{ij}(0)v=0$ if $i<j$. Let
$E=E_{N1}(1)$, 
$F=E_{1N}(-1)$ and $H=[E,F]=E_{NN}(0) - E_{11}(0) + \ell$. Thus
$H^*=H$, $E^*=F$, $[H,E]=2E$ and $[H,F]=-2F$. Moreover $Ev=0$ and
$Hv=\lambda v$ with $\lambda=f_N-f_1+\ell$. By induction on $k$, we have
$[E,F^{k+1}]=(k+1)F^k (H-kI)$ for $k\ge 0$. Hence
$(F^{k+1}v,F^{k+1} v)=(F^*F^{k+1}v, F^k v)=(EF^{k+1}v,
F^kv)=(k+1)(\lambda-k) (F^kv,F^kv)$.
For these norms to be non--negative for all $k\ge 0$, $\lambda$ has to
be non--negative, so that $f_1-f_N\le \ell$ as required.

\noindent (4) We have ${\cal F}_V^{\otimes\ell}(0)=(\Lambda
V)^{\otimes \ell}$. By the results of section~6, the $LG$--module
generated by any irreducible 
summand $V_f$ of ${\cal F}_V(0)$ gives an irreducible positive energy
representation $H$ with $H(0)\cong V_f$. So certainly any irreducible
summand in ${\Lambda V}^{\otimes \ell}$ appears as an $H(0)$. From
the tensor product rules with the $\lambda^kV$'s, these
representations are precisely those with $f_1-f_N\le \ell$.

\noindent (5) Any monomial $A$ in operators from $L^0\g$ is a sum of
monomials $RDL$ with $R$ a monomial in energy raising operators, $D$ a
monomial in constant energy operators and $L$ a monomial in energy
lowering operators. As in section 2, if $v,w\in H(0)$ the inner products 
$(A_1v,A_2w)$ are uniquely determined by $v,w$ and the monomials $A_i$:
for $A_2^*A_1$ is a sum of terms $RDL$ and $(RDLv,w)=(DLv,R^*w)$ with
$R^*$ an energy lowering operator. Hence, if $H^\prime$ is another
irreducible positive energy representation with $H^\prime(0)\cong
H(0)$ by a unitary isomorphism 
$v\mapsto v^\prime$, $U(Av)=Av^\prime$ defines a unitary map of $H^0$
onto $(H^\prime)^0$ intertwining $L^0\g$. This induces a unique unitary
isomorphism $H\rightarrow H^\prime$ which intertwines the one parameter
subgroups corresponding to the skew--adjoint elements in $L^0\g$,
since $H^0$ and ${H^\prime}^0$ are cores for the corresponding
skew--adjoint operators. But these subgroups generate a dense subgroup
of $LG$, so that $U$ must intertwine the actions of $LG$, i.e.
$\pi^\prime(g) =U\pi(g) U^*$ in $PU(H^\prime)$ for $g\in LG$. Thus $H$
and $H^\prime$ are isomorphic as projective representations of $LG$.
From section~6,  $H$ and $H^\prime$ are therefore
unitarily equivalent as projective representations of $LG\rtimes {\rm
Rot}\, S^1$. 

\vskip .1in

\noindent \bf Corollary. \it The irreducible positive energy
representations $H$ of $LG$ of level $\ell$ are uniquely determined by
their lowest energy subspace $H(0)$, an irreducible $G$--module. Only
finitely many irreducible representations of $G$ occur at level
$\ell$: their signatures must satisfy the quantisation condition
$f_1-f_N\le \ell$. The action of $L^0\g\rtimes {\Bbb R}$ on $H^0$ is
algebraically irreducible.

\vfill\eject
\noindent \bf CHAPTER II.~LOCAL LOOP GROUPS AND THEIR VON NEUMANN
ALGEBRAS. 
\vskip .2in
\noindent \bf 10. von Neumann algebras. \rm Let $H$ be a Hilbert
space. The commutant of ${\cal 
S}\subset B(H)$ is 
defined by ${\cal S}^\prime=\{T\in B(H): Tx=xT\,\hbox{for all $x\in
{\cal S}$}\}$. If ${\cal S}^*={\cal S}$, for example if ${\cal S}$ is a
*--algebra or a subgroup of $U(H)$, then ${\cal S}^\prime$ is a unital
*--algebra, closed in the weak or strong operator topology. Such an
algebra is called a von Neumann algebra. von Neumann's double
commutant theorem states that ${\cal S}^{\prime\prime}$ coincides with
the von Neumann algebra generated by ${\cal S}$, i.e.~the weak
operator closure of the unital *--algebra generated by ${\cal S}$.
Thus a *--subalgebra $M\subseteq B(H)$ is a von Neumann algebra iff
$M=M^{\prime\prime}$. By the spectral theorem, the spectral
projections (or more generally bounded Borel functions) of any self--adjoint
or unitary operator in $M$ must also lie in $M$. This implies in
particular that $M$ is generated both by its
projections and its unitaries. Note that, if
$M={\cal S}^\prime$, the projections in $M$ correspond to 
subrepresentations for ${\cal S}$, i.e.~subspaces invariant under
${\cal S}$. 

The centre of a von Neumann algebra $M$ is given by $Z(M)=M\cap
M^\prime$. A von Neumann algebra is 
said to be a {\it factor} iff $Z(M)={\Bbb C}I$. A unitary
representation of a group or a *--representation of a *--algebra is
said to be a factor representation if the commutant is a factor. If
$H$ is a representation with commutant $M$, then two 
subrepresentations $H_1$ and $H_2$ of $H$ are unitarily equivalent iff
the corresponding projections $P_1,P_2\in M$ are the initial and final
projections of a partial isometry $U\in M$, i.e.~$U^*U=P_1$ and
$UU^*=P_2$. $P_1$ and $P_2$ are then said to be equivalent in the
sense of Murray and von Neumann [25]. We shall only need the following
elementary result, which is an almost immediate consequence of the
definitions. 
\vskip .1in
\noindent \bf Proposition. \it If $(\pi,H)$ is a factor 
representation of a set ${\cal S}$ with ${\cal S}^*={\cal S}$ and
$(\pi_1,H_1)$ and $(\pi_2,H_2)$ are subrepresentations, then
\vskip .05in
\item{(1)} there is a unique *--isomorphism $\theta$ of
$\pi_1({\cal S})^{\prime\prime}$ onto $\pi_2({\cal S})^{\prime\prime}$
such that $\theta(\pi_1(x))=\pi_2(x)$ for $x\in {\cal S}$;

\item{(2)} the intertwiner space ${\cal X}={\rm Hom}_{\cal
S}(H_1,H_2)$ satisfies $\overline{{\cal X}H_1}=H_2$, so in particular
is non--zero;

\item{(3)} $\theta(a)T=Ta$ for all $a\in \pi_1({\cal
S})^{\prime\prime}$ and $T\in {\cal X}$;

\item{(4)} if ${\cal X}_0\subseteq {\cal X}$
with $\overline{{\cal X}_0H_1}=H_2$, 
then $\theta(a) $ is the unique $b\in
\pi_2({\cal S})^{\prime\prime}$ such that $bT=Ta$ for all $T\in {\cal
X}_0$. 

\vskip .05in
\noindent \bf Proof. \rm Let $M=\pi({\cal S})^{\prime\prime}$ and
$M_i=\pi_i({\cal S})^{\prime\prime}$. Then
$\overline{M^\prime H_i}$ is invariant under both $M$ and
$M^\prime$. Hence the corresponding projection lies in $M\cap
M^\prime={\Bbb C}$ (since $M$ is a factor). So $\overline{M^\prime
H_i}=H$.  Let $p_i$ be the projection onto $H_i$, so 
that $p_i\in M^\prime$. Clearly $M_i=Mp_i$.
Moreover, the map $\theta_i:M\rightarrow M_i$,
$a\mapsto ap_i$ must be a *--isomorphism: for $ap_i=0$ implies
$aM^\prime H_i=(0)$ and hence $a=0$. 
By definition $\theta_i(\pi(x))=\pi_i(x)$ for $x\in {\cal S}$. 
Now set $\theta=\theta_2\theta_1^{-1}$;
$\theta$ is unique because $M_1$ is generated by $\pi_1({\cal S})$. 
 
Since ${\cal X}={\rm Hom}_{\cal
S}(H_1,H_2)=p_2M^\prime p_1$, we have $T\theta_1(x) =\theta_2(x)T$ for
all $x\in M$. Hence $\theta(a) T=Ta$ for $a\in M_1$ and $T\in {\rm
Hom}_{\cal S}(H_1,H_2)$. Moreover $\overline{{\cal
X}H_1}=\overline{p_2M^\prime H_2}=\overline {p_2H}=H_2$.
Conversely suppose that ${\cal X}_0 \subset {\rm Hom}_{\cal S}(H_1, H_2)$ is a
subspace such that ${\cal X}_0H_1$ is dense in $H_2$ and $a\in M_1$,
$b\in B(H_2)$ satisfy $bT=Ta$ for all $T\in {\cal X}_0$. Let
$c=b-\theta(a)$. Then $c{\cal X}_0=(0)$ and hence $cH_2=(0)$, so that
$c=0$. Thus $b=\theta(a)$ as required.

\vskip .1in
\bf \noindent 11. Abstract modular theory. \rm Let $H$ be a complex
Hilbert space, and $K\subset H$ a 
closed real subspace with $K\cap iK=(0)$ and $K+iK$ dense in $H$. Let
$e$ and $f$ be the projections onto $K$ and $iK$ respectively and set
$r=(e+f)/2$, $t=(e-f)/2$. Then $K^\perp$, $iK^\perp$ and $iK$ satisfy 
the same conditions as $K$, where $\perp$ is taken with respect to the
real inner product ${\rm Re}\,(\xi,\eta)$. 
\vskip .05in
\noindent \bf Proposition~1. \it (1) \it $0\le r\le I$, $t$, $r$ are
self--adjoint, $t$ is conjugate--linear, $r$ is linear, and $t$, $I-r$,
$r$ have zero kernels. 

\noindent (2)  $t^2=r(I-r)$, $rt=t(I-r)$, $(I-r)t=tr$. 

\noindent  (3)  $et=t(I-f)$, $ft=t(I-e)$. 

\noindent (4) If $t$ has polar decomposition $t=|t|j=j|t|$, then
$j^2=I$, $ej=j(I-f)$ and $fj=j(I-e)$. 

\noindent (5) $jK=iK^\perp$ and $(j\xi,\eta)\in {\Bbb R}$ for
$\xi,\eta\in K$.

\noindent (6) Let $\delta^{it}=(I-r)^{it}r^{-it}$. Then
$j\delta^{it} =\delta^{it} j$ and $\delta^{it}K=K$.

\vskip .05in
\noindent \bf Proof. \rm (1), (2) and (3) are straightforward.  (4)
follows from (3), because $e$ and $f$ 
commute with $t^2=(e-f)^2/4$, hence with $|t|$, and $|t|$ has zero kernel.
(4) implies (5), since $jej=I-f$. Finally  
since $jrj=I-r$ and $j$ is conjugate--linear, $j$ commutes with
$\delta^{it}$.  So $\delta^{it}$ commutes with $j$, $r$,
$|t|=\sqrt{r(I-r)}$ and hence $t$. So $\delta^{it}$ commutes with $e$
and $f$.

\vskip .1in

\noindent \bf Proposition~2 (characterisation of modular operators). \it
(1) (Kubo--Martin--Schwinger condition) For each $\xi\in K$,
the function 
$f(t)=\delta^{it}\xi$ on ${\Bbb R}$ extends (uniquely) to a continuous
bounded function $f(z)$ on
$-1/2\le {\rm Im}\, z\le 0$, holomorphic in $-1/2 <{\rm Im}\, z<0$.
Furthermore $f(t-i/2) = jf(t)$ for $t\in{\Bbb R}$.

\noindent (2) (KMS uniqueness) Suppose that $u_t$ is a one--parameter
unitary 
group on $H$ and $j_1$ is a conjugate--linear involution such that
$u_t K=K$ and $j_1u_t=u_t j_1$. Suppose that there is a dense subspace
$K_1$ 
of $K$ such that for each $\xi\in K_1$ the function $g(t)=u_t\xi$
extends to a bounded continuous
function $g(z)$ on the strip $-1/2 \le {\rm Im}\, z \le 0$ 
into $H$, holomorphic in $-1/2<{\rm Im}\, z < 0$, such that
$f(t-i/2)=j_1f(t)$ for $t\in {\Bbb R}$. Then $u_t=\delta^{it}$ and
$j_1=j$. 
\vskip .05in
\noindent \bf Proof. \rm (1) (cf [31]) If $\xi\in K$, then
$\xi=p\xi=(r+t)\xi=r^{\half}(r^{\half}+(I-r)^{\half}j)\xi$. Thus
$\xi=r^{\half}\eta$, where $\eta=(r^{\half}+(I-r)^{\half}j)\xi$. Set $f(z) =
(I-r)^{iz}r^{\half-iz}\eta$ for $-1/2 \le {\rm Im}\, z \le 0$.

\noindent (2) For $\xi\in K_1$, set $h(z)=(g(z),g(\overline{z}-i/2))$.
Then $h$ is continuous and bounded on $-1/2\le {\rm Im}\, z\le 0$,
holomorphic on $-1/2 <{\rm Im}\, z<0$. By uniqueness of analytic
extension, $u_t  f(z) = f(z+t)$ since they 
agree for $z$ real. Hence $h(z+t)=h(z)$, so that $h$ is constant on
lines parallel to the real axis and hence constant everywhere.
Since $h(-i/4)=\|g(-i/4)\|^2\ge 0$, it follows that
$h(0)\ge 0$, i.e.~$(j_1\xi,\xi)\ge 0$. Polarising, we get
$(j_1\xi,\eta)\in {\Bbb R}$ for all $\xi,\eta\in K$. 
Since $u_t$ leaves $K$ and $iK$ invariant, it follows that $u_t$
commutes with $e$ and $f$ and hence $\delta^{it}$. 
Now let $f(z)$ be the function corresponding to $\xi$ and
$\delta^{it}$. Define $k(z)=(g(z),jf(z))$ for $-\half \le {\rm Im}\,
z\le 0$. Then $k(t)=(u_t\xi,j\delta^{it}\xi)$ is real for $t\in {\Bbb
R}$ and
$k(t-i/2)=(j_1u_t\xi,j^2\delta^{it}\xi)=(j_1u_t\xi,\delta^{it}\xi)$ is
real for $t\in {\Bbb R}$. $k$ is   
bounded and continuous on $-\half \le {\rm Im}\, z\le 0$ and
holomorphic on $0<{\rm Im}\,z<\half$. By Schwartz's reflection
principle, $k$ extends to a holomorphic function on ${\Bbb C}$
satisfying $k(z+i)=k(z)$. This extension is bounded 
and therefore constant by Liouville's theorem. Hence
$k(t)=k(0)=k(-i/2)$. Thus $(u_t\delta^{-it}\xi,j\xi)
=(\xi,j\xi)=k(-i/2)=(j_1\xi,\xi)$. By polarisation it follows that 
$u_t=\delta^{it}$ and $j=j_1$, as required.

\vskip .1in
\bf \noindent 12. Modular operators and Takesaki devissage for von
Neumann algebras. \rm  The main application of the modular theory for 
a closed real subspace is when the subspace arises from a von Neumann
algebra with a vector cyclic for the algebra and its commutant.
Let $M\subset B(H)$ be a von Neumann algebra and let
$\Omega\in H$ (the ``vacuum vector'')  satisfy
$\overline{M\Omega}=H=\overline{M^\prime 
\Omega}$. The condition $\overline{M^\prime \Omega}=H$ is clearly
equivalent to the condition 
that $\Omega$ is {\it separating} for $M$, i.e.~$a\Omega=0$ iff $a=0$
for $a\in M$. If in addition $M$ and $H$ are ${\Bbb
Z}_2$--graded, then the graded commutant $M^q$ equals $\kappa M^\prime
\kappa^{-1}$ where the Klein transformation $\kappa$ is given by
multiplication by $1$ on the even part of $H$ and by $i$ on the odd
part; in this case we will always require that $\Omega$ be even.
Let $K=\overline{M_{{\rm sa}}\Omega}$, a closed real subspace of $H$.
\vskip .1in
\noindent \bf Lemma~1. \it $K+iK$ is dense in
$H$ and $K\cap iK = (0)$. 
\vskip .05in
\noindent \bf Proof. \rm $K+iK \supseteq
M\Omega = M_{{\rm 
sa}}\Omega + i M_{{\rm sa}}\Omega$, so $K+iK$ is dense. Now
$K^{\perp}\supseteq i M^\prime_{{\rm sa}}\Omega$, since for
$a\in M_{{\rm sa}}$, $b\in M^\prime_{{\rm sa}}$, we have
${\rm Re}\, (a\Omega, ib\Omega)
= {\rm Re}\, -i (ab\Omega,\Omega) =0$, 
because $(ab)^*=ab$ implies that $(ab\Omega,\Omega)$ is real. Hence 
$iK^{\perp}\supseteq M^\prime_{{\rm sa}} \Omega$. Thus
$K^{\perp} + iK^{\perp} \supseteq M^\prime\Omega$, so $K^{\perp}
 + iK^{\perp}$ is dense. So $K\cap iK=(K^{\perp} + iK^{\perp})^{\perp}
= (0)$. 
\vskip .1in
Let $\Delta^{it}$ and $J$ be the modular operators on $H$ associated
with $K=\overline{M_{\rm sa} \Omega}$. 
The main theorem of Tomita--Takesaki asserts that $JMJ=M^\prime$ and
$\Delta^{it} M \Delta^{-it} =M$. (General proofs can be
found in [7] or [31] for example; for hyperfinite von Neumann algebras
an elementary proof is given 
in [40], based on [31] and [15].) Once the theorem is known, the map
$x\mapsto Jx^*J$ gives an isomorphism between $M^{\rm op}$ ($M$ 
with multiplication reversed) and
$M^\prime$ and $\sigma_t(x)=\Delta^{it} x\Delta^{-it}$ gives a
one--parameter group of automorphisms of $M$. Our development,
however, does not logically 
require any form of the main theorem of Tomita--Takesaki; instead we 
verify it directly for fermions and deduce it for subalgebras
invariant under the modular group using a crucial result of Takesaki
(``Takesaki devissage''). 
\vskip .05in
\noindent \bf Lemma~2. \it If $JMJ\subseteq M^\prime$, then
$JMJ=M^\prime$.
\vskip .05in
\noindent \bf Proof (cf [31]). \rm Clearly
$J\Omega=\Omega$. If $A,B\in M^\prime_{\rm sa}$, then
$(JB\Omega,A\Omega)$ is real since $A\Omega, B\Omega $ lie in
$iK^\perp$ and $J$
is also the modular conjugation operator for $iK^\perp$.
Thus $(AJBJ\Omega,\Omega)=(JB\Omega,A\Omega)
=(A\Omega,JB\Omega)=(JBJA\Omega,\Omega)$.
By complex linearity in $A$ and conjugate--linearity in $B$, it follows that
$(AJBJ\Omega,\Omega)=(JBJA\Omega,\Omega)$ for all $A,B\in M^\prime$. Now
take $a,b\in M^\prime$, $x,y\in M$ and set $A=a$ and
$B=Jy^*JbJxJ.$ Since $JxJ,JyJ\in M^\prime$, $B$ lies in $M^\prime$. 
Hence $(JbJa x\Omega,y\Omega) =(aJbJ x\Omega,y\Omega)$.
Since $\overline{M\Omega}=H$, this implies that $aJbJ=JbJ a$.
Thus $JM^\prime J\subseteq M^{\prime\prime}=M$ and so
$JMJ=M^\prime$. 
\vskip .1in
\noindent \bf Corollary. \it If $A\subset B(H)$ is an Abelian von
Neumann algebra and $\Omega$ a cyclic vector for $A$, then
$\Delta^{it}=I$, $Ja\Omega =a^*\Omega$ and $JaJ=a^*$ for $a\in A$, and 
$A=JAJ=A^\prime$. 
\vskip .05in
\noindent \bf Proof. \rm Since $A\subset A^\prime$, $\Omega$ is
separating for $A$. Thus $Ja\Omega = a^*\Omega$ extends by continuity
to an antiunitary. If $a\in A_{\rm sa}$, the map $f(z)=a$ satisfies 
the KMS conditions for the trivial group and $J$, so they must be the
modular operators. Since $JAJ=A\subseteq A^\prime$, the last assertion
follows from the lemma. 
\vskip .1in

\noindent \bf Theorem (Takesaki devissage [35]). \it Let $M\subset
B(H)$ be a von Neumann algebra, $\Omega\in H$ cyclic for $M$ and
$M^\prime$ and $\Delta^{it}$, $J$ the corresponding modular operators.
Suppose that $\Delta^{it} M\Delta^{-it}=M$ and $JMJ=M^\prime$. 
If $N\subset M$ is a von
Neumann subalgebra such that $\Delta^{it}N\Delta^{-it}=N$, then
\vskip .02in
\noindent (a) $\Delta^{it}$ and
$J$ restrict  to the modular automorphism group $\Delta_1^{it}$ and
conjugation 
operator $J_1$ of $N$ for $\Omega$ on the closure $H_1$ of $N\Omega$.

\noindent (b) $\Delta_1^{it}N\Delta_1^{-it} =N$ and
$J_1NJ_1=N^\prime$.

\noindent (c) If $e$ is the
projection onto $H_1$, then $eMe=Ne$ and
$N=\{x\in M: xe=ex\}$ (the Jones relations [17]).

\noindent (d) $H_1=H$ iff $M=N$.

\vskip .05in

\noindent \bf Proof. \rm  (a) By KMS uniqueness,
$\Delta^{it}$ and $J$ restrict to $\Delta_1^{it}$ and $J_1$ on $H_1=eH$.

\noindent (b) It is clear that ${\rm Ad}\Delta_1^{it}$ normalises
$Ne=N_1$ on $H_1$. 
Now $J_1 NeJ_1=eJNJ e\subseteq eJMJe=eM^\prime e\subseteq eN^\prime
e=(eN)^\prime$. Thus $J_1 N_1J_1\subseteq N_1^\prime$. By Lemma~2,
$J_1N_1 J_1=N_1^\prime$. 

\noindent (c) Since $M^\prime \subset
N^\prime$ and$M^\prime = J M J$, this implies that
$M\subset JN^\prime J$. Compressing by $e$ we get
$eMe \subseteq eJ N^\prime J e = J e N^\prime eJ = J_1 e
N^\prime e J_1 = J_1(N\cdot e)^\prime J_1 = N\cdot e$.
But trivially $Ne \subseteq eMe$, so that $eMe = Ne$.
Clearly $N \subset \langle e\rangle^\prime$. Now suppose that $x\in
M$ commutes with $e$. Then $xe=ye$ for some $y\in N$. But then
$(x-y)e=0$, so that $(x-y)\Omega=0$. Since $\Omega$ is separating for
$M$, $x=y$ lies in $N$.

\noindent (d) Immediate from (c).
\vskip .1in

\vskip .1in

\noindent \bf 13. Araki duality and modular theory for Clifford
algebras. \rm We develop the abstract results implicit in
the work of Araki on the canonical commutation and anticommutation
relations ([1],[2]). This reduces the computation of the modular
operators for Clifford algebras to ``one particle
states'', i.e.~to the prequantised Hilbert space. We first recall that the 
assignment $H\mapsto \Lambda(H)$ 
defines a {\it 
functor} from the additive theory of Hilbert spaces and contractions 
to the multiplicative theory of Hilbert spaces and contractions. 
A {\it contraction} $A:H_1\rightarrow H_2$ between two Hilbert spaces is
a bounded linear map with $\|A\|\le 1$. We define $\Lambda(A)$ to be
$A^{\otimes k}$ on $\Lambda^k(H_1)\subset H_1^{\otimes k}$. Then $\Lambda(A)$
gives a bounded linear operator from $\Lambda(H_1)$ to $\Lambda(H_2)$ with
$\|\Lambda(A)\|\le 1$. Clearly if $\|A\|, \|B\|\le 1$,
then $\Lambda(AB)=\Lambda(A)\Lambda(B)$. Also
$\Lambda(A)^*=\Lambda(A^*)$, so if 
$A$ is unitary, 
then so too is $\Lambda(A)$. Similarly, if $H_1=H_2=H$, then if $A$ is
self--adjoint or positive, so too is $\Lambda(A)$. In particular if $A=UP$
is the polar decomposition of $A$ with $U$ unitary, then
$\Lambda(A)=\Lambda(U)\Lambda(P)$ is the polar decomposition of $\Lambda(A)$
by uniqueness. 
Moreover $\Lambda(A^{it}) = \Lambda(A)^{it}$ if $A$ is in addition
positive (note 
that $(A^{it})^{\otimes k} =(A^{\otimes k})^{it}$). Similarly 
every conjugate--linear contraction $T$ induces an operator
$\tilde \Lambda(T) (\xi_1\wedge \xi_2 \wedge \cdots \wedge \xi_n) =
T\xi_n \wedge T\xi_{n-1} \wedge  \cdots \wedge T\xi_1$. Note that
$\tilde{\Lambda}(T)=\kappa^{-1} \Lambda(iT)$, where $\kappa$ is the
Klein transformation. If $T=UP$ is the polar
decomposition of $T$ with $U$ a conjugate--linear unitary, then
$\tilde{\Lambda}(T) 
=\tilde{\Lambda}(U)\Lambda(P)$ is the 
polar decomposition of $\tilde{\Lambda}(T)$. If $U$ is a linear or
conjugate--linear unitary, then it is easy to check that
$\Lambda(U)a(\xi)\Lambda(U)^*=a(U\xi)$ 
and $\Lambda(U)c(\xi)\Lambda(U)^*=c(U\xi)$.

Let $H$ be a complex Hilbert space and $K\subset H$ a closed real
subspace of $H$ such that $K\cap
iK=(0)$ and $K+iK$ is dense in $H$. For $\xi\in H$ let  $a(\xi)$
denote exterior multiplication by $\xi$ and let $c(\xi)=a(\xi) +
a(\xi)^*$ denote Clifford multiplication. Thus
$c(\xi)c(\eta)+c(\eta)c(\xi)= 
2{\rm Re}\,(\xi,\eta)$. Since the *--algebra generated by the
$a(\xi)$'s acts irreducibly on $\Lambda H$ and since $a(\xi) = (c(\xi) - i
c(i\xi))/2$, the $c(\xi)$' s act irreducibly on $\Lambda H$. 
\vskip .05in 
\noindent \bf Lemma. \it If $M(K)$ is the von Neumann algebra
generated by the $c(\xi)$'s ($\xi\in K$), then $\Omega$ is cyclic
for $M(K)$.  
\vskip .05in
\noindent \bf Proof. \rm Let $H_0=\overline{M(K)\Omega}$  and assume by
induction that all forms of 
degree $N$ or less lie in $H_0$. Let $\omega$ be an $N$--form and take
$f\in K$. Then $f\wedge \omega =c(f)\omega - a(f)^*\omega$ , so that
$f\wedge \omega\in H_0$. Since $K+iK$ is dense in $H$ and $H_0$ is a
complex subspace of $\Lambda H$, it follows that $\xi\wedge\omega\in
H_0$ for all $\xi\in H$. Hence $H_0$ contains all $(N+1)$--forms. 
\vskip .05in
Since $\Omega$ is cyclic for $M(K^\perp)$, which lies in the graded
commutant of $M(K)$, it follows that $\Omega$ is cyclic and separating for
$M(K)$. Let $R$, $T$, $\Delta^{it}=(I-R)^{it}R^{-it}$ and $J$ be the
corresponding modular operators for $M(K)$ and $\Omega$. 
\vskip .05in
\noindent \bf Theorem. \it (i) $J=\tilde\Lambda(j)=\kappa^{-1}\Lambda(ij)$,
$\Delta^{it}=\Lambda(\delta^{it})$, where $j$ and $\delta^{it}$ are
the modular operators for $K$. 

\noindent (ii) For $\xi\in H$,
$\Delta^{it}c(\xi)\Delta^{-it} = 
c(\delta^{it}\xi)$ and $\kappa J c(\xi) J \kappa^{-1} = c(ij\xi)$, 
where $\kappa$ is the Klein transformation.

\noindent (iii) $M(K^\perp)$ is the graded commutant of $M(K)$ and
$M(K)^\prime =JM(K)J$ (Araki
duality). 
\vskip .05in
\noindent \bf Remark. \rm For another proof, analogous to that of [23] for
bosons and the canonical commutation relations, see [40]. 

\vskip .05in
\noindent \bf Proof (cf [2]). \rm Let $\delta^{it}$ and $j$ be
the modular operators associated with the closed real subspace
$K\subset H$. Let $S$ be the conjugate--linear operator on
$\pi_P({\rm Cliff}_{\Bbb R}(K))\Omega$ defined by
$Sa\Omega = a^*\Omega$ for $a\in M=\pi_P({\rm Cliff}_{\Bbb R}(K))$.
This is well--defined, because $\Omega$ is separating for 
$M$. Thus
$Sc(\xi_1)\cdots c(\xi_n)\Omega = c(\xi_n)\cdots c(\xi_1)\Omega$ for
$\xi_i\in  K$. If the $\xi_i$'s are orthogonal, it follows that
$S\xi_1\wedge \cdots \wedge \xi_n = \xi_n\wedge \cdots \wedge \xi_1$.
Since any finite dimensional subspace of $K$ admits an orthonormal
basis, this formula holds by linearity for arbitrary
$\xi_1,\dots,\xi_n\in K$. Since $S$ is conjugate--linear, it follows that
for $\xi_i,\eta_i\in K$ we have
$S\,(\xi_1+i\eta_1)\wedge \cdots \wedge (\xi_n+i\eta_n) =
(\xi_n-i\eta_n) \wedge \cdots \wedge (\xi_1-i\eta_1)$.

Let $J=\tilde {\Lambda}(j) =\kappa^{-1} \Lambda(ij)$ and
$\Delta^{it}=\Lambda(\delta^{it})$. Clearly
$\Delta^{it}J=\Delta^{it}J$ and $\Delta^{it}$ preserves
$\overline{M_{\rm sa}\Omega}$. To check the KMS condition, it 
suffices to show that for $x\in M\Omega$, the function
$F(t)=\Delta^{it}x$ extends to a 
bounded continuous function on $-\half \le {\rm Im}\, z\le 0$,
holomorphic on the interior, with $F(t-i/2)=JSF(t)$. We may assume
that $x=(\xi_1+i\eta_1)\wedge \cdots \wedge (\xi_n+i\eta_n) $ with
$\xi_i,\eta_i\in K$. For each $i$, let $f_i(z)$ be continuous bounded
function on $-\half\le {\rm Im}\, z\le 0$, holomorphic in the
interior, $f_i(t) =\delta^{it}(\xi_i +i\eta_i)$ and $f_i(t-i/2)
=j\delta^{it}(\xi_i-i\eta_i)$.
Let $F(z) = f_1(z) \wedge \cdots \wedge f_n(z)$. Then $F(z)$ is
bounded and continuous on $-\half \le {\rm Im}\, z\le 0$, holomorphic
in the interior, and $F(t)=\Delta^{it} x$. Now 
$F(t-i/2) =f_1(t-i/2) \wedge \cdots \wedge f_n(t-i/2) 
= j \delta^{it}(\xi_1-i\eta_1) \wedge \cdots \wedge
j\delta^{it}(\xi_n-i\eta_n) = \tilde{\Lambda}(j) SF(t)=JSF(t)$.
Thus $F(t-i/2)=JSF(t)$ as required. This proves (i) and (ii) follows
immediately. To prove (iii), note that $ij(K)=K^\perp$, so that 
$M(K^\perp)=\kappa JM(K)J\kappa^{-1}$ by this covariance relation.
But $M(K^\perp)\subseteq M(K)^q=\kappa M(K)^\prime \kappa^{-1}$. Thus
$JM(K)J\subseteq M(K)^\prime$, so the result follows from Lemma~2 in
Section~12.

\vskip .1in

\noindent \bf 14. Prequantised geometric modular theory. \rm In this
section we compute the prequantised modular operators corresponding to
fermions on the circle, using an analytic continuation argument
obtained jointly with Jones. This argument is
prompted by the KMS condition and may be regarded as a prequantised
analogue of the 
computations of Bisognano and Wichmann [4]. (For another approach, see
[40].) Let $H$ be the  
complex Hilbert space $L^2(S^1,V)$ where $V={\Bbb C}^N$. We give $H$ 
a new complex structure by defining multiplication by $i$ as $i(2P-I)$, where
$P$ is the orthogonal projection onto Hardy space $H^2(S^1,V)$. Let
$I$ be the upper semicircle  
and let $K=L^2(I,V)$, a real closed subspace of $H_P$. The real 
orthogonal projection onto $K$, regarding $H$ as a real inner product space, 
is given by $Q$, multiplication by $\chi_I$.
\vskip .05in
\noindent \bf Theorem. \it (a) $K\cap iK=(0)$ and $K+iK$ is dense in
$H_P$.

\noindent (b) $K^\perp =L^2(I^c,V)$. 

\noindent (c) $j =-i(2P-I)F$ where $Ff(z)=z^{-1}f(z^{-1})$ is the flip,
and $\delta^{it}=u_t$, where $(u_tf)(z)=\break (z \sinh\pi t
+\cosh\pi t)^{-1}  
f(z\cosh\pi t + \sinh\pi t/z\sinh\pi t +\cosh\pi t)$.
\vskip .05in

\noindent \bf First proof. \rm (a) It 
suffices to show that $P$ and $Q$  
are in general position. Now conjugation by $r_\pi$ takes $Q$ onto 
$I-Q$ and fixes $P$ while conjugation by the flip $Vf(z) =z^{-1}f(z^{-1})$
takes $Q$ onto $I-Q$ and $P$ onto $I-P$. Thus it will suffice to show that
$PH\cap QH=(0)$. Suppose that the negative Fourier coefficients of 
$f\in L^2(I,V)$ all vanish. Then so do those of $\psi\star f$
for any $\psi\in C^\infty(S^1)$. But $\psi\star f\in C^\infty(S^1,V)$
is the boundary value of a holomorphic function. If $\psi$ is
supported near $1$, $\psi\star f$ vanishes in a subinterval of $I^c$
and therefore must vanish identically (since $\psi\star f$ can be
extended by reflection across this subinterval). Since $\psi\star f$
and $f$ can 
be made arbitrarily close in $L^2(S^1,V)$, we must have $f=0$.

\vskip .02in
\noindent (b) The real orthogonal complement of $L^2(I,V)$ in $L^2(S^1,V)$ 
is clearly $L^2(I^c,V)$.
\vskip .02in

\noindent (c) Let
$K_1\subset K$ be the dense subset of $QH$ 
consisting of functions $Qp$ where $p$ is the restriction of a
polynomial in $e^{i\theta}$. We must show that the
map $f(t)=u_tQp$ extends to a bounded continuous function $f(z)$ on
the closed strip $-1/2\le {\rm Im}\, z\le 0$, holomorphic in the open
strip with $f(t-i/2)=jf(t)$ for $t\in {\Bbb R}$. 
Now $f(t) = Pu_t
Qp+(I-P)u_tQp$. Because of  
the modified complex structure on $H=PH\oplus (I-P)H$, we have to
extend $f_1(t)=Pu_t Qp$ to a holomorphic function with values in $PH$ and 
$(I-P)u_t Qp$ to an antiholomorphic function with values in $(I-P)H$. 
Note that if
$\theta\in [0,\pi]$ and $-3/4 < {\rm Im}\, z< 1/2$, the function $s_z
e^{i\theta} + c_z$ is non--zero, where $s_z=\sinh \pi z$ and
$c_z=\cosh \pi z$. For $-3/4 < {\rm Im}\, z< 1/2$, let
$p_z(e^{i\theta}) = (s_ze^{i\theta} + c_z)^{-1}p(c_z e^{i\theta}
+s_z/s_z e^{i\theta} +c_z)$. Then $Qp_z$ is holomorphic for such $z$,
so $f_1(z)=PQp_z$ gives a holomorphic extension of $f_1$ to 
$-3/4 < {\rm Im}\, z< 1/2$. 
Next note that $f_2(t) = -(I-P)u_t(I-Q)p$,
since $(I-P)p=0$. Set $f_2(z) = -(I-P)(I-Q)p_{\overline{z}}$. This
gives an
antiholomorphic extension of $f_2$ to $-3/4 < {\rm Im}\, z <1/4$, because
$s_{\overline{z}} e^{i\theta} +c_{\overline{z}}$ does not vanish for
$\theta\in [-\pi,0]$. Thus $f(z)=f_1(z)+f_2(z)$ is a holomorphic
function from $-3/4 < {\rm Im}\, z< 1/2$ into $H$. It equals $f(t)$
for $t\in {\Bbb R}$. 
If we show that $f(t-i/2)=jf(t)$, then $f(z)$
will be bounded for ${\rm Im}\,z=0$ or $-1/2$ and hence, by the
maximum modulus principle, on the strip $-1/2\le {\rm Im}\, z\le 0$.
Now $jf(t)
=-i(2P-I)Ff(t)=-iPQFp_t +i(I-P)(I-Q)Fp_t$. Since $s_{t\pm i/2} =\pm i
c_t$ and $c_{t\pm i/2}=\pm is_t$, we have $p_{t\pm
i/2}=\mp i Fp_t$. Hence
$f_1(t-i/2) = -iPQFp_t$ and $f_2(t-i/2)= i(I-P)(I-Q)Fp_t$, so that
$f(t-i/2)=jf(t)$ as required. 
\vskip .05in
\noindent \bf Second proof. \rm  Let $U:L^2(S^1,V)\rightarrow L^2({\Bbb
R},V)$, $Uf(x)=(x-i)^{-1} f(x+i/x-i)$ be the unitary induced by the
Cayley transform. Let $V:L^2({\Bbb R},V)\rightarrow L^2({\Bbb
R},V)\oplus L^2({\Bbb R},V)$ be the unitary defined by
$Vf=(\widehat{f_+},\widehat{f_-})$, where $\widehat{g}$ denotes the
Fourier transform of $g$ and $f_\pm(t) =e^{t/2}f(\pm e^t)$. Let
$W=VU:L^2(S^1,V)\rightarrow L^2({\Bbb R},V)\oplus L^2({\Bbb R},V)$.
If $e_n(\theta)=e^{in\theta}$, it is easy to check that
$We_0 =(g_+,g_-)$ and $We_{-1}=(-g_-,-g_+)$ where $g_\pm(x)=\pi^\half (i\pm
1) e^{\pm \pi x/2} (1+e^{\pm 2\pi x})^{-1}$.

Clearly $WQW^*$ is the projection onto the first summand $L^2({\Bbb R},V)$.
Now $Uu_t U^*=v_{2\pi t}$, where $(v_s f)(x)=e^{s/2} f(e^s x)$; and
$Vv_s V^*=m(e_s)$, where $e_s(t)=e^{ist}$ and $m(e_s)$ is the
corresponding multiplication operator (acting diagonally). Hence
$Wu_tW^*=m(e_{2\pi t})$. These operators generate a copy of
$L^\infty({\Bbb R})$ on $L^2({\Bbb R})$, which by the corollary to
Lemma~2 in section~12 equals its own commutant on $L^2({\Bbb R})$. On the
other hand $P$ commutes with $u_t$ and ${\rm End}\, V$, so that
$WPW^*$ lies in the commutant of the $m(e_{2\pi t})$'s and ${\rm
End}\, V$. Hence $WPW^*=\pmatrix{m(a) & m(b) \cr m(c) & m(d) \cr}$
with $a,b,c,d\in L^\infty({\Bbb R})$. But $Pe_0=e_0$ and $Pe_{-1}=0$.
Transporting these equations by $W$, we get equations for $a,b,c,d$
which can be solved to yield $a(x)=(1+e^{2\pi x})^{-1}$, $b(x)=-c(x) =
ie^{\pi x} (1+e^{2\pi x})^{-1}$ and $d(x)=e^{2\pi x}(1+e^{2\pi
x})^{-1}$.

These formulas show that $WQW^*$ and $WPW^*$ are in general position,
so (a) follows. (b) is clear, since $L^2(I,V)^\perp = L^2(I^c,V)$.
To prove (c), note that $e=Q$ and 
$f=(2P-I)Q(2P-I)$, so that $r=PQP\oplus P^\perp Q P^\perp$ and $I-r=
PQ^\perp P \oplus P^\perp Q P^\perp$. 
Remembering that $r^{it}$ and $(I-r)^{it}$ must be defined using the
complex structure $i(2P-I)$, we get 
$(I-r)^{it}r^{-it}=(I-A)^{it} A^{-it}$, where $A=PQP\oplus P^\perp
Q^\perp P^\perp = QPQ \oplus Q^\perp P^\perp Q^\perp$. 
Hence $WAW^*=m(a)$ and $W\delta^{it}W^*=m((1-a)^{it}a^{-it})=m(e_{2\pi
t}) =Wu_tW^*$, 
so that $\delta^{it}=u_t$. Finally  
$t=(e-f)/2=(2P-I)(QP-PQ)$.
Now $W(QP-PQ)W^*= W(QPQ^\perp - Q^\perp PQ)W^*=\pmatrix{0 &
m(b)\cr 
m(b) & 0\cr}$ so that $j=-i(2P-I)F_1$ where $WF_1W^*=\pmatrix{ 0 & -I\cr
-I & 0\cr}$.  Now $UFU^*=F^\prime$, where $(F^\prime f)(x)=-f(-x)$,
so that $WFW^*=VF^\prime V^*=\pmatrix{0 & -I \cr -I & 0\cr}$. Hence
$F_1=F$, as required.

\vskip .1in

\noindent \bf 15. Haag--Araki duality and geometric modular theory
for fermions on the circle. \rm 
Let $H=L^2(S^1)\otimes V$ with $V={\Bbb C}^N$
and let $P$ be the orthogonal projection onto the Hardy space
$H^2(S^1)\otimes V$. Let $\pi_P$ denote the corresponding irreducible
representation of ${\rm Cliff}(H)$ on fermionic Fock space ${\cal
F}_V$. For any interval $J\subset S^1$, let
$M(J)\subset B({\cal F}_V)$ be
the von Neumann algebra generated by the operators $\pi_P(a(f))$ with $f\in
L^2(J,V)$. Our main result was obtained jointly with Jones ([18],
[40]); it follows almost immediately from the previous sections. 

\vskip .1in

\noindent \bf Theorem. \it Let $I$ denote the upper semicircle with
complement $I^c=S^1\backslash\overline{I}$.  

\noindent (a) The vacuum vector $\Omega$ is cyclic
and separating for $M(I)$.

\noindent (b) (Haag--Araki duality) $M(I^c)$ is the graded commutant
of $M(I)$ and $JM(I)J=M(I)^\prime$, where $J$ is the modular
conjugation with respect to $\Omega$.

\noindent (c) (Geometric modular group) Let $I\subset S^1$ be the
upper semi--circle. The
modular group $\Delta^{it}$ of $M(I)$ with respect to the
vacuum vector $\Omega$ is implemented by $u_t$, where $(u_t f)(z)
=(z \sinh\pi t
+\cosh\pi t)^{-1}  
f(z\cosh\pi t + \sinh\pi t/z\sinh\pi t +\cosh\pi t)$ is the M\"obius
flow fixing the endpoints of $I$. In particular
$\Delta^{it}\pi_P(a(f)) \Delta^{-it} =\pi_P(a(u_t f))$ for $f\in H$. 

\noindent (d) (Geometric modular conjugation) If $\kappa$ is the Klein
transformation, then the antiunitary $\kappa J$ is implemented by $F$,
where $Ff(z)=z^{-1}f(z^{-1})$ is the flip. In particular
$J\pi_P(a(f))J=\kappa^{-1} \pi_P(a(Ff))\kappa$ for $f\in H$.

\vskip .05in
\noindent \bf Remark. \rm Analogous results hold when $I$ is replaced
by an arbitrary interval $J$. This follows immediately by transport of
structure using the canonically quantised action of $SU(1,1)$. 

\vskip .05in
\bf \noindent Proof. \rm If $H_P= PH \oplus \overline{P^\perp H}$ ($H$
with multiplication by $i$ given by $i(2P-I)$), then
${\cal F}_V=\Lambda H_P$ and $\pi_p(a(f))=a(Pf) + a(\overline{P^\perp
f})^*$ on
$\Lambda H_P$ for $f\in H$. Hence $\pi_P(a(f)+a(f)^*)=c(Pf) +
c(\overline{P^\perp f})=c(f)$ for $f\in H$. Now $M(I)$ coincides with
the von Neumann algebra generated by $\pi_P(a(f)+a(f)^*)$ for $f\in
L^2(I,V)$. It therefore may be identified with the von Neumann algebra
generated by the $c(f)$ with $f\in K=L^2(I,V)$, a closed real subspace
of $H_P$. From Section~13, the vacuum vector $\Omega$ is
cyclic for $M(I)$ and $JM(I)J=M(I)^\prime=\kappa^{-1}M(I^c)\kappa$, since
$L^2(I,V)^\perp =L^2(I^c,V)$. From Section~14, we see that
$\Delta^{it}$ is the canonical quantisation 
of $u_t$ and the antiunitary $\kappa J$ is the canonical quantisation of
$F$. Finally the relations
$\Delta^{it}c(f)\Delta^{-it}=c(u_t f)$ and $\kappa J c(f)J\kappa^{-1}
= c(Ff)$ for $f\in H_P$ immediately imply that $\Delta^{it}\pi_P(a(f))
\Delta^{-it} =\pi_P(a(u_t f))$ and $J\pi_P(a(f))J=\kappa^{-1}
\pi_P(a(Ff))\kappa$ for $f\in H$.

\vskip .1in

\noindent \bf 16. Ergodicity of the modular group. \rm  
\vskip .1in
\noindent \bf Proposition. \it The
action $\Lambda(u_t)^{\otimes k}$ of ${\Bbb R}$ on $(\Lambda
H_P)^{\otimes k}$ is ergodic,  
i.e.~has no fixed vectors apart from multiples of the vacuum vector
$\Omega^{\otimes k}$.
\vskip .05in
\noindent \bf Proof. \rm First note that the action $u_t$ of ${\Bbb R}$ on
$L^2({\Bbb T})$  is unitarily
equivalent to the direct sum of two copies of the left regular
representation. In fact the unitary equivalence between
$L^2({\Bbb T})$ and $L^2({\Bbb R})$ 
induced by the Cayley transform 
$Uf(x) =(x-i)^{-1} f(x+i/x-i)$ carries $u_t$ onto the scaling action
$v_{2\pi t}$ of
${\Bbb R}$ on $L^2({\Bbb R})$, where $(v_s f)(x)=e^{s/2}f(e^sx)$. 
For $f\in L^2({\Bbb R})$ define $f_\pm\in
L^2({\Bbb  
R})$ by $f_\pm(t) = e^{t/2} f(\pm e^t)$ and set $W(f) = (f_+,
f_-)$. Thus $W$ is an unitary between $L^2({\Bbb R})$ and $ L^2({\Bbb R})
\oplus L^2({\Bbb R})$.  This unitary carries the scaling action of
${\Bbb R}$ onto the direct sum of two copies of the
regular representation. 

Thus $L^2({\Bbb T})\cong L^2({\Bbb R})\oplus L^2({\Bbb R})$ as a
representation of 
${\Bbb R}$. Now $H=L^2({\Bbb T},V)$ is a direct sum of copies
of $L^2({\Bbb T})$. On the other hand
$\overline{L^2({\Bbb R})} \cong L^2({\Bbb R})$ (by conjugation), it
follows that both $H$ and $\overline{H}$ are subrepresentations of a
direct sum of copies of $L^2({\Bbb R})$. But
$H_P=PH\oplus (I-P)H$ is a subrepresentation of $H\oplus
\overline{H}$, so that $H_P$ is unitarily equivalent to a
subrepresentation of  
$L^2({\Bbb R})\otimes {\Bbb C}^n$ for some $n$.

Thus the action of ${\Bbb
R}$ on $(\Lambda H_P)^{\otimes k}=\Lambda(H_P\otimes {\Bbb C}^k)$ is
unitarily equivalent 
to a subrepresentation of
${\Bbb R}$ on $\Lambda H_1$, where $H_1=L^2({\Bbb R})\otimes {\Bbb C}^m$
for some $m\ge 2$. It therefore suffices to check that ${\Bbb
R}$ has no fixed vectors in $\lambda^k H_1$ for $k\ge 1$,
since the action of ${\Bbb R}$ preserves degree. 

Now $\lambda^k
H_1\subset H_1^{\otimes k}$. On the other hand if $t\mapsto \pi(t)$ is any
unitary representation of ${\Bbb 
R}$ on $H$ and $\lambda(t)$ is the left regular representation on
$L^2({\Bbb R})$, then $\lambda\otimes \pi$ and $\lambda\otimes I$ are
unitarily equivalent: the unitary $V$, defined by
$Vf(x)=\pi(x)f(x)$ for $f\in 
L^2({\Bbb R},H)=L^2({\Bbb R})\otimes H$, gives an intertwiner. 
It follows that $H_1^{\otimes k}$ is
unitarily equivalent to a direct sum of copies of the left regular
representation. Hence $\lambda^k H_1$ is unitarily equivalent to a
subrepresentation of a direct sum of copies of the left regular
representation. Since the Fourier transform on $L^2({\Bbb R})$
transforms $\lambda(t)$ into multiplication by $e_t(x)=e^{itx}$, 
no non--zero vectors in $L^2({\Bbb R})$ are fixed by $\lambda$.
Hence there are no non--zero vectors in
$\lambda^k H_1$ fixed by ${\Bbb R}$ for $k\ge 1$, as claimed.

\vskip .1in
\noindent \bf Corollary. \it The modular group acts ergodically on the
local algebra $M(I)=\pi_P({\rm Cliff}(L^2(I,V)))^{\prime\prime}$,
i.e.~it fixes only the scalar operators. In particular $M(I)$ must be
a factor [in fact a type III$_1$ factor]. 
\vskip .05in
\noindent \bf Proof. \rm Suppose that $x\in M(I)$ is fixed by the
modular group. Then $x\Omega$ is fixed by the modular group, so that
$x\Omega=\lambda\Omega$ for $\lambda\in {\Bbb C}$. Since $\Omega$ is
separating for $M(I)$, this forces $x=\lambda I$. Since the modular
group fixes the centre, $M(I)$ must be a factor.

\vskip .1in

\noindent \bf 17. Consequences of modular theory for local loop groups. \rm
Using only Haag--Araki duality for 
fermions and Takesaki devissage, we
establish several important properties of the von Neumann algebras
generated by local loop groups in positive energy representations.
These include Haag 
duality in the vacuum representation, local equivalence, the fact that
local algebras are factors and
a crucial irreducibility 
property for local loop groups. This irreducibility result will be
deduced from a von Neumann density result, itself a consequence of a
generalisation of Haag duality; it can also be deduced from a 
careful study of the topology on the loop group induced by its
positive energy representations.

Let $L_IG$ be the 
local loop group consisting of loops concentrated in 
$I$, i.e.~loops equal to $1$ off $I$, and let ${\cal L}_I G$ be the
corresponding subgroup of ${\cal L}G$. We need to know in what sense
these subgroups generate $LG$.
\vskip .1in
\noindent \bf Covering lemma. \it If $S^1=\bigcup_{k=1}^n I_k$, then
$LG$ is 
generated by the subgroups $L_{I_k}G$. 
\vskip .05in
\noindent \bf Proof. \rm By the exponential lemma we just have to prove
that every exponential $\exp(X)$ lies in the group generated by
$L_{I_k}G$. Let $(\psi_k)\subset C^\infty(S^1)$ be a smooth partition
of the identity subordinate to $(I_k)$. Then $X=\sum \psi_k\cdot X$, so
that $\exp(X)=\exp(\psi_1 \cdot X) \cdots \exp(\psi_n\cdot X)$ with
$\exp(\psi_k\cdot X)\in L_{I_k}G$.
\vskip .1in
Let $\pi:LSU(N)\rightarrow PU({\cal F}_V)$ be the 
basic representation of $LSU(N)$, so that
$\pi(g) \pi_P(a(f))
\pi(g)^* =  \pi_P(a(g\cdot f))$ for $f\in L^2(S^1,V)$
and $g\in LSU(N)$. Let $\pi_i$ be an irreducible
positive energy representation of level $\ell$. Haag--Araki duality
and the fermionic construction of $\pi_i$ imply
that operators in $\pi_i(L_I G)$ and $\pi_i(L_{I^c}G)$, defined up to a
phase, actually commute (``locality''):
\vskip .1in
\noindent \bf Proposition (locality). \it For any positive energy
representation $\pi_i$, we have
$\pi_i(g)\pi_i(h)\pi_i(g)^*\pi_i(h)^*=I$ for $g\in {\cal L}_ISU(N)$ and $h\in
{\cal L}_{I^c}SU(N)$. 
\vskip .05in

\noindent \bf Proof. \rm As above let $M(I)\subset B({\cal F}_V)$ be 
the von 
Neumann algebra generated by fermions $a(f)$ with $f\in L^2(I,V)$. 
Since $\pi(g)$ commutes with $M(I^c)$ and
is even, it must lie in $M(I)$ by Haag--Araki duality. Similarly
$\pi(h)$ lies in $M(I)$. Since they are both even operators they must
therefore commute. Clearly this result holds also with $\pi^{\otimes
\ell}$ in place of $\pi$ and passes to any subrepresentation $\pi_i$
of $\pi^{\otimes \ell}$. 
\vskip .1in
The embedding of
$LSU(N)$ in $LSU(N\ell)$ gives a projective representation $\Pi$ on ${\cal
F}_W$ where $W=({\Bbb C}^N)^{\otimes \ell}$. Now ${\cal F}_W$ is
can naturally be identified with ${\cal F}_V^{\otimes \ell}$ and under this
identification $\Pi=\pi^{\otimes \ell}$. Let $M=\pi_P({\rm
Cliff}(L^2(I,W)))^{\prime\prime}$ and let
$N=\pi^{\otimes\ell}({\cal L}_ISU(N))^{\prime\prime}=\Pi({\cal L}_I
SU(N))^{\prime\prime}$, so that $N\subset M$. The operators $u_t$ and
$F$ lie in $SU_\pm (1,1)$ so are compatible
with the central extension ${\cal L}G$ introduced in section~5. 
It follows immediately that
$N$ is invariant under the modular group of $M$.
In order to identify $\overline{N\Omega}$ we need a
preliminary result. 
\vskip .1in

\noindent \bf Reeh--Schlieder theorem. \it Let $\pi$ be an irreducible
positive energy 
projective representation of $LG$ on $H$ and let $v$ be a finite energy
vector(i.e.~an eigenvector for rotations). Then the linear span of
$\pi({\cal L}_IG)v$ is dense in $H$.
\vskip .1in
\noindent \bf Proof (cf [30]). \rm It suffices to show that if $\eta\in H$
satisfies $( \pi(g)v, \eta)=0$ for all $g\in {\cal L}_I
G$, then $\eta=0$. We start by using the positive energy condition to
show that this identity holds for all $g\in LG$. For $z_1,\dots,z_n\in
{\Bbb T}$ and $g_1,\dots,g_n \in {\cal L}_J G$, where $J\subset\subset I$,
consider $F(z_1,\dots,z_n)=(  R_{z_1} \pi(g_1) R_{z_2} \pi(g_2)
\cdots R_{z_n} \pi(g_n)v,\eta)$.
This vanishes if all the $z_i$'s are sufficiently close to $1$. Now
freeze $z_1,\dots,z_{n-1}$. As a function of $z_n$, the positive
energy condition implies that the function $F$
extends to a continuous function on the closed unit disc, holomorphic
in the interior and vanishing on the unit circle near $1$. By the Schwarz
reflection principle, $F$ must vanish identically in $z_n$. Now freeze
all values of $z_i$ except $z_{n-1}$. The same argument shows that $F$
vanishes for all values of $z_{n-1}$, and so on. After $n$ steps, we
see that $F$ vanishes for all values of $z_i$ on the unit circle. Thus
$(\pi(g)v,\eta) =0$ for all $g$ in the group
generated by ${\cal L}_JG$ and its rotations, i.e.~the whole group
${\cal L}G$. Therefore, since $\pi$
is irreducible, $\eta=0$ as required.
\vskip .1in
We may now apply Takesaki devissage with the following
consequences. 
\vskip .1in
\noindent \bf Theorem~A (factoriality). \it
$N=\pi^{\otimes \ell}({\cal L}_IG)^{\prime\prime}$, and hence each
isomorphic 
$\pi_i({\cal L}_I 
G)^{\prime\prime}$, is a factor.

\vskip .05in

\noindent \bf Proof. \rm By Takesaki devissage, $N$
has ergodic modular group and therefore must be a factor.
 If $p_i$ is a projection in
$\pi^{\otimes\ell}(LG)^{\prime}
\subset \pi^{\otimes\ell}({\cal L}_IG)^{\prime}$ corresponding to the 
irreducible positive energy representation $H_i$, then 
$\pi_i({\cal L}_IG)^{\prime\prime}$ is isomorphic to
$\pi^{\otimes \ell}({\cal L}_IG)^{\prime\prime}p_i\cong N$ and is therefore also a
factor. 
\vskip .1in

\noindent \bf Theorem~B (local equivalence). \it For every positive
energy representation $\pi_i$ of level $\ell$, there is a unique
*--isomorphism $\pi_i:\pi_0({\cal L}_I 
G)^{\prime\prime} \rightarrow \pi_i({\cal L}_IG)^{\prime\prime}$ sending
$\pi_0(g)$ to $\pi_i(g)$ for all $g\in {\cal L}_IG$. 
If ${\cal X}={\rm Hom}_{{\cal L}_IG}(H_0,H_i)$, then ${\cal X}\Omega$ is
dense in $H_i$ and $\pi_i(a)T=Ta$ for all $T\in {\cal X}$ and $a\in
\pi_0({\cal L}_IG)^{\prime\prime}$. If ${\cal X}_0$ is a subspace of
${\cal X}$ with  ${\cal X}_0H_0$ dense in $H_i$, then $\pi_i(a)$ is
the unique operator $b\in B(H_i)$ such $bT=Ta$ for all $T\in {\cal
X}_0$.
\vskip .05in
\noindent \bf Proof. \rm This is immediate from the proposition in
Section~10, since $\pi_0$ and $\pi_i$ are subrepresentations of the factor
representation $\pi^{\otimes \ell}\otimes I$. Since ${\cal X}={\cal
X}\pi_0({\cal L}_{I}G)$ and $\Omega$ is 
cyclic for $\pi_0({\cal L}_I G)$, it follows that $\overline{{\cal
X}\Omega }=\overline{{\cal X}H_0}=H_i$.
\vskip .1in
\noindent \bf Remarks. \rm Note that, if $p_i,p_j$ are projections onto
copies of $H_i,H_j$ in ${\cal F}_W$, explicit intertwiners
$H_j\rightarrow H_i$ are given by compressed fermi fields $p_i a(f)
p_j$ with $f$ supported in $I^c$; these are essentially the smeared
vector primary fields that we study in Chapter~IV. Theorem B is a
weaker version of the much 
stronger result that the restrictions of $\pi_0$ and $\pi_i$ to ${\cal
L}_IG$ are unitarily equivalent. This follows because $\pi^{\otimes
\ell}$ restricts to a type III factor representation of ${\cal L}_I
G$ (because the modular group is ergodic). Thus any non--zero
subrepresentations are unitarily equivalent. Local equivalence may 
also be proved more directly using an argument of
Borchers [6] to show that the local algebras are ``properly infinite''
instead of type~III (see [40] and [41]).

\vskip .1in
\noindent \bf Theorem~C (Haag duality). \it If $\pi_0$ is the vacuum
representation at level $\ell$, then $\pi_0({\cal L}_IG)^{\prime\prime}
=\pi_0({\cal L}_{I^c}G)^\prime$. The 
corresponding modular operators are geometric. Analogous results hold
when $I$ is replaced by an arbitrary interval.
\vskip .05in
\bf \noindent Remark. \rm Locality leads immediately to the canonical
so--called ``Jones--Wassermann'' inclusion
$\pi_i({\cal L}_IG)^{\prime\prime} \subseteq 
\pi_i({\cal L}_{I^c}G)^\prime$ ([18], [40]). This inclusion {\it
measures the failure of Haag duality} in non--vacuum representations.
\vskip .05in 
\noindent \bf Proof. \rm By the Reeh--Schlieder theorem, the vacuum
vector is  
cyclic for $\pi_0({\cal L}_IG)^{\prime\prime}$, and hence
$\pi_0({\cal L}_IG)^{\prime}$ (since it contains $\pi_0({\cal
L}_{I^c}G)^{\prime\prime}$) . Let $e$ be the projection onto
$\overline{N\Omega}$. 
Then $N\rightarrow Ne$, $x\mapsto xe$ is an isomorphism. Clearly $Ne$
may be identified with $\pi_0({\cal L}_IG)^{\prime\prime}$. Its commutant is
$JNJe$, so $\pi_0({\cal L}_{I^c}G)^{\prime\prime}$. The identification of the
modular operators is immediate. Now $SU(1,1)=SU_+(1,1)$ acts on the
vacuum representation fixing the vacuum vector and carries $I$
onto any other interval of the circle. Since the modular operators lie
in $SU_\pm (1,1)$, the results for an arbitrary interval follow by
transport of structure. 

\vskip .1in

\vskip .1in
\noindent \bf Theorem~D (generalised Haag duality). \it Let $e$ be the
projection onto the vacuum 
subrepresentation of $\pi^{\otimes\ell}$. Then
$\pi_P({\rm Cliff}(L^2(I,W)))^{\prime\prime} \bigcap ({\Bbb
C}e)^\prime =\pi^{\otimes\ell}({\cal L}_IG)^{\prime\prime}$. 
Moreover $\pi^{\otimes\ell}({\cal L}_IG)^{\prime\prime}$ is the subalgebra of
the ``observable algebra'' $\pi^{\otimes \ell}(LG)^{\prime\prime}$
commuting with all fields $\pi_P(a(f))$ with $f$ localised in $I^c$.
\vskip .05in
\noindent \bf Proof. \rm The first assertion is just the second of the
Jones relations $N=\{x\in M:ex=xe\}$ and therefore a consequence of
Takesaki devissage. To prove the second, note that
$$\pi^{\otimes \ell}({\cal L}_IG)^{\prime\prime}\subseteq \pi_P({\rm
Cliff}(L^2(I,W)))^{\prime\prime}\bigcap \pi^{\otimes\ell} 
(LG)^{\prime\prime} \subseteq \pi_P({\rm
Cliff}(L^2(I,W)))^{\prime\prime}\bigcap ({\Bbb C}e)^\prime
=\pi^{\otimes\ell}({\cal L}_IG)^{\prime\prime}.$$  
Thus we obtain $\pi^{\otimes \ell}({\cal L}_IG)^{\prime\prime} =  \pi_P({\rm
Cliff}(L^2(I,W)))^{\prime\prime}\bigcap \pi^{\otimes\ell} 
(LG)^{\prime\prime}$. But $\pi_P({\rm
Cliff}(L^2(I,W)))^{\prime\prime}$ is equal to the 
graded commutant of $\pi_P({\rm Cliff}(L^2(I^c,W)))$. Since all
operators in $\pi^{\otimes\ell}(LG)^{\prime\prime}$ are even, it
follows that 
$\pi_P({\rm Cliff}(L^2(I^c,W)))^{\prime}\bigcap \pi^{\otimes\ell}
(LG)^{\prime\prime} =\pi^{\otimes \ell}({\cal L}_IG)^{\prime\prime}$, as required.

\vskip .1in
\noindent \bf Theorem~E (von Neumann density). \it Let $I_1$ and $I_2$
be {\it touching} 
intervals obtained by removing a point from the proper interval $I$.
Then if $\pi$ is a positive energy representation of $LG$ (not
necessarily irreducible), we have
$\pi({\cal L}_{I_1}G)^{\prime\prime} \vee \pi({\cal L}_{I_2}G)^{\prime\prime} =
\pi({\cal L}_IG) ^{\prime\prime}$ (``irrelevance of points''). 
\vskip .05in
\noindent \bf Proof. \rm By local equivalence, there is an isomorphism
$\pi$ between $\pi_0({\cal L}_{I}G)^{\prime\prime}$ and
$\pi({\cal L}_{I}G)^{\prime\prime}$ taking $\pi_0(g)$ onto $\pi(g)$ for
$g\in {\cal L}_IG$. Thus $\pi$ carries 
$\pi_0({\cal L}_{I_1}G)^{\prime\prime} \vee \pi_0({\cal L}_{I_2}G)^{\prime\prime}$
onto $\pi({\cal L}_{I_1}G)^{\prime\prime} \vee \pi({\cal L}_{I_2}G)^{\prime\prime}$.
It therefore suffices to prove the result for the vacuum
representation $\pi_0$. Let $J_1=I_1^c$ and $J_2=I_2^c$. Now for
$k=1,2$ we have
$\pi^{\otimes\ell}({\cal L}_{J_k}G)^{\prime\prime}=\pi_P({\rm
Cliff}(L^2(I_k,W)))^\prime \cap ({\Bbb C}e)^\prime$. So
$$\eqalign{\pi^{\otimes\ell} ({\cal L}_{J_1}G)^{\prime\prime}\cap \pi^{\otimes\ell}
({\cal L}_{J_2}G)^{\prime \prime} &=\pi_P({\rm Cliff}(L^2(I_1,W)))^\prime\cap
\pi_P({\rm Cliff}(L^2(I_2,W)))^\prime \cap ({\Bbb C}e)^\prime\cr
&=\pi_P({\rm Cliff}(L^2(I,W)))^\prime \cap ({\Bbb C}e)^\prime
=\pi^{\otimes\ell}({\cal L}_{I^c}G)^{\prime\prime}.\cr}$$
Here we have used Theorem~C and the equality
$L^2(I,W)=L^2(I_1,W)\oplus L^2(I_2,W)$. Taking
commutants, we get 
$\pi^{\otimes\ell} ({\cal L}_{J_1}G)^{\prime}\vee \pi^{\otimes\ell}
({\cal L}_{J_2}G)^{\prime} =\pi^{\otimes \ell}({\cal L}_{I^c}G)^\prime$.
Compressing by $e$, this yields
$\pi_0({\cal L}_{J_1}G)^{\prime}\vee \pi_0
({\cal L}_{J_2}G)^{\prime} =\pi_0({\cal L}_{I^c}G)^\prime$.
Using Haag duality in the vacuum representation to identify these
commutants, we get \break
$\pi_0({\cal L}_{I_1}G)^{\prime\prime}
\vee\pi_0({\cal L}_{I_2}G)^{\prime\prime}=\pi_0({\cal L}_IG)^{\prime\prime}$, 
as required.
\vskip .1in
\noindent \bf Theorem~F (irreducibility). \it Let $A$ be finite subset
of $S^1$ and let $L^AG$ be the subgroup of $LG$ consisting of loops
trivial to all orders at points of $A$. Let ${\cal L}^AG$ be the
corresponding subgroup of ${\cal L}G$. If $\pi$ is a positive energy
representation of $LG$ (not necessarily irreducible), we have
$\pi(L^AG)^{\prime\prime}=\pi(LG)^{\prime\prime}$. In particular the
irreducible positive energy representations of $LG$ stay irreducible 
and inequivalent when restricted to $L^AG$. 
\vskip .05in
\noindent \bf Proof. \rm Clearly ${\cal L}^AG={\cal L}_{I_1}G \cdot
\cdots \cdot {\cal L}_{I_n}G$, if $S^1\backslash A$ is the disjoint
union of the consecutive intervals $I_1$,\dots , $I_n$. Let
$J_k$ be the interval obtained by adding the common endpoint to
$I_k\cup I_{k+1}$ (we set
$I_{n+1}=I_1$). By von Neumann density, 
$\pi({\cal L}_{I_k}G)^{\prime\prime}\vee \pi({\cal
L}_{I_{k+1}}G)^{\prime\prime}=\pi({\cal L}_{J_k}G)^{\prime\prime}$.
Hence $\pi({\cal L}^AG)^{\prime\prime} =\bigvee \pi({\cal
L}_{J_k}G)^{\prime\prime}$. But the subgroups ${\cal L}_{J_k}G$ generate
${\cal L}G$ algebraically. Hence $\pi({\cal L}^AG)^{\prime\prime}
=\pi({\cal L}G)^{\prime\prime}$. Taking commutants, we get
$\pi({\cal L}^AG)^{\prime}=\pi({\cal L}G)^{\prime}$. By Schur's lemma,
this implies that the irreducible positive energy
representations of $LG$ stay irreducible 
and inequivalent when restricted to $L^AG$. 
\vskip .1in
\noindent \bf Remark. \rm Direct proofs of Haag duality (Theorem~C) have been
discovered since the announcement in [18] that do not use Takesaki
devissage from fermions. Theorems A, B and F can also be proved without using
Takesaki devissage. In fact Jones and I proved in 
[40] that the topology on ${\cal L}G$ induced by pulling
back the strong operator topology on $U({\cal F}_P)$ makes ${\cal L}^AG$
dense in ${\cal L}G$. Since any level $\ell$ representation $\pi$ is
continuous for this 
topology, it follows that $\pi({\cal L}^AG)$ is dense in $\pi({\cal
L}G)$ in the strong operator topology. So $\pi({\cal
L}^AG)^{\prime\prime}= \pi({\cal L}G)^{\prime\prime}$ and Theorem~F
follows. The reader is warned that several incorrect proofs of these
results have appeared in published articles.
\vfill\eject

\noindent \bf CHAPTER III.~THE BASIC ORDINARY DIFFERENTIAL EQUATION. 
\vskip .2in
\noindent \bf 18. The basic ODE and the transport problem. \rm 
Consider the ODE
$${df\over dz} = {Pf\over z} +{Qf\over 1-z}\eqno{(1)}$$
where $f(z)$ takes values in $V={\Bbb C}^N$ and $P,Q\in {\rm End}\, V$.
Suppose that $P$ has distinct
eigenvalues $\lambda_i$ with corresponding eigenvectors $\xi_i$, none
of which differ by positive integers, 
and $Q$ is a non--zero multiple of a rank one idempotent in general
position with respect to $P$. Thus $Q^2=\delta Q$, ${\rm
Tr}(Q)=\delta$ with $\delta\ne 0$, so that $Q(x)=\phi(x)v$ for
$v\in V$, $\phi\in V^*$ with $\phi(v)=\delta$. ``General position'' means
that
$v=\sum \delta_i\xi_i$ with $\delta_i\ne 0$ for all $i$ and
$\phi(\xi_i)\ne 0$ for all $i$; the eigenvectors can therefore be
normalised so that $\phi(\xi_i)=1$. Let
$R=Q-P$ and suppose that $R$ satisfies the same conditions as $P$ with 
respect to $Q$. Let $(\zeta_j,-\mu_j)$ be the normalised eigenvectors and
eigenvalues of $R$. Let $f_i(z)=\sum \xi_{i,n}
z^{\lambda_i+n}$ be the formal power series solutions of (1) expanded
about $0$ with $\xi_{i,0}=\xi_i$.  The $f_i(z)$'s are defined and
converge in $\{z: |z|<1, z\notin [0,1)\}$. If $g(z)=f(z^{-1})$, then 
$${dg\over dz} ={ Rg\over z} +{ Qg\over 1-z},\eqno{(2)}$$ 
so we can look for formal power series solutions
$h_j(z)=\sum \zeta_{j,n} z^{\mu_j-n}$ of (1) expanded about $\infty$
with $\zeta_{j,0}=\zeta_{j}$. The $h_j(z)$'s are defined and converge
in $\{z: |z|>1, z\notin [1,\infty)\}$. The solutions $f_i(z)$ and
$h_j(z)$ extend analytically to single--valued holomorphic functions
on ${\Bbb C}\backslash [0,\infty)$. 
\vskip .1in
\noindent \bf Problem. \it Compute the transport coefficients
$c_{ij}$ for which $f_i(z) = \sum c_{ij} h_j(z)$ for $z\in {\Bbb
C}\backslash [0,\infty)$.
\vskip .1in
\noindent \rm This problem will be solved by finding a
rational canonical form for the matrices $P,Q,R$ which links the ODE
with the generalised hypergeometric equation, first studied by Thomae.
It can be seen directly that the projected solutions
$(1-z)\phi(f_i(z))$ can be represented by multiple Euler integrals.
This allows one coefficient of the transport matrix $(c_{ij})$ to be
computed when the $\lambda_i$'s and $\mu_j$'s are real and $\delta$ is
negative. The rational canonical form shows that the transport 
matrices are holomorphic functions of the $\lambda_i$'s and $\mu_j$'s
alone, symmetric in an obvious sense. So the computation of the
$c_{ij}$'s follows by analytic continuation and symmetry from the
particular solution:
\vskip .1in
\noindent \bf Theorem. \it The coefficients of the transport matrix
are given by the formula
$$c_{ij} =e^{i\pi(\lambda_i-\mu_j)}{\prod_{k\ne i} 
\Gamma(\lambda_i -\lambda_k+1)  \,
\prod_{\ell\ne j} \Gamma( \mu_j - \mu_\ell) \over
\prod_{\ell\ne j} \Gamma( \lambda_i -\mu_\ell  +1)\,
\prod_{k\ne i} \Gamma(\mu_j - \lambda_k)}.$$
\vskip .1in

\rm \noindent For applications it will be convenient to have a
slightly 
generalised version of this result. Let $B$ be a
matrix of the form $-\alpha I + \beta Q$ ($\beta\ne 0$) where $Q$ is a 
rank one idempotent. Let $A$ be a matrix such that both $A$ and $B-A$
are in 
general position with respect to $Q$ and have distinct eigenvalues not
differing by integers (so distinct). Around $0$ the ODE
$${df\over dz}= {Af\over z} + {Bf\over 1-z}\eqno{(3)}$$
has a canonical basis of solutions $f_i(z)=\xi_i z^{\lambda_i} +
\xi_{i,1} z^{\lambda_i +1} + \cdots$, where $A\xi_i=\lambda_i\xi_i$ and
$\phi(\xi_i)=1$ if $Q(\xi)=\phi(\xi)v$. Similarly around $\infty$, the
ODE has a canonical basis of solutions $h_j(z) =\zeta_i z^{\mu_i}
+\zeta_{i,1} z^{\mu_i-1}+\cdots$ where $(A-B)\zeta_i=\mu_i \zeta_i$
and $\phi(\zeta_i)=1$. 
\vskip .1in
\noindent \bf Corollary. \it In ${\Bbb C}\backslash [0,\infty)$ we have
$f_i(z) = \sum c_{ij} h_j(z)$, where
$$c_{ij} =e^{i\pi(\lambda_i-\mu_j)}{\prod_{k\ne i} 
\Gamma(\lambda_i -\lambda_k+1)  \,
\prod_{\ell\ne j} \Gamma( \mu_j - \mu_\ell) \over
\prod_{\ell\ne j} \Gamma( \lambda_i -\mu_\ell +\alpha +1)\,
\prod_{k\ne i} \Gamma(\mu_j - \lambda_k-\alpha)}.$$
\vskip .05in
\noindent \bf Proof. \rm By a gauge transformation $f(z)\mapsto
(1-z)^\gamma f(z)$, the ODE $(3)$ is changed into the ODE considered
before. It is then trivial to check that the transport
relation for that ODE implies the stated transport relation for $(3)$.
\vskip .1in

\noindent \bf 19. Analytic transformation of the ODE (cf [16]). \rm
Consider   
the ODE $f^\prime(z) = A(t,z)f(z)$ where $A(t,z) = \sum_{n\ge
0}  A_n(t) z^{n-1}$ with each matrix $A_n(t)\in {\rm End}\, V$ a
polynomial (or holomorphic function) in $t\in W=
{\Bbb C}^m$ and $A(t,z)$ is convergent in $0<|z|<R$ for all $t\in
{\Bbb C}^m$. 
\vskip .1in

\noindent \bf Proposition. \it Let $U=\{t\in {\Bbb C}^m : 
\hbox{$A_0(t)$ has no eigenvalues differing by positive 
integers}\}$. For $t\in U$, there is a unique gauge transformation
$g(t,z)\in GL(V)$, holomorphic on $U\times \{z:|z|<R\}$, such that
$g(t,z)^{-1} A(t,z) 
g(t,z) - g(t,z) ^{-1} {\partial}g(t,z)/\partial z= 
A_0(t)/z$.
\vskip .05in
\noindent \bf Proof. \rm If we write $g(t,z) =\sum_{n\ge 0}
g_n(t) z^n$ with $g_0(t)=I$, then the
$g_n(t)$'s are given by the recurrence relation
$$ng_n(t) =n \,(n-{\rm ad}\,A_0(t))^{-1}
\sum _{m=1}^n A_m(t) g_{n-m}(t).$$
Let $\overline {B}$ be a closed ball in $U$. Then 
$\sup_n \|n(n-{\rm ad} A_0(t))^{-1}\|$ is bounded by $M<\infty$ on
$\overline {B}$. So $\|g_n(t)\|$ is bounded on $\overline{B}$
by the solutions $f_n$ of the recurrence relation
$$nf_n = \sum_{m=1}^n b_m f_{n-m},$$
where $b_m= M\sup_{t\in\overline{B}} \|A_m(t)\|$ and $\sum_{m\ge 1}
b_m z^m$ is convergent in $|z|<R$. 
But then $f(z)=\sum_{n\ge 0} f_n z^n$ is the formal power series
solution of $zf^\prime(z) =(\sum_{m\ge 1} b_m z^m) f(z)$ with
$f(0)=1$, i.e. 
$f^\prime(z) =b(z) f(z)$ where $b(z)=\sum_{m\ge 0} b_{m+1} z^m$. This
has the unique solution $f(z)=\exp \int_0^z b(w)\, dw$ so that in
particular $f(z) =\sum f_n z^n$ is convergent in $|z|<R$. Since
$\|g_n(t)\|\le f_n$, it follows that $\sum g_n(t)
z^n$ converges uniformly on $\{(t,z): t\in \overline{B},\,|z|\le r\}$
for any $r<R$. Since 
$t\mapsto g_n(t)$ is holomorphic in $t$, for fixed $z$, $g(z,t)$ is
the uniform limit on compacta of holomorphic functions in $t$. Since
the uniform limit on compacta of holomorphic functions 
is holomorphic, it follows that $t\mapsto
g(t,z)$ is holomorphic on $U$ for fixed $z$. 

To show that $g(t,z)$ is invertible for fixed $t$, note that
${\partial_z}g = Ag-gA_0/z$. Replacing $A$ by $-A^t$, we find $f$
such that ${\partial_z}f=-fA+A_0f/z$. Hence ${\partial_z}(fg)=
[A_0,fg]/z$. The only formal power series solution $h$ of this
equation with $h(0)=I$ is $h\equiv I$. Hence $fg\equiv I$ as
required. 
\vskip .1in
\bf \noindent Remarks. \rm This argument applies also when $A_0(t)=0$.
Clearly we may apply the proposition to the basic ODE. The argument with
$A_0(t)=0$ near points 
$z\ne 0,1$ shows that the gauge transformation $g(z)$ extends to a
holomorphic map ${\Bbb C}\backslash 
[1,\infty)\rightarrow GL(N,{\Bbb C})$ such that 
$g(z)^{-1} A(z)g(z) - g(z) ^{-1} g^\prime(z) = 
A_0/z$ for $z\notin [1,\infty)$. The gauge transformation reduces the
basic ODE about $0$ to the ODE $f^\prime(z) = 
z^{-1} A_0 f(z)$ which has solutions $z^{A_0}v=\exp(A_0\log z)v$
defined in ${\Bbb C}\backslash [0,\infty)$ say. Applying the gauge
transformation, it follows that any formal power 
series solution of the original ODE is automatically convergent in
$|z|<1$ and extends to 
a single--valued holomorphic function on ${\Bbb C}\backslash [0,\infty)$.

\vskip .1in

\noindent \bf 20. Algebraic transformation of the ODE. \rm Let $P$ be a 
matrix with distinct eigenvalues $\lambda_i$ and corresponding
eigenvectors $v_i$. Let $Q$ be proportional to a rank one idempotent
on $V$ so that 
$Q(x) =\phi(x) v$ with $\phi\in V^*, v\in V$ and $\phi(v)=\delta\ne 0$.
We assume that 
$P$ is in general position with respect to $Q$. This means that the
eigenvectors $\xi_i$ satisfy $\phi(v_i)\ne 0$ and that $v=\sum \alpha_i
\xi_i$ with $\alpha_i\ne 0$ for all $i$. The next result gives a
rational canonical form for the matrices $P$, 
$Q$ and $R$. 
\vskip .1in
\noindent \bf Proposition (Rational Canonical Form). \it  If $P$ has
distinct eigenvalues and
$Q$ is a non--zero multiple of a rank one idempotent in general
position with respect to $P$, there is a 
(non--orthonormal!) basis of 
$V$ such that  
$$P=\pmatrix{ 0 & 1 &0  & & & 0\cr
              0 & 0 & 1 & & & 0\cr
              0 & 0 & 0 & 1 && 0\cr
              \cdot & \cdot & \cdot & \cdot & \cdot & \cdot\cr
            \cdot & \cdot & \cdot & \cdot & \cdot & \cdot\cr
              0 & &&&& 1\cr
              a_1& a_2& &&& a_N\cr}, 
       Q=\pmatrix{ 0 & 0&0  & & & 0\cr
               &  &  & & & \cr
               &  &  &  && \cr
              \cdot & \cdot & \cdot & \cdot & \cdot & \cdot\cr
            \cdot & \cdot & \cdot & \cdot & \cdot & \cdot\cr
               & &&&& \cr
              b_1& b_2& &&& b_N\cr},
-R=P-Q=\pmatrix{ 0 & 1 &0  & & & 0\cr
              0 & 0 & 1 & & & 0\cr
              0 & 0 & 0 & 1 && 0\cr
              \cdot & \cdot & \cdot & \cdot & \cdot & \cdot\cr
            \cdot & \cdot & \cdot & \cdot & \cdot & \cdot\cr
              0 & &&&& 1\cr
              c_1& c_2& &&& c_N\cr}, $$
where $b_N={\rm Tr}(Q)\ne 0$ and $c_i=a_i-b_i$. Conversely if $P$ and
$Q$ are of the above form and the 
roots of $a(t)=t^N-\sum a_it^{i-1}$ (the characteristic polynomial of
$P$) are distinct, then $P$ and $Q$ are in general
position iff $b(t) =\sum b_i t^{i-1}$ and $a(t)$ have no common roots
iff $c(t)=a(t)-b(t)$ and $a(t)$ have no common roots.
(Here $c(t)$ is the characteristic polynomial of $P-Q$.) 
\vskip .05in
\noindent \bf Remark. \rm This gives a unique canonical form for
$P,Q,R=Q-P$ with equivalence given by conjugation by matrices in
$GL(N,{\Bbb C})$: for $a(t)$ and $c(t)$ are the characteristic
polynomials of $P$ and $P-Q$, so that
the constants $a_i,b_j$ are invariants (since $b(t)=a(t)-c(t)$). 
Moreover the orbit space
of the pairs $(P,R)$ under the action by conjugation of $GL(N,{\Bbb
C})$ can naturally be identified with the space of rational
canonical forms.

\vskip .05in
\noindent \bf Proof. \rm Let $Q(x)=\phi(x)v$, with
$\phi(v)\ne 0$. Since $Q$ and $P$ are in general position, the elements
$\phi, \phi\circ P, \cdots,\phi\circ P^{N-1}$ form a basis of $V^*$.
In particular there is a unique solution $w$ of
$\phi(w)=\phi(Pw)=\cdots = \phi(P^{N-2}w)=0$, $\phi(P^{N-1}w)=1$. The
set $w, Pw,\dots , P^{N-1}w$ must be linearly independent, because
otherwise $P^{N-1}w$ would have to be a linear combination of $w, Pw,
\dots, P^{N-2}w$ contradicting $\phi(P^{N-1}w)=1$. Thus $(P^jw)$ is a
basis of $V$. Clearly $P$ and $Q$ have the stated form with respect
to this basis. Furthermore $b_N={\rm Tr}(Q)$. 
\vskip .05in
We next must check that if $P$ and $Q$ have the stated form, then
no eigenvector $u\ne 0$ of $P$ can satisfy $Qu=0$ and no eigenvector
$\psi$ of $P^t$ can satisfy $Q^t\psi=0$. For $\psi$, the condition
$Q^t\psi=0$ means that $\psi=(x_1,x_2,\dots,x_{N-1},0)$ with $x_i\in
{\Bbb C}$. The condition $P^t\psi=\lambda\psi$ forces $x_1=\lambda
x_2$, $x_2=\lambda x_3$, \dots, $x_{N-1}=0$. Hence $x_i=0$ for all $i$
and $\psi=0$. Now suppose that $Pu=\lambda u$ and $Qu=0$. Then it is
easily verified that $u$ is proportional to
$(1,\lambda,\lambda^2,\dots,\lambda^{N-1})^t$. Thus $Qu=
(0,0,\dots,0,b(\lambda))^t$, so that $Qu\ne 0$ iff $b(\lambda)\ne 0$. 
Finally the characteristic polynomial of $R$ is $c(t)=a(t) -b(t)$.
Clearly $a(t)$ and $b(t)$ have no common roots iff $c(t)$ and $b(t)$
have no common roots, so the last assertion follows.

\vskip .1in 

\noindent \bf 21. Symmetry and analyticity properties of transport
matrices. \rm 
\vskip .1in
\bf \noindent Proposition. \it The transport matrix $c_{ij}$ from $0$ to
$\infty$ of the basic ODE depends only on the eigenvalues
$\lambda_i$ of $P$ and 
$\mu_j$ of $P-Q$. This dependence is holomorphic. Moreover the
coefficients $c_{ij}$, indexed by the eigenvalues $\lambda_i$ and
$\mu_j$, have the symmetry property
$c_{ij}(\lambda_1,\dots,\lambda_N,\mu_1,\dots,\mu_N) =
c_{\sigma i,\tau j} ( \lambda_{\sigma 1},\dots,
\lambda_{\sigma N},\mu_{\tau 1},\dots,\mu_{\tau N})$
for
$\sigma,\tau\in S_N$.
\vskip .05in
\noindent \bf Proof. \rm We can conjugate by a matrix in $GL(N,{\Bbb
C})$ so that $P$, $Q$ and $R$ are in rational canonical form. The
transport matrix from $0$ to $\infty$ is invariantly defined, so does
not change under such a conjugation. Thus the assertions are invariant
under conjugation, so it suffices to prove them when $P,Q,R$ are in
rational canonical form. Setting $g(z)=f(z/(z-1))$, where $f(z)$ is a
solution of the basic ODE, we get the ODE
$${dg\over dz} = {Pg\over z} + {Rg\over
z-1}\eqno{(4)}$$
where $R=Q-P$. Thus we have to compute the transport matrices for (4)
from $0$ to $1$ where the solutions at $0$ are labelled by the
eigenvalues $\lambda_i$ of $P$ and at $1$ by the eigenvalues of
$\mu_j$ of $-R$. We shall consider variations of 
$P$, $Q$, and $R$ within rational canonical form. $P$ and $R$ can be
specified by 
prescribing the eigenvalues $(\lambda_i)$ of $P$ and $(\mu_j)$ of
$-R$. This completely determines the $a_i$'s and $c_i$'s and hence the
$b_i$'s. The $\lambda_i$'s and $\mu_j$'s should be distinct and no
two $\lambda_i$'s or $\mu_j$'s should differ by a 
positive integer. We also impose the linear constraint 
that $\sum \lambda_i-\mu_i \ne 0$. Thus we obtain an open
path--connected subset $U_0$ of 
the $2N$--dimensional linear space $W=\{(\lambda,\mu)\}={\Bbb C}^{2N}$.
Applying the proposition in section~19 with
$t=(\lambda,\mu)\in W$ and
$A(t,z)=z^{-1} P + (z-1)^{-1}R$, we deduce that the
gauge transformations $g(t,z)$, $h(t,z)$ transforming $A(t,z)$ into
$z^{-1} P$ and $(z-1)^{-1} R$ respectively
depend holomorphically on $t\in U$ for a fixed $z\in (0,1)$. We
already saw in section~20 that
the normalised eigenvectors of $P$ and $R$ are given by
$$\xi_i(t)=b(\lambda_i)^{-1} (1,\lambda_i,\lambda_i^2, \cdots,
\lambda_i^{N-1})^t \quad
\zeta_j(t)=b(\mu_j)^{-1} (1,\mu_j,\mu_j^2,\cdots,\mu_j^{N-1})^t.$$
Thus the normalised solutions at $0$ are
$z^{\lambda_i}g(t,z)\xi_i(t)$ and the normalised solutions at
$1$ are given by $(z-1)^{\mu_j} h(t,z)\zeta_j(t)$. So the
transport matrix $c_{ij}(t)$ (independent of $z$) is specified by the equation
$$z^{\lambda_i}g(t,z)\xi_i(t) =\sum c_{ij}(t)
(z-1)^{-\mu_j} h(t,z)\zeta_j(t)$$
for $|z-1/2|<1/2$. Fix such a value of $z$ (say $z=1/2$) and let
$(\psi_j(t))$ be the dual basis to $(\zeta_j(t))$. Clearly $\psi_j(t)$
is a rational function of $(\lambda,\mu)$ so is holomorphic on $U$.
Moreover
$$c_{ij}(t)
=(z-1)^{\mu_j}z^{\lambda_i}
(\psi_j(t),h(t,z)^{-1}g(t,z)\xi_i(t)).$$
This equation shows that $c_{ij}(t)$ depends holomorphically on $t\in
U_0$ and has the stated symmetry properties.
\vskip .1in

\noindent \bf 22. Projected power series solutions.
\rm Let
$\lambda=\lambda_i$ be an eigenvalue of $P$ and consider the
corresponding (formal)
power series solution $f_i(z)=\sum \xi_{i,n}z^{\lambda_i
+n}$ of the basic ODE. Dropping the index $i$ for clarity, we have
$$zf^\prime(z) =  P f +  Q( z + z^2 + z^3 +
\cdots) f,$$
with $f(z) = \sum \xi_n z^{\lambda +n}$ and $P\xi_0
=\lambda \xi_0$. Substituting in the formal power series and
dividing out by $z^{\lambda}$, we get
$$\sum_{n\ge 0} (n+\lambda) \xi_n z^n
=\sum_{n\ge 0}  P\xi_n z^n + Q(z + z^2 + z^3
+\cdots) \sum_{n\ge 0} \xi_n z^n.$$
Thus for $n\ge 1$ we get
$$(n+ \lambda - P) \xi_n = Q(\xi_0 +
\cdots +\xi_{n-1})$$
and hence
$$Q\xi_n = Q(n+ \lambda -P)^{-1} Q
(\xi_0+\cdots +\xi_{n-1}).$$
Let $Q(\xi_0+\cdots +\xi_n) =\alpha_n v$, where $\alpha_n\in {\Bbb
C}$. Thus we obtain the recurrence relation
$\alpha_n -\alpha_{n-1} = \chi(\lambda + n) \alpha_n$, so
that $\alpha_n =\chi_P(\lambda +n)\alpha_{n-1}$, where the rational
function 
$\chi_P(t)$ is defined by $Q + Q(tI-P)^{-1}Q=\chi_P(t)Q$. Thus,
reintroducing the index $i$, we have
$$\alpha_{i,n}=\alpha_{i,0} \prod_{m=1}^n \chi_P(\lambda_i
+m),\eqno{(5)}$$
where $\alpha_{i,0}=\phi(\xi_i)$. We now must compute $\chi_P(t)$.
Bearing in mind that equation (2) gives the corresponding power series
expansions about $\infty$, we define $\chi_R(t)$ by
$Q+Q(tI-R)^{-1}Q=\chi_R(t)Q$. 
\vskip .1in
\noindent \bf Inversion lemma. \it $\chi_R(t)=\chi_P(-t)^{-1}$.
\vskip .05in
\noindent \bf Proof. \rm Let $A$ be an invertible matrix with
$QA^{-1}Q = (1-\alpha)Q$, where $\alpha\ne 0$. 
Expanding $(A-Q)^{-1}=(I-A^{-1}Q)^{-1} A^{-1}$, we
find that $Q(A-Q)^{-1} Q= (\alpha^{-1}-1)Q$. Hence
$$\chi_R(t) Q = Q + Q(t-R)^{-1} Q = Q+Q(t+P-Q)^{-1}Q
=\alpha^{-1}Q,$$
if $Q(t+P)^{-1}Q=(1-\alpha)Q$. But
$Q(t+P)^{-1}Q=-Q(-t-P)^{-1}Q=(1-\chi_P(-t))Q$, so that
$\alpha=\chi_P(-t)$ and hence $\chi_R(t)=\alpha^{-1}=\chi_P(-t)^{-1}$
as required.
\vskip .1in
\bf \noindent Corollary. \it $\chi_P(t) = \prod (t-\mu_i)/\prod
(t-\lambda_j)$ 
where the $\mu_j$'s are the eigenvalues of $P-Q$.
\vskip .05in
\noindent \bf Proof. \rm $X_P(t)$ has the form $p(t)/\prod
(t-\lambda_i)$, where $p(t)$ is a monic polynomial of degree $N$.
Similarly $X_R(t)$ has the form $q(t)/\prod(t+\mu_i)$ where the
$\mu_i$'s are the eigenvalues of $-R=P-Q$. Since
$X_R(t)=X_P(-t)^{-1}$, we see that $p(t)=\prod(t-\mu_i)$ and
$q(t)=\prod(t+\lambda_i)$, as required.

\vskip .1in
\noindent \bf Corollary. \it $\sum \lambda_i - \sum \mu_i =\delta$.
\vskip .05in
\noindent \bf Proof. \rm This follows by taking the trace of the
identity $P+R=Q$.
\vskip .1in
From (5) and the formula for $\chi_P(t)$, we have for $n\ge 1$ 
$$\alpha_{i,n} = \alpha_{i,0}\prod_{j=1}^N \prod_{m=1}^n {m +
\lambda_i -\mu_j \over
m + \lambda_i -\lambda_j} ,$$
where $\alpha_{i,0} = \phi(\xi_i)$. 
\vskip .1in

\noindent \bf 23. Euler--Thomae integral representation of projected
solutions (cf [36],[45]).
\rm We assume here that the eigenvalues $\lambda_i$ of $P$ are real
with $\lambda_1 > 
\lambda_2 > \cdots >\lambda_N$; that the eigenvalues $\mu_i$ of $P-Q$
are real with $\mu_1 > \mu_2 >\cdots >\mu_N$; and that
$\lambda_1 + 1 >\mu_j>\lambda_1$ for all $j$. In particular this implies that
$\delta={\rm Tr}(Q)$ must be {\it negative}. We start by obtaining
an integral  
representation of the projected solutions $(1-z)\phi(f_i(z))$ around $0$.
Recalling that the eigenvectors $\xi_i$ and $\zeta_i$ of $P$ and $P-Q$
are normalised so that $\phi(\xi_i)=1=\phi(\eta_i)$, where
$Q(x)=\phi(x)v=\phi(x)\eta$, we have already shown that
$$(1-z)^{-1}z^{-\lambda_i} \phi(f_i(z))=\sum_{n\ge
0} \alpha_{i,n} z^n=\sum_{n\ge
0}z^n\cdot \prod_{j=1}^N \prod_{m=1}^n {m+\lambda_i
-\mu_j\over m+\lambda_i-
\lambda_j}.$$
Using the formula
$(a)_n\equiv a(a+1)\cdots (a+n-1) =\Gamma(a+n)/\Gamma(a)$, we get
$$(1-z)^{-1}z^{-\lambda_i} \phi(f_1(z))= \sum_{n\ge
0} {( \lambda_1-\mu_1 +1)_n\over n!}
\prod_{j\ne 1} {\Gamma(
\lambda_1-\mu_j+n+1)\Gamma(\lambda_1-
\lambda_j+1) 
\over \Gamma( \lambda_1 -
\mu_j+1)\Gamma(\lambda_1 - \lambda_j+n+1)}.$$
Using the beta function identity
$\Gamma(a)\Gamma(b)/\Gamma(a+b) =\int_0^1 t^{a-1} (1-t)^{b-1}\, dt$
for $a,b>0$, we obtain
$$\phi(f_1(z))=(1-z)z^{\lambda_1} K \int_0^1 \int_0^1
\cdots \int_0^1 (1-z t_2  \cdots t_N)^{
\mu_1- \lambda_1-1} \prod_{j\ne 1} t_j^{
\lambda_1 -\mu_j} (1-t_j)^{ \mu_j -
\lambda_j -1} \, dt_j,\eqno{(6)}$$
where
$$K=\prod_{j\ne 1} {\Gamma(
\lambda_1- \lambda_j+1)\over \Gamma(\lambda_1
- \mu_j+1) \Gamma(\mu_j -
\lambda_j)}.$$
(The inequalities $\mu_i> \lambda_i$ and
$\lambda_1-\mu_j>-1$ guarantee that this summation by
integrals is valid.) Note that this Euler type integral representation
is also valid for 
$z$ real and negative, since it 
is analytic in $z$ where defined. The solutions about
$\infty$ have a Laurent expansion (for $|z|$ large)
$g_j(z) =\zeta_j z^{\mu_j} + \zeta_{j,1} z^{\mu_j-1} +\cdots$
where $\zeta_j$ are the eigenvectors of $P-Q$ with
$(P-Q)\zeta_j=\mu_j\zeta_j$. Hence the projected solution
$\phi(g_j(z))$ satisfies $\phi(g_j(z))\sim (\zeta_j,\eta)
z^{\mu_j}$ because of the normalisation $\phi(\zeta_j)=1$. In
particular if $x$ is large and negative 
$\phi(g_j(x))\sim |x|^{\mu_j} e^{ \pi i\mu_j}$.
Let $c_{ij}$ be the transport matrix connecting the solutions at $0$
and $\infty$, so that
$f_1(z) = \sum c_{1j} g_j(z)$.
Since $Q$ and $P$ are in general position, we lose no information by
writing the above equation as
$\phi(f_1(z)) =\sum c_{1j} \phi(g_j(z))$. 
Since $\mu_1$ is the largest of the $\mu_j$'s, we find that for $x$
large and negative,  
$$\phi(f_1(x)) \sim c_{11}
|x|^{\mu_1} e^{i\pi \mu_1}.\eqno{(7)}$$
On the other hand by (6) we have for $x<<0$
$$\phi(f_1(x))\sim K e^{i\pi\lambda_1} |x|^{\mu_1}
\prod_{j\ne i} \int_0^1 t_j^{\mu_1-
\mu_j-1}(1-t_j)^{ \mu_j- \lambda_j -1}\,
dt_j.\eqno{(8)}$$ 
Comparing (7) and (8), we obtain
$$c_{11} = e^{ i\pi(\lambda_1-\mu_1)} K
\prod_{j\ne 1} \int_0^1 t_j^{\mu_1-
\mu_j-1}(1-t_j)^{ \mu_j- \lambda_j -1}\,
dt_j=K e^{i\pi (\lambda_1-\mu_1)}\prod_{j\ne 1} {\Gamma( \mu_1-
\mu_j) \Gamma(\mu_j 
-\lambda_j) 
\over\Gamma(\mu_1 -\lambda_j)}.$$
Substituting in the value of $K$, we get the fundamental formula:
$$c_{11} = e^{i\pi(\lambda_1-\mu_1)} \prod_{j\ne 1} 
{\Gamma(\lambda_1 - \lambda_j+1)
\Gamma( \mu_1 - \mu_j) \over
\Gamma(\lambda_1 - \mu_j +1) \Gamma(
\mu_1 - \lambda_j)}.\eqno{(9)}$$

\vskip .1in

\noindent \bf 24. Computation of transport matrices. \rm
\vskip .1in

\noindent \bf Theorem. \it The transport matrix $c_{ij}$ from the
solutions at 
$0$ to the solutions at $\infty$ of the basic ODE is given by
$$c_{ij} = e^{i\pi(\lambda_i-\mu_j)}{\prod_{k\ne i} 
\Gamma(\lambda_i - \lambda_k+1)  \,
\prod_{\ell\ne j} \Gamma( \mu_j - \mu_\ell) \over
\prod_{\ell\ne j} \Gamma( \lambda_i - \mu_\ell +1)\,
\prod_{k\ne i} \Gamma( \mu_j - \lambda_k)}.$$
\vskip .1in
\noindent \bf Proof. \rm We obtained this formula in section~23 for
$c_{11}$ when $\lambda_i,\mu_j$ took on special values.
On the other hand $c_{11}$ and the right hand side are
analytic 
functions of $\lambda_i,\mu_j$. The special values sweep
out an open subset of the real part of the parameter space $U_0$, so by
analytic continuation we must have equality for all parameters in $U_0$.
The formula for $c_{ij}$ now follows immediately from the symmetry
property of the $c_{ij}$'s.
\vfill\eject

\bf \noindent CHAPTER IV.~VECTOR AND DUAL VECTOR PRIMARY FIELDS.
\vskip .1in
\noindent \bf 25. Vector and dual vector primary fields. \rm Let $V$
be an irreducible representation of $SU(N)$. Then ${\cal
V}=C^\infty(S^1,V)$ has an action of $LG\rtimes {\rm Rot}\, S^1$ with
$LG$ acting by multiplication and ${\rm Rot}\, S^1$ by rotation,
$r_\alpha f(\theta) =f(\theta+\alpha)$. There is corresponding
infinitesimal action of $L^0\g\rtimes {\Bbb R}$ which leaves invariant
the finite energy subspace ${\cal V}^0$. We may write ${\cal V}^0=\sum
{\cal V}(n)$ where ${\cal V}(n)=z^{-n}\otimes V$. Set $v_n=z^n v$ for 
for $v\in V$. Thus $dv_n =-n v_n$ (so that $d=-i\, d/d\theta$)  and
$X_nv_m=(Xv)_{m+n}$. Let $H_i$ and $H_j$ be
irreducible positive energy representations at level $\ell$. A map
$\phi:{\cal V}^0 \otimes H_j^0 \rightarrow H_i^0$ commuting with the
action of $L^0\g\rtimes {\rm Rot}\, S^1$ is called a {\it primary field
with charge $V$}. For $v\in V$ we define
$\phi(v,n)=\phi(v_n):H^0_i\rightarrow H^0_j$: these are called the
{\it modes} of $\phi$. The intertwining property of $\phi$ is expressed
in terms of the modes through the commutation relations:
$$[X(n),\phi(v,m)]=\phi(X\cdot v,m+n),\quad [D,\phi(v,m)]=-m\phi(v,m).$$
\vskip .1in
\noindent \bf Uniqueness Theorem. \it If $\phi:{\cal
V}^0\otimes H_j^0\rightarrow H_i^0$ is a primary field, 
then $\phi$ restricts to a $G$--invariant map $\phi_0$ of
${\cal 
V}(0)\otimes H_j(0)=V\otimes H_j(0)$ into $H_i(0)$. Moreover $\phi$ is uniquely
determined by $\phi_0$, the {\rm initial term} of $\phi$. 
\vskip .05in
\noindent \rm \bf Proof. \rm ${\cal V}(0)\otimes H_j(0)$ is fixed by
${\rm Rot}\, S^1$ and hence so is its image under $\phi$. It therefore
must lie in $H_i(0)$. Since $\phi$ is $G$--equivariant (or equivalently
$\g$--equivariant), the restriction of $\phi$ is $G$--equivariant. To
prove uniqueness, we must show that if the initial term $\phi_0$
vanishes then so too does $\phi$. It clearly suffices to show that
$(\phi(\xi\otimes f),\eta)=0$ for all $\xi\in H^0_j$, $f\in {\cal
V}^0$ and $\eta\in H_i^0$. By assumption this is true for $\xi\in
H_j(0)$, $v\in {\cal V}(0)$ and $\eta\in H_i(0)$. By ${\rm Rot}\,
S^1$--invariance, this is also true if $v\in {\cal V}(n)$ for $n\ne
0$ and hence for any $v\in {\cal V}^0$.

Now we assume by induction on $n$ that
$(\phi(a_n a_{n-1} \cdots a_1\xi \otimes v),\eta)=0$
whenever $\xi\in H_j(0)$, $\eta\in H_i(0)$, $v\in {\cal V}^0$ and $a_k
=X_k(m_k)$ with $m_k<0$. Then
$$(\phi(a_{n+1} a_n \cdots a_1 \xi \otimes v),\eta) 
=- (\phi(a_n\cdots a_1 \xi\otimes a_{n+1}v),\eta) 
+(\phi(a_n \cdots a_1\xi\otimes v),a_{n+1}^*\eta),$$
and both terms vanish, the first by induction and the second because
$$a_{n+1}^*\eta =X_{n+1}(m_{n+1})^*\eta=-X_{n+1}(-m_{n+1})\eta=0.$$
Finally we prove by induction on $n$ that
$(\phi(\xi\otimes v),b_n\cdots b_1\eta)=0$
for all $\xi\in H_j^0$, $v\in {\cal V}^0$, $\eta\in H_i(0)$ and
$b_k=X_k(m_k)$ with $m_k<0$. In fact
$$(\phi(\xi\otimes v), b_{n+1}b_n\cdots b_1\eta)
=(\phi(b_{n+1}^*\xi\otimes v + \xi\otimes b_{n+1}^*v),b_n\cdots
b_1\eta),$$
which vanishes by induction.
\vskip .1in
\bf \noindent Adjoints of primary fields. \rm Let
$\phi(v,n):H_j^0\rightarrow H_i^0$ be a primary field of charge $V$. Thus
$\phi(v,n)$ takes $H_j(m)$ into $H_i(m-n)$ and satisfies
$[X(m),\phi(v,n)]=\phi(X\cdot
v,n+m)$, $[D,\phi(v,n)]=-n\phi(v,n)$.
Hence the adjoint operator $\phi(v,n)^*$ carries $H_i(m)$ into
$H_j(m+n)$. Let $\psi(v^*,n)=\phi(v,-n)^*$ where $v^*\in V^*$ is
defined using the inner product: $v^*(w)=(w,v)$. Thus
$\psi(v^*,n):H_i(m)\rightarrow H_j(m-n)$, so that $\psi(v^*,n)$ takes
$H_i^0$ into $H_j^0$. Taking adjoints in the above equation, we get
$[D,\psi(v^*,n)]=-n\psi(v^*,n)$ and $[X(m),\psi(v^*,n)]=\psi(X\cdot
v^*, n+m)$. Thus $\psi(v^*,z)$ is a primary field of charge $V^*$
called the {\it 
adjoint} of $\phi(v,z)$. Note that the initial terms of $\psi$ and
$\phi$ are related by the simple formula $\psi(v^*,0)=\phi(v,0)^*$.
Moreover for $\xi\in H_j^0$, $\eta\in H_i^0$ we have
$(\phi(v,n)\xi,\eta) =(\xi,\psi(v^*,-n)\eta)$.
\vskip .1in
\bf \noindent Fermionic initial terms. \rm Let $V=V_\square={\Bbb C}^N$ and
$W=V_\square\otimes {\Bbb C}^\ell$. The irreducible
summands of $\Lambda W=(\Lambda V)^{\otimes \ell}$ are precisely the
permissible lowest energy spaces at level $\ell$. Note that $\Lambda
W$ can naturally be identified with the lowest energy subspace of
${\cal F}_W={\cal F}_V^{\otimes \ell}$. 
\vskip .1in
\noindent \bf Lemma. \it Each non--zero intertwiner
$T\in {\rm Hom}_G(V_\square \otimes V_f,V_g)$ 
arises by taking the composition of the exterior multiplication map
$S:W\otimes \Lambda(W)\rightarrow \Lambda(W)$
with projections onto irreducible summands of the three factors, i.e.
$T=p_g S(p_\square\otimes p_f)$. 
\vskip .05in
\noindent \bf Proof. \rm Let
$e_f=e_1^{\otimes f_1 -f_2} \otimes (e_1\wedge e_2)^{\otimes f_2-f_3}
\otimes \cdots \otimes (e_1\wedge \cdots \wedge e_{N-1})^{\otimes
f_{N-1}-f_N}\otimes I^{\otimes \ell-f_1+f_N}$
be the highest weight vector for a copy of $V_f$ in $(\Lambda
V)^{\otimes\ell}$. Let $g_i=f_i$ if $i\ne k$ and $g_k=f_k+1$ so that
$g$ is a permissible signature obtained by adding one box to $f$.
Clearly the corresponding highest weight vector $e_g$ is obtained by
exterior multiplication by $e_k$ in the $f_1-f_k$ copy of $\Lambda V$
in $(\Lambda V)^{\otimes \ell}$. Let $S:W\otimes \Lambda( W)\rightarrow
\Lambda (W)$ be the map $w\otimes x\mapsto w\wedge x$. Let $p_\square$
be the projection onto the $f_1-f_k$ copy of $V$ in $W=V\otimes {\Bbb
C}^\ell$. Then, up to a sign, $S(p_\square \otimes I):V\otimes
(\Lambda V)^{\otimes \ell} 
\rightarrow (\Lambda V)^{\otimes \ell}$ is the operation of
exterior multiplication by elements of $V$ on the $f_1-f_k$ copy of
$\Lambda V$. Let $p_f$, $p_g$ be the projections onto the irreducible
modules $V_f$, $V_g$ generated by $e_f$ and $e_g$. Then
$T=p_gS(p_\square\otimes p_f):V\otimes V_f \rightarrow V_g$ satisfies
$T(e_k\otimes e_f) = \pm e_g$. Hence $T$ is non--zero. Since $S$ and
the three projections are $SU(N)$--equivariant, it follows that $T$ is
also, as required. 

\vskip .1in
\noindent \bf Construction of all vector primary fields. \it Any
$SU(N)$--intertwiner $\phi(0):V_\square \otimes H_j(0)\rightarrow
H_i(0)$ is the initial term of a vector primary field. All vector
primary fields arise as compressions of fermions so satisfy
$\|\phi(f)\|\le A\|f\|_2$ for $f\in C^\infty(S^1,V_\square)$. The map
$f\mapsto \phi(f)$ extends continuously to $L^2(S^1,V)$ and satisfies
the global covariance relation $\pi_j(g)\phi(f)\pi_i(g)^*=\phi(g\cdot
f)$ for $g\in {\cal L}G\rtimes {\rm Rot}\, S^1$.
\vskip .05in
\noindent \bf Proof. \rm By the result on initial terms, it is
possible to find an $SU(N)$--equivariant map $V\rightarrow W$,
$v\mapsto \overline{v}$ and projections $p_i$
and $p_j$ onto $SU(N)$--submodules of $\Lambda W$ isomorphic
to $V_i$ and $V_j$ such that $p_ia(\overline{v}_0) p_j:V_j\rightarrow
V_i$ is the given initial term. But $V_i$ and $V_j$ generate $LG$
modules $H_i$ and $H_j$ with corresponding projections $P_i$ and
$P_j$. The required primary field is $\phi_{ij}(v,n)=P_i
a(\overline{v}_n)P_j$ which 
clearly has all the stated properties.

\vskip .1in
\noindent \bf Dual vector primary fields. \rm Since the adjoint of a
vector primary field is a dual vector primary field, we immediately
deduce the following result.
\vskip .1in
\noindent \bf Theorem. \it Any
$SU(N)$--intertwiner $\phi(0):V_{\overline{\square}} \otimes
H_j(0)\rightarrow 
H_i(0)$ is the initial term of a dual vector primary field. All vector
dual primary fields arise as compressions of adjoints of fermions so
satisfy 
$\|\phi(f)\|\le A\|f\|_2$ for $f\in
C^\infty(S^1,V_{\overline{\square}})$. The map 
$f\mapsto \phi(f)$ extends continuously to
$L^2(S^1,V_{\overline{\square}})$ and satisfies 
the global covariance relation $\pi_j(g)\phi(f)\pi_i(g)^*=\phi(g\cdot
f)$ for $g\in {\cal L}G\rtimes {\rm Rot}\, S^1$. 

\vskip .1in
\noindent \bf 26. Transport equations for four--point functions and
braiding of primary fields. \rm We now establish the braiding
properties of primary fields. We divide the circle up into two
complementary open intervals $I$, $I^c$ 
with $I$ the upper semicircle, $I^c$ the lower semicircle say. Let $f,
g$ be test functions with $f$ supported in $I$ and $g$ in
$I^c$, so that $f\in C_c^\infty(I)$ and $g\in C^\infty_c(I^c)$. In 
general the braiding relations for primary fields will have the
following form 
$$\phi_{ik}^U(u,f)\phi_{kj}^V(v,g)=\sum c_{k,h} 
\phi_{ih}^V(v,e_{\mu_{kh}}\cdot g) \phi_{h
j}^U(u,e_{-\mu_{kh}}\cdot f),$$
where the braiding matrix $(c_{kh})$ and the phase corrections
$\mu_{kh}$ also depend on $i$, $k$, $h$ and $j$. For
$f\in C_c^\infty(S^1\backslash\{1\})$, the expression $e_\mu f$ is
defined (unambiguously) by cutting the circle at $1$, so that $e_\mu\cdot
f(e^{i\theta}) = e^{i\mu\theta} f(e^{i\theta})$ for $\theta\in
(0,2\pi)$. To prove the braiding relation we introduce the formal
power series
$$F_k(z)=\sum_{n\ge 0} z^n
(\phi_{ik}^U(u,n)\phi_{kj}^V(v,-n)\xi,\eta),
\quad G_h(z)=\sum_{n\ge 0} z^n
(\phi_{ih}^V(v,n)\phi_{h j}^U(u,-n)\xi,\eta),$$
where $\xi$ and $\eta$ range over lowest energy vectors. These power
series are called (reduced) four--point functions and take values in
${\rm Hom}_G (U\otimes V \otimes V_j, V_i)$. Since the modes
$\phi_{ij}^U(n)$ and $\phi_{pq}^V(n)$ are uniformly bounded in norm,
they define 
holomorphic functions for $|z|<1$. We start by showing how the matrix
coefficients of products of primary fields can be recovered from
four--point functions.
\vskip .1in
\noindent \bf Proposition~1. \it Let $F_k(z)=\sum_{n\ge 0}
(\phi_{ik}^U(u,n)\phi_{kj}^V(v,-n)\xi,\eta)z^n = 
\sum F_nz^n$, convergent in $|z|<1$. If $f\in C^\infty_c(I)$, $g\in
C_c^\infty(I^c)$ and $\widetilde f(e^{i\theta})=f(e^{-i\theta})$, then 
$$(\phi^U_{ik}(u,f)\phi_{kj}^V(v,g)\xi,\eta)=\lim_{r\uparrow 1}{1\over
2\pi} \int_0^{2\pi}
\widetilde{f}\star g(e^{i\theta}) F_k(re^{i\theta})\, d\theta.$$
\vskip .05in
\noindent \bf Proof. \rm  If $f(z)=\sum f_nz^n$ and
$g(z)=g_nz^n$, then
$$(\phi^U_{ik}(u,f)\phi^V_{kj}(v,g)\xi,\eta) 
=\sum_{n\ge 0} f_n g_{-n} (\phi_{ik}^U(u,n)\phi_{kj}^V(v,-n)\xi,\eta)=
\lim_{r\uparrow 1} {1\over 2\pi} 
\int_0^{2\pi} \widetilde{f}\star g(e^ {i\theta}) F_k(re^{i\theta}) \,
d\theta.$$
\vskip .1in
\noindent \bf Corollary. \it Suppose that $f\in C_c^\infty(I)$, $g\in
C_c^\infty(I^c)$ and suppose further that $F_k(z)$ extends to a
continuous function on $S^1\backslash\{1\}$. Then
$$(\phi_{ik}^U(u,f)\phi_{kj}^V(v,g)\xi,\eta) ={1\over 2\pi}
\int_{0+}^{2\pi-} \widetilde{f}\star g(e^{i\theta}) F_k(e^{i\theta})
\, d\theta.$$
\vskip .05in
\noindent \bf Proof. \rm The assumptions on $f$ and $g$ imply
that the support of $\widetilde{f}\star g(e^{i\theta})$
is contained in 
$[\delta,2\pi-\delta]$ for some $\delta>0$, so the result follows. 
\vskip .1in
The next result explains how to translate from transport
equations for four point functions to braiding relations for
smeared primary fields. It is the analogue of the
Bargmann--Hall--Wightman theorem in axiomatic quantum field theory
([19],[34]). 
\vskip .1in 
\noindent \bf Proposition~2. \it Suppose that $U$ and $V$ are the
vector representation or its dual. Let
$$F_k(z)=\sum (\phi^U_{ik}(u,n)\phi^V_{kj}(v,-n)\xi,\eta) z^n,\quad
G_h(z) =\sum (\phi^V_{ih}(v,n) \phi^U_{hj}(u,-n)\xi,\eta) z^n,$$
where
$\xi$ and $\eta$ are lowest energy vectors. If $F_k(z),G_h(z^{-1})$ extend
to continuous functions on $S^1\backslash\{1\}$ with
$$F_k(e^{i\theta}) = \sum c_{kh} e^{i\mu_{kh}\theta} G_h(e^{-i\theta}),$$
where $\mu_{kh}\in {\Bbb R}$, then for $f\in C_c^\infty(0,\pi)$, $g\in
C_c^\infty(\pi,,2\pi)$ we
have
$$(\phi_{ik}^U(u,f)\phi_{kj}^V(v,g)\xi,\eta) = \sum c_{kh}
(\phi_{ih}^V(v,e_{\mu_{kh}}\cdot g) \phi_{hj}^U(u,e_{-\mu_{kh}}\cdot
f)\xi,\eta),$$
where $e_\mu(e^{i\theta})=e^{i\mu\theta}$ for $\theta\in (0,2\pi)$.
\vskip .05in
\noindent \bf Proof. \rm  For $\theta\in (0,2\pi)$ we have
$F_k(e^{i\theta}) = \sum c_{kh} e^{i\mu_{kh}\theta} G_h(e^{-i\theta})$.
Substituting in the equation of the corollary and changing variables
from $\theta$ to 
$2\pi-\theta$, we obtain 
$$(\phi_{ik}^U(u,f)\phi_{kj}^V(v,g)\xi,\eta) 
 =\sum c_{kh} {1\over 2\pi} \int_{0+}^{2\pi -}
e^{2i\mu_{kh}\pi} e^{-i\mu_{kh}\theta}\widetilde{g}\star 
f(e^{i\theta}) G_k(e^{i\theta})\, d\theta.$$
It can be checked directly that $e_{-\mu}\cdot (\widetilde{g}\star f)
= e^{-2\pi i \mu} 
\widetilde{e_{\mu}\cdot g} \star 
(e_{-\mu}\cdot f)$ (the corresponding identity is trivial for point measures
supported in $(0,\pi)$ and $(\pi,2\pi)$ and follows in general by weak
continuity); this implies the braiding relation.
\vskip .1in
A standard argument with lowering and raising operators allows us to
extend this braiding relation to arbitrary finite energy vectors $\xi$
and $\eta$ and hence arbitrary vectors.
\vskip .1in
\noindent \bf Proposition~3. \it If 
$$(\phi_{ik}^U(u,f)\phi_{kj}^V(v,g)\xi,\eta) =\sum c_{kh} 
(\phi_{ih}^V(v,e_{\mu_{kh}} \cdot g) \phi_{hj}^U
(u,e_{-\mu_{kh}}\cdot f) \xi,\eta),$$
for $\xi,\eta$ lowest energy vectors, then the relation holds for
all vectors $\xi,\eta$. 
\vskip .05in
\noindent \bf Proof. \rm By bilinearity and continuity, it will
suffice to prove the braiding relation for finite energy vectors
$\xi,\eta$. Suppose that $\eta$ is a lowest energy vector. We start by
proving that the braiding relations holds for $\xi,\eta$ by induction
on the energy of $\xi$. When $\xi$ has lowest energy, the relation is
true by assumption. Now suppose that the relation holds for
$\xi_1,\eta$. Let us prove it for $\xi,\eta$ where $\xi=X(-n)\xi_1$,
where $n>0$. Then
$$\eqalign{(\phi_{ik}^U&(u,f)\phi_{kj}^V(v,g)\xi,\eta)
 = (\phi_{ik}^U(u,f)\phi_{kj}^V(v,g)X(-n)\xi_1,\eta)\cr
&=-(\phi_{ik}^U(u,f)\phi_{kj}^V(Xv,e_{-n}\cdot g)\xi_1,\eta)
-(\phi_{ik}^U(Xu,e_{-n}\cdot f)\phi_{kj}^V(v,g)\xi_1,\eta)\cr
&=-\sum_h c_{kh}
(\phi_{ih}^V(Xv,e_{\mu_{kh}}e_{-n}g)\phi_{hj}^U(u,e_{-\mu_{kh}}f)\xi_1,\eta) 
-\sum _h c_{kh}
(\phi_{ih}^V(v,e_{\mu_{kh}}g)\phi_{hj}^U(u,e_{-\mu_{kh}}e_{-n}f)
\xi_1,\eta)\cr
&= \sum_h c_{kh}
(\phi_{ih}^V(v,e_{\mu_{kh}}g)\phi_{hj}^U(u,e_{-\mu_{kh}}f)\xi,\eta).\cr}$$
This proves the braiding relation for all $\xi$ and all lowest energy
vectors $\eta$. A similar inductive argument shows the braiding
relation holds for all $\xi$ and all $\eta$.

\vskip .1in
\noindent \bf Corollary~1. \it If $f$ and $g$ are supported in
$S^1\backslash\{1\}$ and the support of $g$ is anticlockwise after the
support of $f$, then 
$$\phi_{ik}^U(u,f)\phi_{kj}^V(v,g)=\sum c_{kh} 
\phi_{ih}^V(v,e_{\mu_{kh}} \cdot g) \phi_{hj}^U
(u,e_{-\mu_{kh}}\cdot f).$$
\vskip .05in
\noindent \bf Proof. \rm This result follows immediately from the
proposition, using a partition of unity and rotating if
necessary so that neither the support of $f$ nor $g$ pass $1$.
\vskip .1in
\noindent \bf Corollary~2. \it  If $f$ and $g$ are supported in
$S^1\backslash\{1\}$ and the support of $g$ is anticlockwise after the
support of $f$, then 
$$\phi_{ik}^U(u,f)\phi_{kj}^V(v,g)=\sum d_{kh} 
\phi_{ih}^V(v,e_{\mu_{kh}} \cdot g) \phi_{hj}^U
(u,e_{-\mu_{kh}}\cdot f),$$
where $d_{kh}=e^{2\pi i\mu_{kh}} c_{kh}$.
\vskip .05in
\noindent \bf Proof. \rm This follows by applying a rotation of
$180^{\rm o}$ in the proposition and then repeating the reasoning in the proof
of corollary~1. 
\vskip .1in
\bf \noindent 27. Sugawara's formula. \it Let $H$ be a positive
energy irreducible representation at level $\ell$ and let $(X_i)$ be
an orthonormal basis of $\g$. Let $L_0$ be the
operator defined on $H^0$ by
$$L_0={1\over N+\ell} \left(-\sum_{i} {1\over 2} X_i(0)X_i(0)
-\sum_{n>0}\sum_i  X_i(-n)X_i(n)\right).$$
Then $L_0 = D + \Delta/2(N+\ell)$ if $-\sum_i X_i(0)X_i(0)$ acts
on $H(0)$ as multiplication by $\Delta$. 
\vskip .1in
\bf \noindent Remark. \rm Note that the operator $C=-\sum X_i X_i=\sum
E_{ij}E_{ji} -(\sum E_{ii})^2/N$ acts
in $V_f$ as the constant 
$$\Delta_f=[\sum f_i^2 +f_i(N-2i +1)] - (\sum f_i)^2/N.$$
In particular, for the adjoint representation on $\g$ ($f_1=1$,
$f_2=f_3=\cdots= f_{N-1}=0$, $f_N=-1$) we have $\Delta=2N$.
\vskip .1in

\noindent \bf Proof (cf [28]). \rm Since $\sum_i X_i(a) X_i(b)$ is
independent of the orthonormal basis $(X_i)$, it commutes with $G$ and
hence each $X(0)$ for $X\in \g$. Thus $\sum_i[X,X_i](a) X_i(b) + X_i(a)
[X,X_i](b) =0$ for all $a,b$. If
$A=\sum_i{1\over 2}  X_i(0)X_i(0) +\sum_{n>0} X_i(-n) X_i(n)$, then
using the above relation we get
$$\eqalign{[X(1),A]
&=N\ell X(1) +\sum_i {1\over 2}([X,X_i](1) X_i(0) + X_i(0)
[X,X_i](1)) \cr
& \quad +\sum_n[X,X_i](-n+1) X_i(n) + X_i(-n) [X,X_i](n+1)\cr
& = N\ell X(1) +{1\over 2} \sum _i [[X,X_i](1),X_i(0)]= N\ell X(1)
+{1\over 2} \sum_i [[X,X_i],X_i](0),\cr}$$ 
since $([X,X_i],X_i)=0$ by invariance of $(\cdot,\cdot)$. 
Hence $[X(1),A]=(N+\ell) X(1)$, since $-\sum_i {\rm ad}(X_i)^2 =2N$.
Now formally $X(1)^*=-X(-1)$ and 
$A^*=A$, so taking adjoints we get $[X(-1),A]=-(N+\ell)X(-1)$, so that
$(N+\ell)D+A$ commutes with all $X(\pm 1)$'s. Since
$[\g,\g]=\g$, these generate $L^0\g$, and hence $(N+\ell)D +A
=\lambda I$ for some $\lambda\in {\Bbb C}$. Evaluating on $H(0)$, we
get $\lambda=-\Delta/2$.

\vskip .1in

\noindent \bf 28. The Knizhnik--Zamolodchikov ODE (cf [22]). \rm Let
$\phi(a,n):H_j^0 \rightarrow H_k^0$ and $\phi(b,m):H_k^0 
\rightarrow H_i^0 $ be primary fields of charges $V_2$ and $V_3$
respectively.
Let $a_{nm}$ be the matrix coefficient $a_{nm}=\break (\phi(v_2,n) \phi(v_3,m)
v_4,  v_1)$, where $V_4=H_j(0)$ and $V_1=H_i(0)$. Since $Dv_4=0=Dv_1$ and 
$[D,\phi(v_2,n)]=-n\phi(v_2,n)$, $[D,\phi(v_3,m)]=-m\phi(v_3,m)$, it
follows immediately that $a_{n,m} = 0$ unless $n+m=0$. Moreover
$\phi(a,m)v=0$ if $n>0$, so that $a_{nm}=0$ if $m>0$. We define
four commuting actions of $SU(N)$ on ${\rm Hom}(V_2\otimes V_3\otimes
V_4,V_1)$ by
$\pi_1(g) T=gT$, $\pi_2(g)T=T(g^{-1}\otimes I\otimes I)$,
$\pi_3(g)T=T(I\otimes g^{-1} \otimes I)$ and $\pi_4(g) T=T(I\otimes
I\otimes g^{-1})$. Thus  $\pi_1(g)\pi_2(g)\pi_3(g)\pi_4(g)T=T$ 
if $T$ is $G$--equivariant.

Now let $(X_i)$ be an orthonormal basis of $\g$ and
define operators $\Omega_{ij}$ on $W={\rm Hom}_G(V_2\otimes V_3\otimes
V_4,V_1)$ as $- \sum \pi_i(X_k)\pi_j(X_k)$.
Thus $\Omega_{ij}=\Omega_{ji}$. Moreover, if $i,j,k$ are
distinct, then 
$\Omega_{ij} + \Omega_{jk} + \Omega_{ki}=h$ on $W$,
where $h$ is a constant. In fact, if $m$ is the missing index,
$$\eqalign{\Omega_{ij} + \Omega_{jk} + \Omega_{ki}
 & = -{1\over 2}\left[\sum_p (\pi_i(X_p) + \pi_j(X_p) +\pi_k(X_p))^2 -
 \pi_i(X_p)^2 -\pi_j(X_p)^2 -\pi_k(X_p)^2\right]\otimes I\cr
 &=-{1\over 2}[\sum (-\pi_m(X_p))^2 +\Delta_i
 +\Delta_j+\Delta_k]=(\Delta_m-\Delta_i-\Delta_j-\Delta_k)/2,\cr}$$
since $\g$ acts trivially on $W$.

\vskip .1in
\noindent \bf Theorem. \it The formal power series 
$f(v,z)=\sum_{n\ge 0} (\phi(v_2,n) \phi(v_3,-n) v_4, v_1)z^n$,
taking values in $W$, satisfies the Knizhnik--Zamolodchikov ODE
$$(N+\ell) {df\over dz} = \left(
{\Omega_{34}-(\Delta_k -\Delta_3-\Delta_4)/2\over z} +{\Omega_{23}\over
z-1}\right) f(z).$$
\vskip .05in
\noindent \bf Proof. \rm This is proved by inserting $D$ in
the 4--point function $f(z)$ and comparing it with the Sugawara
formula $D=L_0-h$. In fact
$zf^\prime(z)
 =\sum_{n\ge 0} (\phi(v_2,n)D\phi(v_3,-n)v_4, v_1) z^n$, 
since $[D,\phi(v,m)]=-m\phi(v,m)$ and $Dv_4=0$. Now on $H_k^0$ we have
$D=L_0-h$ where $h=\Delta_k/2(N+\ell)$, so that
$$\eqalign{zf^\prime(z)
 & =-h\cdot f(z) -(N+\ell)^{-1} \sum_{n\ge 0,i}[\sum_{m>0} 
 (\phi(v_2,n)X_i(-m)X_i(m)\phi(v_3,-n)v_4, v_1) z^n \cr
  & \qquad+{1\over 2}
   (\phi(v_2,n)X_i(0)X_i(0)\phi(v_3,-n)v_4, v_1) z^n].\cr}$$ 
Now $[X(n),\phi(v,m)]=\phi(X\cdot v,n+m)$, so that
$\phi(v_2,n)X_i(m) = X_i(m) \phi(v_2,n) -\phi(X_i\cdot v_2,n+m)$ and
$X_i(m) \phi(v_3,n)=\phi(v_3,n) X_i(m) + \phi(X_i\cdot v_3,n+m)$.
Substituting in these expressions, we get
$$\eqalign{zf^\prime(z) 
 & = -h\cdot f(z) +(N+\ell)^{-1} \sum_{n\le 0,i} 
    \sum_{m>0}  (\phi(X_iv_2,n-m)\phi(X_iv_3,-n+m)v_4, v_1) z^n\cr
 &  \qquad- (2(N+\ell))^{-1} \sum_{n\ge 0,i} ((X_i(0)\phi(v_2,n)
  -\phi(X_iv_2,n)) (\phi(v_3,-n)X_i(0)+\phi(X_iv_3,-n))v_4,v_1)
  z^n\cr
 &=( N+\ell)^{-1} (-\Delta_k/2 - {1\over 2}\Omega_{23} {z\over 1-z}
  -{1\over 2}(\Omega_{23} 
  +\Omega_{13} + \Omega_{14} +\Omega_{24}))f(z)\cr
 &  =(N+\ell)^{-1}({\Omega_{34}-{1\over 2}(\Delta_k-\Delta_3-\Delta_4)}
   +\Omega_{23}{z\over z-1})f(z). \cr}$$

\vskip .1in

\noindent \bf 29. Braiding relations between vector and dual vector
primary fields. \rm Consider the four--point functions
$F_k(z)=\sum_{n\ge 0}
(\phi_{ik}^U(u,-n)\phi_{kj}^V(v,n)\xi,\eta)z^n$ and
$G_h(z)=\sum_{n\ge 0}
(\phi_{ih}^V(v,-n)\phi_{hj}^U(u,n)\xi,\eta)z^n$,
where the charges $U$ and $V$ are either ${\Bbb C}^N$ or
its dual. Thus any $V_k$ appears with multiplicity one in
the tensor product $V\otimes V_j$ or $U\otimes V_j$, and all but
possibly one of these summands will be permissible at level $\ell$.
\vskip .1in
\noindent \bf Proposition. \it (a) $f_k(z)=z^{\lambda_k}F_k(z)$
satisfies the KZ ODE
$$(N+\ell){df\over dz} ={\Omega_{vj}\over z}f(z)
+{\Omega_{uv}\over z-1}f(z),$$
where $\lambda_k=(\Delta_k
-\Delta_v-\Delta_j)/2(N+\ell)$ is the eigenvalue of $(N+\ell)^{-1}
\Omega_{vj}$ corresponding to the summand $V_k\subset V\otimes V_j$.

\noindent (b) $g_h(z) = z^{\mu_h} G_h(z^{-1})$ satisfies the same ODE,
where $\mu_h=(\Delta_i-\Delta_v
-\Delta_h)/2(N+\ell)$ is the eigenvalue of $(N+\ell)^{-1} (\Omega_{vj} +
\Omega_{uv})$ corresponding to the summand $V_h\subset U\otimes V_j$. 
\vskip .1in
\noindent \bf Proof. \rm (a) Since
$$\Omega_{vj} =-\sum \pi_v(X_i)\pi_j(X_i) = -{1\over 2} \sum(\pi_v(X)
+ \pi_j(X))^2 + {1\over 2} \sum \pi_q(X_i)^2 +{1\over 2}
\sum\pi_j(X_i)^2,$$
$(N+\ell)^{-1}\Omega_{vj}$ acts as the scalar $\lambda_k=(\Delta_k
-\Delta_v-\Delta_j)/2(N+\ell)$ on the subspace $V_k\subset V\otimes V_j$.
Thus the result follows from the previous section.

\noindent (b) Similarly $\nu_h=(\Delta_h -\Delta_u-\Delta_j)/2(N+\ell)$ is the
eigenvalue of $(N+\ell)^{-1}\Omega_{uj}$ corresponding to the summand
$V_h$ of $U\otimes V_j$. Let
$\mu=(\Delta_i-\Delta_u-\Delta_v-\Delta_j)/2(N+\ell)$. It is easy to
verify that $h(z)=z^{\mu-\nu_h}G_h(z^{-1})$ satisfies the same ODE,
since $(N+\ell)^{-1} (\Omega_{uv}+\Omega_{vj} +\Omega_{ju})=\mu$ on
${\rm Hom}_G(U\otimes V\otimes V_j,V_i)$. Here
$\mu_h=\mu-\nu_h=(\Delta_i-\Delta_v-\Delta_h)/2(N+\ell)$ is the
eigenvalue of $(N+\ell)^{-1} (\Omega_{vj} +\Omega_{uv})$ corresponding
to the summand $V_h\subset U\otimes V_j$. 

\vskip .1in

Thus the solutions $f_k(z)$ form part of a complete set of
solutions about $0$ of the KZ ODE; and the solutions $g_h(z)$
form part of a solution set about $\infty$ of the same ODE. They may
only form part, because 
one of the summands $V_k$ or $V_h$, and hence eigenvalues $\lambda_k$
or $\mu_h$, might
correspond to a representation not permissible 
at level $\ell$; there can be at most one such summand. Let $f_k(z)$
and $g_h(z)$ denote the two complete sets of solutions, regardless of
whether the eigenvalues $\lambda_k$ or $\mu_h$ are permissible. They
define holomorphic functions in ${\Bbb C}\backslash [0,\infty)$. 
Let $c_{kh}$ be the transport matrix 
relating the solutions at $0$ to the solutions around $\infty$, so
that $f_k(z) =\sum c_{kh} g_h(z)$ for $z\in {\Bbb C}\backslash
[0,\infty)$. Thus
$F_k(z)=\sum c_{kh} z^{\mu_{kh}} G_h(z^{-1})$, for $z\in {\Bbb C}\backslash
[0,\infty)$ where $\mu_{kh}=\mu_h-\lambda_k
=(\Delta_i+\Delta_j-\Delta_h-\Delta_k)/2(N+\ell)$. Whenever
an $F_k$ or $G_h$ does not correspond to a product of primary fields 
(because $V_k$ or $V_h$ is not permissible at level $\ell$), we will
find that the corresponding transport coefficient $c_{kh}$ is zero.
(This is not accidental. As explained in [41], there is
an algebraic boundary condition which picks out the solutions that
arise as four--point functions.) All the examples we will
consider will be those for which the theory of the previous chapter is
applicable.  

\vskip .1in

\noindent \bf Theorem~A (generalised hypergeometric braiding). \it Let
$F\in L^2(I,V)$ and $G\in L^2(J,V^*)$ where $I$ and $J$ are intervals in
$S^{1}\backslash\{1\}$ with $J$ anticlockwise after $I$. Then
$$\phi_{gf}^{\square}(F)\phi_{f g}^{\overline{\square}} (G)=\sum
\nu_{fh} \phi^{\overline{\square}}_{gh}(e_{\mu_{fh}}
G)\phi^\square_{hg} (e_{-\mu_{fh}}F)$$
with $\nu_{fh}\ne 0$, if $h> g$ and $\mu_{fh}= (2\Delta_{g}
-\Delta_f - \Delta_h)/2(N+\ell)$. 
\vskip .05in
\noindent \bf Proof. \rm The KZ ODE reads
$$(N+\ell) {df\over dz} = {\Omega_{\overline{\square}f}f(z)\over
z} + {\Omega_{\square\overline{\square}}f(z)\over z-1},$$
where $f(z)$ takes values in $W={\rm 
Hom}_G(V_{\square}\otimes V_{\overline{\square}}\otimes V_g,V_g)$. 
Now the eigenvalue of $\Omega_{\square\overline{\square}}$
corresponding to the trivial representation is $(0-\Delta_\square
-\Delta_{\overline{\square}})/2= N^{-1}-N$ and has multiplicity one,
while that corresponding 
to the adjoint representation is $(\Delta_{{\rm Ad}} - \Delta_\square
-\Delta_{\overline{\square}})/2= N^{-1}$ with multiplicity at most
$N-1$. Thus $\Omega_{\square\overline{\square}}=N^{-1} -NQ $, if $Q$
is the rank one projection in $W$ corresponding to the trivial
representation. So
$$-(N+\ell)^{-1}\Omega_{\square\overline{\square}}={N\over N+\ell}
Q-{1\over N(N+\ell)}.$$
Thus $\alpha=1/N(N+\ell)$ and $\beta=N/(N+\ell)$ (in the notation of
section~18). 

We next check that $A=(N+\ell)^{-1}\Omega_{\overline{\square}f}$ and $Q$
are in general position. In fact if we identify $W$ with ${\rm
End}_G(V_g\otimes V_\square)$, then the inner product becomes ${\rm
Tr}(xy^*)$. The identity operator $I$ is the generator of the range of
$Q$ with $Q(x)$ proportional to ${\rm Tr}(x)$. The eigenvectors of $A$
are just given by the orthogonal projections $e_g$ onto the
irreducible summands $V_g$ of $V_f\otimes V_\square$. Since ${\rm
Tr}(e_g)>0$, it follows that
$A$ and $Q$ are in general position.

The eigenvalues of $A$ are given by $\lambda_f= (\Delta_f -\Delta_{\overline{\square}} -
\Delta_g)/2(N+\ell)$, so that $|\lambda_f-\lambda_{f_1}|=|\Delta_f -
\Delta_{f_1}|/2(N+\ell)$. This has the form $|g_i-g_j-i +j|/(N+\ell)$
for $i\ne j$, if $f$ and $f_1$ are obtained by removing boxes from the
$i$th and $j$th rows of $g$. Since $g_i + N-i$ is strictly increasing
and $g_1-g_N\le 
\ell$, the maximum possible difference is $|g_N-g_1-N+1|/(N+\ell)=1 -
(N+\ell)^{-1} <1$. Hence $0<|\lambda_{f}-\lambda_{f_1}|<1$ if $f\ne
f_1$.  Similarly $\mu_h=(\Delta_g
-\Delta_h-\Delta_{\overline{\square}})/2(N+\ell)$ and the difference
$|\mu_h-\mu_{h_1}|$ has the form $|g_i-g_j-i +j|/(N+\ell)$
for $i\ne j$, if $h$ and $h_1$ are obtained by adding boxes to the
$i$th and $j$th rows of $g$. Hence $0<|\mu_h-\mu_{h_1}|<1$ if $h\ne h_1$. 
\vskip .05in
\noindent \bf Caveat. \rm The indexing set for the $f_j$' and $h_k$'s
are distinct, even though they have the same cardinality. This is 
easy to see if one draws $f$ as a Young diagram. The $f_j$'s
correspond to corners pointing north--west while the $h_k$'s
correspond to corners pointing south--east.
\vskip .05in
The anomaly $\mu_{fh}$ is given by the stated formula by our preamble,
so it only remains to check that permitted terms $c_{fh}$ are
non--zero and forbidden terms zero. In fact the numerator
is always non--zero
because $\Gamma(x)\ne 0$ for all $x\notin -{\Bbb N}$. Thus the only 
way $c_{fh}$ can vanish is if one of the arguments of $\Gamma$
in the denominator
$\prod_{\ell\ne j} \Gamma( \lambda_i -\mu_\ell +\alpha +1)\,
\prod_{k\ne i} \Gamma(\mu_j - \lambda_k-\alpha)$
is a non--positive integer. Now $\mu_h=(\Delta_g
-\Delta_h-\Delta_{\overline{\square}})/2(N+\ell)$ and 
$\lambda_f= (\Delta_f -\Delta_{\overline{\square}} -
\Delta_g)/2(N+\ell)$. Suppose that $h$ is obtained by adding a box to
the $i$th row of $g$ and $f$ is obtained by removing a box from the
$j$th row of $g$. Then $\lambda_f-\mu_h=(N+\ell)^{-1}(g_j-g_i + 1 + i
-j-N^{-1})$. Thus
$$\lambda_f -\mu_h +\alpha = (N+\ell)^{-1} (g_j - g_i + 1 + i-j).$$
This has modulus less than $1$ unless $i=1$, $j=N$ and $g_1-g_N=\ell$,
when it gives $-1$. It is then easy to see that if $f$ or $h$ is
non--permissible, the corresponding coefficient vanishes and otherwise
it is non--zero.
\vskip .1in

The next example of braiding could have been done using the classical
theory of the hypergeometric function ([16],[45]);
however, since the 
equation is in matrix form and some knowledge of Young's orthogonal form is
required to translate this matrix equation into the hypergeometric
equation, it is much simpler to use the matrix and eigenvalue
techniques. 
\vskip .1in

\noindent \bf Theorem~B (hypergeometric braiding). \it Let
$F\in L^2(I,V)$ and $G\in L^2(J,V)$ where $I$ and $J$ are intervals in
$S^{1}\backslash\{1\}$ with $J$ anticlockwise after $I$. Then
$\phi_{hg}^{\square}(F)\phi_{gf }^{\square} (G)=\sum
\mu_{gg_1} \phi^{\square}_{hg_1}(e_{\alpha_{gg_1}}
G)\phi^\square_{g_1f} (e_{-\alpha_{gg_1}}F)$
with $\mu_{gg_1}\ne 0$, if $h> g,g_1> f$ and $\alpha_{gg_1}=
(\Delta_h +\Delta_f -\Delta_g -\Delta_{g_1})/2(N+\ell) $. 
\vskip .05in
\noindent \bf Proof. \rm In this case $W={\rm Hom}_G(V_\square\otimes
V_\square\otimes V_f,V_h)$ has dimension $2$. The eigenvalues of
$(N+\ell)^{-1}\Omega_{\square\square}$ correspond to the summands
$V_{\YD\R2,.}$ 
and $V_{\YD\R1,\R1,.}$. We have $\lambda_{\YD\R2,.}=(\Delta_{\YD\R2,.}
-2\Delta_\square)/2(N+\ell)=(N-1)/N(N+\ell)$ and 
$\lambda_{\YD\R1,\R1,.}=(\Delta_{\YD\R1,\R1,.}-2\Delta_\square)/2(N+\ell)=
(-N-1)/N(N+\ell)$. If $Q$ is the projection corresponding to
$V_{\YD\R1,\R1,.}$ and $\beta Q -\alpha
I=-(N+\ell)^{-1}\Omega_{\square\square}$, then $\beta=2/N(N+\ell)$ and
$\alpha= (N-1)/N(N+\ell)$. 

We have $\lambda_g=(\Delta_g -\Delta_f -\Delta_\square)/2(N+\ell)$ and
$\mu_g=(\Delta_h -\Delta_g -\Delta_\square)/2(N+\ell)$. Thus
$|\lambda_g-\lambda_{g_1}|=|\mu_g-\mu_{g_1}|=|\Delta_g-\Delta_{g_1}|/
2(N+\ell)= |f_i -i - f_j+j|/(N+\ell)$, if $g$ and $g_1$ are obtained
by adding boxes to $f$ in the $i$th and $j$th rows. As above, it
follows that $|\lambda_g-\lambda_{g_1}|=|\mu_g-\mu_{g_1}|<1$.

We next check that the operators $A=(N+\ell)^{-1}\Omega_{\square f}$
and $Q$ are in 
general position. The operator $\Omega_{\square\square}$ is a linear
combination of the identity operator ${\rm id}$ and $\sigma$, where
$\sigma(T)=T(S\otimes 
I)$ and $S$ is the flip on $V_\square\otimes 
V_\square$. The operators $T_i$ in $W$ which diagonalise $\Omega_{\square
f}$ are obtained by composing intertwiners $V_\square\otimes
V_f\rightarrow V_{g_i}$ and $V_\square \otimes V_{g_i}\rightarrow
V_h$. These intertwiners are specified by their action on vectors
$e_i\otimes v$ where $(e_i)$ is a basis of $V_\square$ and $v$ is a
highest weight vector. If $g_1$ and $g_2$ are obtained by adding boxes
to $f$ in rows $i$ and $j$ with $i,j$, it is easy to see that
$T_2(e_i\otimes e_j\otimes v_f)$ is a non--zero highest weight vector
in $V_h$ while $\sigma(T_2)(e_i\otimes e_j\otimes v_f)=
T_2(e_j\otimes e_i \otimes v_f)=0$. So $T_2$ is not an eigenvector of
$\sigma$. This proves that $A$ and $Q$ are in general position.
 
The anomaly $\alpha_{gg_1}$ is as stated by our preamble, so it only
remains to check that permitted terms $c_{gg_1}$ are non--zero and
forbidden terms zero. As above, a term can vanish iff one of the
arguments in the denominator $\Gamma(\lambda_{g}-\mu_{g_1^\prime}
+\alpha +1) 
\Gamma(\mu_{g_1^\prime}-\lambda_{g^\prime} -\alpha)$ is a non--positive
integer (where $g^\prime$ denotes the other diagram to $g$ between $f$
and $h$).  Now $\lambda_g-\mu_{g_1}=(\Delta_g+\Delta_{g_1} -
\Delta_f-\Delta_h)/2(N+\ell)$. Hence
$\lambda_{g}-\mu_{g^\prime}=1/N(N+\ell)$, so that
$\lambda_{g}-\mu_{g^\prime}+\alpha+1 =1 +(N+\ell)^{-1}$ and 
$\mu_{g^\prime} -\lambda_{g}-\alpha = -(N+\ell)^{-1}$. This shows that, if
$g$ is permissible, none of the arguments is a non--positive integer
and hence that $c_{gg}\ne 0$. On the other hand
$\lambda_{g}-\mu_g=(f_i-i -f_j +j)/(N+\ell) +1/N(N+\ell)$, if $g$ is
obtained by adding a box to the $i$th row of $f$. Thus
$\lambda_g-\mu_g+\alpha+1=1+ (f_i-i-f_j+j+1)/(N+\ell)$, which can
never be a non--positive integer, while
$$\mu_{g^\prime}-\lambda_{g^\prime}-\alpha=(f_i-i-f_j+j-1)/(N+\ell).$$
This has modulus less than $1$ unless $i=N$, $j=1$ and $f_1-f_N=\ell$,
when it gives $-1$. This is the critical case where $g$ is permissible
(it is obtained by adding a box to the last row of $f$) while
$g^\prime$ is inadmissible (it is obtained by adding a box to the
first row of $f$). In this case therefore $c_{gg^\prime}=0$ while in
all other cases the coefficient is non--zero.

\vskip .1in
\noindent \bf Theorem~C (Abelian braiding). \it Let $F\in L^2(I,V)$ and
$G\in  L^2(I^c,V^*)$. Let $g\ne g_1$ be signatures, permissible at
level $\ell$, 
obtained by adding one box to $f$. Then
$\phi_{gf}^\square(F)\phi_{fg_1}^{\overline{\square}} (G)=\varepsilon
\phi_{gh}^{\overline{\square}}(e_{\mu}G)
\phi_{hg_1}^\square(e_{-\mu}F)$
with $\varepsilon\ne 0$ and $\mu= (\Delta_g +\Delta_{g_1} - \Delta_f
-\Delta_h)/2$. 
\vskip .05in
\noindent \bf Proof. \rm The corresponding ODE takes values in the
one--dimensional space ${\rm Hom}_G(V_\square\otimes
V_{\overline{\square}} \otimes V_{g_1}, V_g)$ so $\varepsilon$ must be
non--zero and $\mu$ is as stated by our preamble.
\vskip .1in
\noindent \bf Theorem~D (Abelian braiding). \it Suppose that $g$ is
the unique signature such that $h> g>f$, so that $h$ is
obtained either by adding two boxes in the same row of $f$ (symmetric case
$+$) or in the same column column (antisymmetric case $-$). Let
$F\in L^2(I,V)$ and $G\in L^2(J,V)$ where $I$ and $J$ are intervals in
$S^{1}\backslash\{1\}$ with $J$ anticlockwise after $I$. Then there
are non--zero constants $\delta_+\ne \delta_-$ depending only on the
case such that 
$$\phi_{hg}^{\square}(F)\phi_{gf }^{\square} (G)=
\delta_{\pm} \phi^{\square}_{hg}(e_{\alpha}
G)\phi^\square_{gf} (e_{-\alpha}F)$$
with $\delta_\pm\ne 0$ and $\alpha=
(\Delta_h +\Delta_f -2\Delta_g)/2 $. In fact $\delta_\pm
=e^{i\pi\nu_\pm}$ where $\nu_\pm=(\pm N -1)/N(N+\ell)$. 
\vskip .05in
\noindent \bf Proof. \rm We use the same reasoning as in the proof of
Theorem~C. The ODE is now a scalar equation $f^\prime=(\lambda_g z^{-1}
+ \nu_{\pm} (z-1)^{-1})f$. The $\nu_+$ and $\nu_-$ are the eigenvalues of
$(N+\ell)^{-1}\Omega_{\square\square}$ corresponding to the summands
$V_{\YD\R2,.}$ 
and $V_{\YD\R1,\R1,.}$ respectively. So $\nu_\pm=(\pm N
-1)/N(N+\ell)$. The normalised solution near $0$ of the ODE is
$z^{\lambda_g}(1-z)^{\nu_\pm}$ while near $\infty$ it is
$z^{\lambda_g+\nu_\pm} (1-z^{-1})^{\nu_\pm}$. Taking $z=-x$, with $x$
real and positive, it follows immediately that the transport
coefficient is $e^{i\pi\nu_{\pm}}$. 
\vskip .1in

\noindent \bf Summary of braiding properties. \rm If we define
$a_{gf}^\square=\phi^\square_{gf}(e_{-\alpha} F)$ where
$\alpha=(\Delta_g-\Delta_f-\Delta_\square)/2(N+\ell)$ and 
$a_{fg}^{\overline{\square}}=\phi_{fg}^{\overline{\square}}
(e_\alpha F^*)$, then the adjoint relation between these two primary
fields implies that
$(a_{gf}^{\square})^*=a_{fg}^{\overline{\square}}$. Incorporating the
anomalies $e_\mu$ into the smeared primary fields in this way, the
braiding properties established above for vector and dual vector
primary fields may be stated in the following form.

\vskip .1in
\noindent \bf Theorem. \it Let $(a_{ij})$, $(b_{ij})$ denote
vector primary fields smeared in intervals $I$ and $J$ in
$S^{1}\backslash\{1\}$ with $J$ anticlockwise after $I$.
\vskip .05in

\noindent (a) $a_{gf}b_{g_1f}^*=\sum \nu_h b_{hg}^*a_{h g_1}$
with $\nu_h\ne 0$, if $h> g,g_1>f$. 

\noindent (b) $a_{gf}b_{fh} =\sum \mu_{f_1} b_{gf_1} a_{f_1h}$ with
$\mu_{f_1}\ne 0$ if $h<f_1<g$.

\noindent (c) $a_{gf}b_{g_1 f}^*=\varepsilon b_{hg}^*a_{hg_{1}}$ with
$\varepsilon \ne 0$.

\noindent (d) $a_{hg}b_{gf} =\delta_{\pm} b_{hg} a_{gf}$ where
$\delta_+\ne \delta_-$ are non--zero, with $+$ if
$h$ is obtained from $f$ by adding two boxes in the same row and $-$
if they are in the same column.
\vskip .1in
\rm Note that (c) and (d) may be regarded as degenerate versions of (a) 
and (b) respectively so may be combined. Rotating through $180^{{\rm
o}}$ as before, or taking adjoints or simply rewriting the above
equations, we obtain our final result.  (For simplicity we have
suppressed the resulting phase changes in the coefficients.) 
\vskip .1in
\noindent \bf Corollary. \it Let $(a_{ij})$, $(b_{ij})$ denote
vector primary fields smeared in intervals $I$ and $J$ in
$S^{1}\backslash\{1\}$ with $J$ anticlockwise after $I$.

\vskip .05in
\noindent (a) $b_{gf}a_{g_1f}^*=\sum \nu_h a_{hg}^*b_{h g_1}$
with $\nu_h\ne 0$ if $h> g,g_1$ is permissible.

\noindent (b) $b_{gf}a_{fh} =\sum \mu_{f_1} a_{gf_1} b_{f_1 h}$ with
$\mu_{f_1}\ne 0$ if $h< f_1< g$.

\noindent (c) $b_{gf}a_{g_1 f}^*=\varepsilon a_{hg}^*b_{hg_1}$ with
$\varepsilon\ne 0$.

\noindent (d) $b_{hg}a_{fh} =\delta^{-1}_\pm a_{hg}b_{gf}$ with
$\delta_+\ne \delta_-$ non--zero.
\vfill\eject
\bf\noindent CHAPTER V.~CONNES FUSION OF POSITIVE ENERGY REPRESENTATIONS.
\vskip.2in
\noindent \bf 30. Definition and elementary properties of Connes
fusion for positive energy representations. \rm In [40] and [41]
we gave a fairly extensive treatment of Connes' tensor
product operation on bimodules over von Neumann algebras. It was then
applied to define a fusion operation on positive energy
representations of ${\cal L}G$. Here we give a simplified direct
treatment of fusion with more emphasis on loop groups than von Neumann
algebras. Let $X$ and $Y$ be positive energy representations of $LG$
at level $\ell$.  To define their fusion, we consider intertwiners (or
fields) $x\in {\cal
X}={\rm Hom}_{{\cal L}_{I^c}G}(H_0,X)$, $y\in {\cal Y}={\rm 
Hom}_{{\cal L}_IG}(H_0,Y)$ instead of the vectors (or states) $\xi=x\Omega$ and
$\eta=y\Omega$ they create from the vacuum. We define an inner product
on the algebraic tensor product ${\cal X}\otimes {\cal Y}$ 
by the {\it four--point formula} 
$\langle x_1\otimes y_1,x_2\otimes y_2\rangle  = (x_2^*x_1
y_2^*y_1\Omega,\Omega)$. 
\vskip .1in
\noindent \bf Lemma. \it The four--point formula defines an (pre--)inner
product on ${\cal X}\otimes {\cal Y}$. The Hilbert space completion
$H=X\boxtimes Y$ naturally admits a continuous unitary
representation $\pi$ of ${\cal L}^{\pm 1}G={\cal 
L}_IG\cdot  {\cal L}_{I^c}G$ of level $\ell$.
\vskip .1in
\noindent \bf Proof. \rm If $z=\sum x_i\otimes y_i\in {\cal X}\otimes
{\cal Y}$, then
$\langle z,z\rangle=\sum (x_i^*x_j y_i^* y_j \Omega,\Omega)$.
Now $x_{ij}=x_i^*x_j$ lies in $M=\pi_0({\cal
L}_{I^c}G)^\prime=\pi_0({\cal L}_I G)^{\prime\prime}$.
The operator $X=(x_{ij})\in M_n(M)$ is non--negative,
so has the form $X=A^*A$ for some $A=(a_{ij})\in M_n(M)$. 
Similarly, if $y_{ij}=y_i^*y_j\in M^\prime$, then $Y=(y_{ij})\in
M_n(M^\prime)$ can be written $Y=B^*B$ for some $B=(b_{ij})\in
M_n(M^\prime)$. Hence
$$\langle z,z\rangle=\sum_{p,q,i,j}
(a_{pi}^*a_{pj}b_{qi}^*b_{qj}\Omega,\Omega) =
\sum_{p,q} \|\sum_i a_{pi}b_{qi}\Omega\|^2\ge 0.$$
We next check that
${\cal L}_IG\cdot {\cal L}_{I^c}G$ acts continuously on ${\cal
X}\otimes {\cal Y}$, preserving the inner product. The action of
$g\cdot h$ on $x\otimes y$ is given by $(g\cdot h)(x\otimes
y)=gx\otimes hy$. It clearly preserves the inner product, so the group
action passes to the Hilbert space completion. Note that
since we have defined things on the level of central extensions, we have to
check that $\zeta\in {\Bbb T}={\cal L}_IG\cap {\cal L}_{I^c}G$ acts by
the correct scalar. This is immediate. Finally we must show that when
the matrix coefficients are continuous for vectors in ${\cal X}\otimes
{\cal Y}$ are continuous on ${\cal L}_IG\cdot {\cal L}_{I^c}G$. 
But 
$$\langle gx_1\otimes hy_1,x_2\otimes y_2\rangle
=(x_2^*gx_1y_2^*hy_1\Omega,\Omega)=(x_1y_2^*hy_1\Omega,g^*x_2\Omega).$$
Since the maps ${\cal L}_IG\rightarrow X$, $g\mapsto g^*x_2\Omega$
and ${\cal L}_{I^c}G\rightarrow Y$, $h\mapsto hy_1\Omega$ are
continuous, the matrix coefficient above is continuous.
\vskip .1in
\noindent \bf Lemma. \it There are canonical unitary isomorphisms
$H_0\boxtimes X \cong X \cong X\boxtimes H_0$. 
\vskip .1in
\noindent \bf Proof. \rm If $Y=H_0$, the unitary $X\boxtimes
H_0\rightarrow X$ is given by $x\otimes y\mapsto xy\Omega$ and the
unitary $H_0\boxtimes X\rightarrow X$ is given by $y\otimes
x\mapsto xy\Omega$.
\vskip .1in
\noindent \bf Lemma. \it If $J$ is another interval of the circle and
the above construction is accomplished
using the local loop groups 
${\cal L}_J G$ and ${\cal L}_{J^c}G$ to give a Hilbert space $K$ with a level $\ell$
unitary representation $\sigma$ of ${\cal L}_JG\cdot {\cal L}_{J^c}G$, then if
$\phi\in SU(1,1)$ carries $I$ onto $J$, there is a natural
unitary $U_\phi:H\rightarrow K$ that
$U_\phi(\pi(g))U_\phi^*=\sigma(g\circ\phi^{-1})$. 
\vskip .1in
\noindent \bf Proof. \rm Take $\phi\in SU(1,1)$ such that $\phi(I)=J$.
If $x\in {\cal X}_I={\rm Hom}_{{\cal L}_{I^c}G}(H_0,X)$ and $y\in{\cal
Y}_I={\rm 
Hom}_{{\cal L}_{I}G}(H_0,Y)$. Choose once and for all unitary representatives
$\pi_X(\phi)$ and $\pi_Y(\phi)$ (there is no choice for
$\pi_0(\phi)$). Define $x^\prime = 
\pi_X(\phi)x\pi_0(\phi)^*$ and $y^\prime=\pi_Y(\phi)y\pi_0(\phi)^*$.
The assignments $x\mapsto x^\prime$, $y\mapsto y^\prime$ give
isomorphisms ${\cal X}_I\rightarrow {\cal X}_J$, ${\cal
Y}_I\rightarrow {\cal Y}_J$ which preserve the inner products since
$\pi_0(\phi)\Omega =\Omega$. Since
$\pi_X(\phi)\pi_X(g)\pi_X(\phi)^*=\pi_X(g\cdot \phi^{-1})$ and 
$\pi_y(\phi)\pi_Y(g)\pi_Y(\phi)^*=\pi_Y(g\cdot \phi^{-1})$ for
$\phi\in SU(1,1)$ and $g\in {\cal L}G$, the map $U_\phi:x\otimes y
\mapsto x^\prime\otimes y^\prime$ extends to a unitary between
$X\boxtimes_I Y$ and $X\boxtimes_J Y$ such that
$U_\phi\pi_I(g)U_\phi^*=\pi_J(g\cdot \phi^{-1})$ for $g\in {\cal
L}_IG\cdot {\cal L}_{I^c}G$.
\vskip .1in
\vskip .1in
\noindent \bf Hilbert space continuity lemma. \it The natural map 
${\cal X}\otimes {\cal Y}\rightarrow X\boxtimes Y$ extends canonically
to continuous maps 
$X\otimes {\cal Y}\rightarrow X\boxtimes Y$ and ${\cal X}\otimes Y\rightarrow
X\boxtimes Y$. In fact $\|\sum x_i\otimes y_i\|^2\le  \|\sum x_ix_i^* \| \, 
\sum \|y_i\Omega\|^2$ and $\|\sum x_i\otimes y_i \|^2 \le \|\sum
y_iy_i^*\| \, \sum \|x_i\Omega\|^2$. 

\vskip .05in
\noindent \bf Proof (cf [24]). \rm If $z=\sum x_i\otimes y_i\in {\cal
X}\otimes {\cal Y}$, then $\sum ((x_i^*x_j) y_i^* y_j
\Omega,\Omega)=\sum y_i^*\pi_Y(x_i^*x_j) y_j\Omega, \Omega)$, 
since $S_{ij}=x_i^*x_j$ lies in $\pi_0({\cal L}_{I^c}G)^\prime
=\pi_0({\cal L}_IG)^{\prime\prime}$. Let
$\eta_j = y_j\Omega$ and $\eta=(\eta_1,\dots,\eta_n) \in H_0^n$. Then
$$\|\sum x_i\otimes y_i\|^2= (\pi_Y(S)\eta,\eta) \le \|S\|\, \|\eta\|^2 =
\|\sum x_ix_i^* \| \, 
\sum \|y_i\Omega\|^2.$$
Here we used the fact that $S=x^*x$ where $x$ is the column vector
with entries $x_i$; this gives $\|S\|=\|x^*x\|=\|xx^*\|=\|\sum
x_ix_i^*\|$. Similarly we can prove that $\|\sum x_i\otimes y_i \|^2
\le \|\sum y_iy_i^*\| \, \sum \|x_i\Omega\|^2$. 
\vskip .1in
\noindent \bf Corollary (associativity of fusion). \it \it There is a
natural unitary isomorphism 
$X\boxtimes (Y\boxtimes Z)\rightarrow (X\boxtimes Y)\boxtimes Z$.
\vskip .05in
\noindent \bf Proof. \rm The assignment $(x\otimes y)\otimes
z\rightarrow x\otimes (y\otimes z)$ makes sense by the lemma and
clearly implements the unitary 
equivalence of bimodules.
\vskip .1in

\noindent \bf 31. Connes fusion with the vector
representation. \rm 
\rm In the previous chapter we proved that if $(a_{ij})$, $(b_{ij})$
are vector primary fields smeared in intervals $I$ and $J$ in
$S^{1}\backslash\{1\}$ with $J$ anticlockwise after $I$, then:
\vskip .1in
\noindent (a) $b_{gf}a_{g_1f}^*=\sum \nu_h a_{hg}^*b_{h g_1}$
with $\nu_h\ne 0$ if $h> g,g_1$ is permissible.
\vskip .05in
\noindent (b) $b_{gf}a_{fh} =\sum \mu_{f_1} a_{gf_1} b_{f_1 h}$ with
$\mu_{f_1}\ne 0$ if $h< f_1<g$.
\vskip .1in
We use these braiding relations to establish
the main technical result required in the computation of
$H_\square\boxtimes H_f$. This 
answers the following natural question. The operator $a_{\square
0}^*a_{\square 0}$ on $H_0$ commutes with ${\cal L}_{I^c}G$, so lies
in $\pi_0({\cal L}_IG)^{\prime\prime}$. Thus, by local equivalence, we
have the right to ask what its image is under the natural isomorphism
$\pi_f:\pi_0({\cal L}_IG)^{\prime\prime} \rightarrow
\pi_f({\cal L}_IG)^{\prime\prime}$.
\vskip .1in
\noindent \bf Theorem (transport formula). \it $\pi_f(a_{\square
0}^*a_{\square 0}) =  
\sum \lambda_g a_{gf}^*a_{gf}$ with $\lambda_g>0$.
\vskip .1in
\noindent \bf Remark. \rm It is possible, using induction or the
braiding computations in [41], to obtain
the precise values of the coefficients. Since the precise numerical
values are not important for us, we have 
preferred a proof which makes it manifest why the right hand side must
have the stated form with strictly positive coefficients
$\lambda_g$. 
\vskip .1in
\noindent \bf Proof. \rm (1) We proceed by induction on $|f|=\sum
f_i$. Suppose that
$\pi_f(a_{\square 0}^*a_{\square 0}) = 
\sum \lambda_g a_{gf}^*a_{gf}$ and 
$\pi_f(b_{\square 0}^*b_{\square 0}) = 
\sum \lambda_g b_{gf}^*b_{gf}$ with $\lambda_g>0$.
Polarising the second identity, we get
$\pi_f(b_{\square 0}^*b^\prime_{\square 0}) = 
\sum \lambda_g b_{gf}^*b_{gf}^\prime$.
In particular if $x\in {\cal L}_J G$, we may take $b^\prime_{ij} =
\pi_i(x)b_{ij} \pi_j(x)^*$ and thus obtain
$$\pi_f(b_{\square 0}^*\pi_\square(x)b_{\square 0}\pi_0(x)^*) = 
\sum \lambda_g b_{gf}^*\pi_g(x)b_{gf}\pi_f(x)^*.$$
Since $\pi_f(\pi_0(x)^*)=\pi_f(x)^*$, we may cancel $\pi_f(x)$ on both
sides to get
$$\pi_f(b_{\square 0}^*\pi_\square(x)b_{\square 0}) = 
\sum \lambda_g b_{gf}^*\pi_g(x)b_{gf}.$$
\noindent (2) Take $x\in {\cal L}_JG$. By the braiding relations and (1), we
have
$$a_{gf}^*\pi_g(b_{\square  0}^*\pi_\square(x)b_{\square 0})
a_{gf} =\pi_f(b_{\square 0}^*\pi_\square(x) b_{\square 0})
a_{gf}^*a_{gf}=
\sum_{g_1} \sum_{h,k} \lambda_{g_1}\nu_h\mu_k b_{g_1 f}^*
a_{hg_1}^*a_{hk}
\pi_k(x)b_{kf}.$$
If $x_i\in {\cal L}_JG$, let $Y=(y_{ij})$ be the operator--valued
matrix with entries $y_{ij}=a_{gf}^*\pi_g(b_{\square
0}^*\pi_\square(x_i^{-1}x_j)b_{\square 0})a_{gf}$. 
Then $Y$ is positive, so that $\sum (y_{ij}\xi_j,\xi_i)\ge 0$ for
$\xi_i\in H_f$. Substituting the expression on the left hand side
above, this gives 
$$\sum_{i,j}\sum_{g_1} \lambda_{g_1} (b_{g_1 f}^*\pi_{g_1}(x_i^{-1})
\left( \sum \nu_h\mu_k
a_{hg_1}^*a_{hk}\right) 
\pi_k(x_i)b_{kf}\xi_j,\xi_i)\ge 0.$$
On the other hand, von Neumann density implies that $\pi({\cal
L}_JG\cdot {\cal 
L}_{J^c}G)^{\prime\prime}=\pi({\cal 
L}G)^{\prime\prime}$ for any positive energy representation at level
$\ell$.  
This implies that vectors of the form $\eta=(\eta_k)$, where
$\eta_k=\pi_k(x) b_{kf}\xi$ with
$\xi\in H_f$ and $x\in {\cal L}_{J}G$, span a dense subset of
$\bigoplus H_k$.  But from the above equation we have
$\sum \lambda_{g_1}\,\nu_h\mu_k (a_{hk}\eta_k ,a_{h
g_1}\eta_{g_1})\ge 0$, and
this inequality holds for all choices of $\eta_k$. In particular,
taking all but one $\eta_{g_1}$ equal to zero, we get $\lambda_{g_1}
\nu_h\mu_{g_1}>0$.  Thus in the expression
$b_{g_1f}a_{gf}^* a_{gf} =\sum_{h,k} \nu_h \mu_k a_{hg_1}^*a_{hk}b_{kf}$,
we have $\nu_h\mu_{g_1}>0$. 
\vskip .05in
\noindent (3) Now for $x \in {\cal L}_JG$, we have
$$b_{g_1 f}^*\pi_{g_1}(a_{\square 0}^*a_{\square 0})
\pi_{g_1}(x) b_{g_1 
f} 
=b_{g_1 f}^*\pi_{g_1}(x)b_{g_1 f} \sum \lambda_g
a_{gf}^*a_{gf}
= \sum \lambda_g\,\nu_h\mu_k b_{g_1 f}^*a_{hg_1}^*
\pi_{h}(x)a_{hk}b_{kf}.$$
If $x_i\in {\cal L}_JG$, let $Z=(z_{ij})$ be the operator--valued
matrix with entries
$z_{ij} =b_{g_1 f}^*\pi_{g_1}(a_{\square 0}^*a_{\square
0})\pi_{g_1}(x_i^{-1}x_j)b_{g_1f}$.
Then $Z$ is positive, so that if $\xi_i\in H_f$, $\sum
(z_{ij}\xi_j,\xi_i)\ge 0$. Let $\eta=(\eta_k)$ where $\eta_k=\sum
\pi_k(x_i)b_{kf}\xi_i$. As above, von Neumann density implies the vectors
$\eta$ are dense in $\bigoplus H_k$. Moreover we have
$$\sum\lambda_g \, \nu_h\mu_k (a_{hk}\eta_k,a_{hg_1}\eta_{g_1}) =
(\pi_{g_1}(a_{\square 0}^*a_{\square 0}) \eta_{g_1} ,\eta_{g_1}).$$
Since this is true for all $\eta_k$'s, all the terms with $k\ne g_1$
must give a zero contribution and
$$(\pi_{g_1}(a_{\square 0}^*a_{\square 0}) \eta_{g_1} ,\eta_{g_1})
=\sum \lambda_g \, \nu_h\mu_{g_1}
(a_{hg_1}\eta_{g_1},a_{hg_1}\eta_{g_1}).$$
But we already saw that $\nu_h\mu_{g_1}>0$ and hence
$\pi_{g_1}(a_{\square 0}^*a_{\square 0})=\sum \Lambda_h
a_{hg_1}^*a_{hg_1}$, with $\Lambda_h>0$, as required.
\vskip .1in
\noindent \bf Corollary. \it If $H_f$ is any irreducible
positive
energy representation of level $\ell$, then as positive energy
bimodules we have
$$H_\square \boxtimes H_f\cong \bigoplus H_g,$$
where $g$ runs over all permissible Young diagrams that can be
obtained by adding a box to $f$. Moreover the action of ${\cal L}_I
G\cdot {\cal 
L}_{I^c}G$ on $H_\square\boxtimes H_f$ extends uniquely to an action
of ${\cal L}G\rtimes {\rm Rot}\, S^1$. 

\vskip .05in
\noindent \bf Proof. \rm Let ${\cal X}_0\subset
{\rm Hom}_{{\cal L}_{I^c}G}
(H_0, H_\square)$ be the linear span of intertwiners $x=\pi_\square(h)
a_{\square 0}$, where $h\in {\cal L}_IG$ and $a$ is a vector primary
field supported in $I$. Since $x\Omega=(\pi_\square(h) a_{\square
0}\pi_0(h)^*) \pi_0(h)\Omega$, it follows from the Reeh--Schlieder
theorem that ${\cal 
X}_0\Omega$ is dense in ${\cal X}_0 H_0$. But then the von
Neumann density argument (for example) implies that ${\cal X}_0\Omega$
is dense in $H_\square$. If $x=\sum
\pi_\square(h^{(j)})a^{(j)}_{\square 0}\in {\cal X}_0$, set  
$x_{gf}=\sum \pi_g(h^{(j)})a^{(j)}_{gf}$. Let $y\in {\rm
Hom}_{{\cal L}_IG}(H_0,H_f)$. By the transport formula
$$(x^*x y^*y\Omega,\Omega)=
(y^* \pi_f(x^*x)y\Omega,\Omega)=\sum_g \lambda_g\,
(x_{gf}^*x_{gf}y\Omega,y\Omega) 
=\sum_g \lambda_g \,\|x_{gf}y\Omega\|^2.$$
This formula shows that $x_{gf}$ is independent of the expression for
$x$. More importantly, by polarising we get an isometry $U$ of the
closure of ${\cal X}_0\otimes {\cal Y}$ in
$H_\square \boxtimes H_f$ into $\bigoplus H_g$, sending $x\otimes y$ to
$\bigoplus \lambda_g^{1/2} x_{gf}y\Omega$. By the Hilbert space
continuity lemma, ${\cal X}_0\otimes {\cal Y}$ is dense in
$H_\square\boxtimes H_f$. Since each of the maps $x_{gf}$ can be
non--zero, Schur's lemma implies that $U$ is surjective and hence a
unitary. The action of ${\cal L}^{\pm 1}G$ extends uniquely to ${\cal
L}G$ by Schur's lemma. The extension to ${\rm Rot}\, S^1$ is uniquely
determined by the fact that ${\rm Rot}\, S^1$ has to fix the lowest
energy subspaces of each irreducible summand of $H_f\boxtimes
H_\square$. 
\vskip .1in

\noindent \bf 32. Connes fusion with exterior
powers of the vector representation. \rm We now extend the methods of
the previous section to the 
exterior powers $\lambda^k V=V_k$. We shall simply write $[k]$ for
the corresponding signature, i.e.~$k$ rows with one box in each. For
$a\in L^2(I,V)$, we shall write $\phi_{gf}(a)$ for
$\phi^\square_{gf}(e_{-\alpha_{gf}}a)$, where $\alpha_{gf}=(\Delta_g-\Delta_f-\Delta_\square)/2(N+\ell)$ is the
phase anomaly introduced in Section~29. For any path $P$ of length
$k$, $f_0<f_1<\cdots <f_k$ with $f_i$ permissible, we define
$a_P=\phi_{f_k f_{k-1}}(a_k) 
\cdots \phi_{f_1f_0}(a_1)$ for $a_i\in L^2(I,V)$. In particular we let $P_0$ be the path 
$0<[1]<[2]<\cdots<[k]$.
\vskip .1in

\noindent \bf Theorem. \it If $a_i,b_i$ are test functions in
$L^2(I,V)$, 
then $$\pi_f(b_{P_0}^*a_{P_0})
=\sum_{g>_k f} (\sum_{P:f\rightarrow g} \lambda_P(g)
b_P)^*(\sum_{P:f\rightarrow g} \lambda_P(g) a_P),$$
where $P$ ranges over all paths $f_0=f<f_1<\cdots <f_k=g$ with each
$f_i$ permissible and where for
fixed $g$ either $\lambda(g)=0$ or $\lambda_P(g)\ne 0$ for all $P$.
\vskip .1in

\noindent \bf Proof. \rm (1) \it The linear span of vectors
$\bigoplus_{f_k>f_{k-1} 
>\cdots >f_1>f} \phi_{f_k 
f_{k-1}}(a_k) \phi_{f_{k-1} 
f_{k-2}}(a_{k-1}) \cdots \phi_{f_1 f}(a_1)\xi$ with $a_j\in
L^2(I_j,V)$ (where $I_j\subseteq I$) and $\xi\in H_f$ is dense in
$\bigoplus_{f_k>f_{k-1} >\cdots >f_1>f} H_{f_k}$.
\vskip .05in
\noindent \bf Proof. \rm  We prove the result by induction on $k$. For
$k=1$, let $H$ denote the closure of this subspace so that $H$ is
invariant under ${\cal L}^{\pm 1}G$ and hence ${\cal L}G$. By Schur's
lemma $H$ must coincide with $\bigoplus_{f_1>f} H_{f_1}$ as required.
By induction the linear span of vectors $\bigoplus_{f_{k-1}
>\cdots >f_1>f} \phi_{f_{k-1} 
f_{k-2}}(a_{k-1}) \cdots \phi_{f_1 f}(a_1)\xi$ with $a_i\in
L^2(I,V)$ and $\xi\in H_f$ is dense in
$\bigoplus_{f_{k-1} >\cdots >f_1>f} H_{f_{k-1}}$. The proof is
completed by noting that if $g$ is fixed and $h_1,\dots ,h_m<g$ (not
necessarily distinct) then the vectors $\bigoplus\phi_{gh_i}(a) \xi_i$
with $a\in L^2(I,V)$ and $\xi_i\in H_{h_i}$ span a dense subspace of
$H_g\otimes {\Bbb C}^m$. Again the closure of the subspace is ${\cal
L}G$ invariant and the result follows by Schur's lemma, because the
$\xi_i$'s vary independently.

\vskip .05in
\noindent (2) \it We have
$$\pi_f(b_{P_0}^*a_{P_0})
=\sum_{g>_k f} \sum_{P,Q:f\rightarrow g} \mu_{PQ}(g)
b_P^*a_Q,$$
where $g$ ranges over all permissible signatures that can be obtained
by adding $k$ boxes to $f$ and $P$, $Q$ range over all permissible
paths $g=f_k>f_{k-1} > \cdots >f_1>f$ and $\mu(g)=(\mu_{PQ}(g))$ is a
non--negative matrix. 
\vskip .05in
\noindent \bf Proof. \rm We assume the result by induction on
$|f|=\sum f_i$. By polarisation, it is 
enough to prove the result with $b_j=a_j$ for all $j$.
If $h>f$, let $x_{hf}=\phi_{hf}(c)$ with $c\in
L^2(I^c,V)$ and $y=a_{P_0}$. Then for $f^\prime >f$ fixed,
$x_{f^\prime f}\pi_f(y^*y) = 
\pi_{f^\prime}(y^*y)x_{f^\prime f}$.  
Substituting for $\pi_f(y^*y)$ and using 
the braiding relations with vector 
primary fields and their duals, $x_{f^\prime f}\pi_f(y^*y)$ can be
rewritten as  
$$x_{f^\prime f}\pi_f(y^*y) =\sum_{g^\prime}\sum_{f_1>f}  \sum_{P,Q}
\mu_{P,Q}(g^\prime) 
a_P^*a_Q 
x_{f_1 f},$$
where $g^\prime$ ranges over signatures obtained by adding $k$ boxes
to $f^\prime$, $P$ ranges over paths $f^\prime<h_1<\cdots
<h_k=g^\prime$ and $Q$ 
ranges over paths $f_1<h_1<\cdots <h_k=g^\prime$. By (1), the vectors
$\bigoplus_{f_1>f} x_{f_1f}H_f$
span a dense subset of $\bigoplus_{f_1>f} H_{f_1}$. Since
$x_{f^\prime f}\pi_f(y^*y) =\pi_{f^\prime}(y^*y)x_{f^\prime f}$, it follows that 
$\pi_{f^\prime}(y^*y) =\sum_{g^\prime}\sum_{f_1>f}  \sum_{P,Q}
\mu_{P,Q}(g^\prime) 
a_P^*a_Q$.
Since $\pi_{f^\prime}(y^*y)$ lies in $B(H_{f^\prime})$, only terms
with $f_1=f^\prime$ appear 
in the above expression so that
$$\pi_{f^\prime}(y^*y) =\sum_{g^\prime}\sum_{P,Q} \mu_{P,Q}(g^\prime)
a_P^*a_Q,$$
where $P$ and $Q$ range over paths from $f^\prime$ to $g^\prime$. Now
suppose $z=y_1+\cdots + y_m$ with $y_i$ having the same form as $y$.
Then 
$$\pi_{f^\prime}(z^*z) =\sum_{g^\prime} \sum_{P,Q}\mu_{P,Q}(g^\prime) 
\sum_{i,j} a_{P,i}^*a_{Q,j}.$$
But $(\pi_{f^\prime}(z^*z)\xi,\xi)\ge 0$ for $\xi\in H_{f^\prime}$ and
the linear span 
of vectors $\bigoplus_Q a_Q \xi$ is dense in $\bigoplus_Q
H_{g^\prime}$. Fixing $g^\prime$, it follows that $\sum
\mu_{P,Q}(g^\prime)(\xi_P,\xi_Q) \ge 0$ for all choices of $\xi_P$ in
$H_{g^\prime}$. Taking all the $\xi_P$'s proportional to a fixed vector
in $H_{g^\prime}$, we deduce that $\mu(g^\prime)$ must be a non--negative
matrix. 
\vskip .05in

\noindent (3) \it If $g>_k f$ is permissible, then $\mu(g)$ has rank
at most one; otherwise
$\mu(g)=0$. If $\mu(g)\ne 0$, then
$\mu_{PQ}(g)=\overline{\lambda_P(g)}\lambda_Q(g)$ with
$\lambda_P(g)\ne 0$ for all $P$.
\vskip .05in
\noindent \bf Proof. \rm We have 
$$\pi_f(b^*a) 
=\sum_{g>_k f} \sum_{P,Q:f\rightarrow g} \mu_{PQ}(g) b_P a_Q^*,$$ 
where $a=a_{P_0}$ and $b=b_{P_0}$. We choose
$a_j$ to be concentrated in disjoint intervals $I_j$ with $I_j$
preceding $I_{j+1}$ going anticlockwise. Fix $i$ and let $a^\prime$,
$a^\prime_Q$ be the intertwiners resulting from swapping $a_i$ and
$a_{i+1}$. Then $a^\prime =\delta_- a$ where $\delta_-\ne 0$
while either $a^\prime_Q=\alpha_Q a_Q +\beta_Q a_{Q_1}$ and
$a^\prime_{Q_1}=\gamma_Q a_Q +\delta_Q a_{Q_1}$, with
$\alpha_Q,\beta_Q,\gamma_Q,\delta_Q\ne 0$, or $a^\prime_Q=\delta_\pm a_Q$. Here if $Q$ is the path $f<f_1<\cdots <f_k=g$,
then $Q_1$ is the other possible path $f<f_1^\prime <\cdots
<f_k^\prime=g$ with $f_j^\prime=f_j$ for $j\ne i$. In the second case,
$\delta_+$ occurs if $f_{i+1}$ is obtained by adding two boxes to the
same row of $f_{i-1}$ while $\delta_-$ occurs if they are added to the same
column. 

Now we still
have
$\pi_f(b^*a^\prime)=\sum_{g>_k f} (\sum_{P:f\rightarrow g} \mu_{PQ}(g)
b_P^*a_P^\prime$.
If $Q$ and $Q_1$ are distinct, it follows that $\delta_- \mu_{PQ}
=\alpha_Q\mu_{PQ}  +\gamma_Q \mu_{PQ_1}$ and 
$\delta_-\mu_{PQ}=\beta_Q \mu_{PQ}
+\delta_{Q_1}\mu_{PQ_1}$ for all $P$. In the case where $Q_1=Q$,
we get $\delta_-\mu_{PQ}=\delta_{\pm}\mu_{PQ}$. 
Now for a vector $(\lambda_Q)$, consider the equations $\delta_-
\lambda_{Q} 
=\alpha_Q\lambda_{Q}  +\gamma_Q \lambda_{Q_1}$ and
$\delta_-\lambda_{Q_1}=\beta_Q \lambda_{Q}
+\delta_{Q}\lambda_{Q_1}$; or $\delta_-\lambda_Q =\delta_{\pm}
\lambda_Q$. These are satisfied when $\lambda_Q=\mu_{PQ}$. We claim
that, if $g>_k f$, these equations have up to a scalar multiple at most
one non--zero solution, 
with all entries non--zero, and otherwise only the
zero solution. This shows 
that $\mu(g)$ has rank at most one with the stated form if $g>_k f$
and $\mu(g)=0$ otherwise. 
 
We shall say that two paths are {\it adjacent} if one is obtained from
the other by changing just one signature. We shall say that two paths
$Q$ and $Q_1$ are {\it connected} if there is a chain of adjacent
paths from $Q$ to $Q_1$. We will show below that any other path $Q_1$ 
from $f$ to $g$ is connected 
to $Q$. This shows on the one hand that if a path $Q$ is obtained 
by successively adding two boxes to  
the same row, we have $\delta_-\lambda_Q=\delta_+\lambda_Q$, so
that $\lambda_Q=0$ since $\delta_+\ne \delta_-$; while on the other
hand if $Q$ and $Q_1$ are adjacent, $\lambda_{Q_1}$ is uniquely 
determined by $\lambda_Q$ and is non--zero if $\lambda_Q$ is.

We complete the proof by showing by induction  
on $k$ that any two paths $f=f_0<f_1<\cdots <f_k=g$ and
$f=f^\prime_0<f_1^\prime <\cdots <f^\prime_k=g$ are connected. The
result is trivial for $k=1$. Suppose the result is known for $k-1$.
Given two paths $f=f_0<f_1<\cdots <f_k=g$ and
$f=f^\prime_0<f_1^\prime <\cdots <f^\prime_k=g$, either
$f_1=f_1^\prime$ or $f_1\ne f_1^\prime$. If $f_1=f_1^\prime=h$, the
result follows because the paths $h=f_1<\cdots <f_k=g$ and
$h=f_1^\prime <\cdots <f^\prime_k=g$ must be connected by the
induction hypothesis. If $f_1\ne f_1^\prime$, there is a unique
signature $f_2^{\prime\prime}$ with
$f_2^{\prime\prime}>f_1,f_1^\prime$. We can then find a path
$f_2^{\prime\prime}<f_3^{\prime\prime}<\cdots <f_k^{\prime\prime}=g$. 
The paths $Q:f<f_1<f_2^{\prime\prime}<\cdots<f_k^{\prime\prime}=g$ and
$Q_1^\prime:f<f_1^\prime<f_2^{\prime\prime}<\cdots<f_k^{\prime\prime}=g$
are adjacent. By induction $Q$ is connected to $Q^\prime$ and $Q_1$ is
connected to $Q_1^\prime$. Hence $Q$ is connected to $Q_1$, as
required. 

\vskip .1in
\noindent \bf Corollary. \it $H_{[k]}\boxtimes H_f= \bigoplus_{g>_k
f,\,\lambda(g)\ne 0} H_g\le \bigoplus_{g>_k f} H_g $.
\vskip .05in

\noindent \bf Proof. \rm If $h\in {\cal L}_IG$, then we have
$$\pi_f(b_{P_0}^* \pi_{[k]}(h) a_{P_0})
=\sum_{g>_k f} (\sum_{P:f\rightarrow g} \lambda_P
b_P)^*\pi_g(h)(\sum_{P:f\rightarrow g} \lambda_P a_P).$$
Now the intertwiners $x=\pi_{[k]}(h)a_{P_0}$ span a subspace ${\cal
X}_0$ of 
${\rm Hom}_{{\cal L}_{I^c}}(H_0,H_{[k]})$. As before the transport
formula shows 
that the assignment $x\otimes y
\mapsto \bigoplus_g \sum \lambda_P(g) \pi_g(h)
a_P y\Omega$ extends to a 
linear isometry  $T$ of ${\cal X}_0\otimes {\cal Y}$ into
$\bigoplus_{\lambda(g)\ne 0}
H_g$. $T$ intertwines ${\cal L}^{\pm 1}G$, so
by Schur's lemma 
extends to an isometry of the closure of ${\cal X}_0\otimes {\cal Y}$
in $H_{[k]}\boxtimes H_f$ onto $\bigoplus_{\lambda(g)\ne 0} H_g$. 
On the other hand, by the argument used in the corollary in the
previous section, ${\cal X}_0\Omega$ is dense in $H_{[k]}$.
Therefore, by the Hilbert  
space continuity lemma, the image of ${\cal X}_0\otimes 
{\cal Y}$  is dense in $H_{[k]}\boxtimes H_f$. Hence $H_{[k]}\boxtimes
H_f =\bigoplus_{\lambda(g)\ne 0} H_g$, as required.

\vskip .1in
\noindent \bf 33. The fusion ring. \rm Our aim now is to show 
that if $H_i$ and $H_j$ are irreducible positive energy
representations, then $H_i\boxtimes H_j = \bigoplus N_{ij}^k H_k$ where
the fusion coefficients $N_{ij}^k$ are finite and to be determined.
\vskip .1in
\noindent \bf Lemma (closure under fusion). \it (1) Each irreducible
positive 
energy representation $H_i$ appears in 
some $H_\square^{\boxtimes n}$.

\noindent (2) The $H_i$'s are closed under Connes fusion.

\vskip .05in
\noindent \bf Proof. \rm (1) Since $H_f\boxtimes H_\square =\bigoplus
H_g$, it follows easily by induction that each $H_g$ is contained in
$H_{\square}^{\boxtimes m}$ for some $m$.
\vskip .05in
\noindent (2) Since $H_f\subset H_{\square}^{\boxtimes m}$ for some
$m$ and $H_g\subset H_{\square}^{\boxtimes n}$ for some $n$, we have
$H_f\boxtimes H_g\subset H_\square^{\boxtimes (m+n)}$. By induction
$H_\square^{\boxtimes k}$ is a direct sum of irreducible positive energy
bimodules. By Schur's lemma any subrepresentation of
$H_\square^{\boxtimes (m+n)}$ must be a direct sum of irreducible
positive energy bimodules. In particular this applies to 
$H_f\boxtimes H_g$, as required.
\vskip .1in
\bf \noindent Corollary. \it If $X$ and $Y$ are positive energy
representations, the action of ${\cal L}_I G\cdot {\cal
L}_{I^c}G$ on $X\boxtimes Y$ extends uniquely to an action
of ${\cal L}G\rtimes {\rm Rot}\, S^1$. 
\vskip .05in
\noindent \bf Proof. \rm The action extends uniquely to ${\cal L}G$ by
Schur's lemma. The extension to ${\rm Rot}\, S^1$ is uniquely
determined by the fact that ${\rm Rot}\, S^1$ has to fix the lowest
energy subspaces of each irreducible summand of $X\boxtimes
Y$. 
\vskip .1in

\noindent \bf Braiding lemma. \it The map $B:{\cal X}\otimes {\cal
Y}\rightarrow Y\boxtimes X$, $B(x\otimes y) =
R_\pi^*[R_\pi(y)R_\pi^*\otimes R_\pi(x)R_\pi^*]$ extends to a unitary of
$X\boxtimes Y$ onto $Y\boxtimes X$ intertwining the actions of ${\cal
L}G$. 
\vskip .05in
\noindent \bf Proof. \rm Note that the $B$ is well--defined, for
rotation through $\pi$ interchanges ${\cal L}_I G$ and ${\cal L}_{I^c}G$. Hence
$R_\pi xR_\pi^*$ lies in ${\rm Hom}_{{\cal L}_I G}(H_0,X)$ and $R_\pi y
R_\pi^*$ lies in ${\rm Hom}_{{\cal L}_{I^c}G}(H_0,Y)$. So $R_\pi
yR_\pi^* \otimes R_\pi xR_\pi^*$ lies in ${\cal Y}\otimes {\cal X}$.
Since $R_\pi \Omega=\Omega$, the map $B$ preserves the inner product.
Interchanging the r\^oles of $X$ and $Y$, we get an inverse of $B$
which also preserves the inner product. Hence $B$ extends by
continuity to a unitary of $X\boxtimes Y$ onto $Y\boxtimes X$.
Finally, we check that $B$ has the correct intertwining property.
Let $g\in {\cal L}_I G$ and $h\in {\cal
L}_{I^c}G$. Then 
$$\eqalign{B(gx\otimes hy)& =R_{\pi}^*[R_\pi(hy)R_\pi^* \otimes
R_\pi(gx)R_\pi^*]
=R_\pi^*[(h\circ r_\pi)(g\circ r_\pi) (R_\pi y R_\pi^*\otimes R_\pi x
R_\pi^*)] \cr
&=R_\pi^* (h\circ r_\pi)(g\circ r_\pi) R_\pi^* R_\pi [R_\pi y
R_\pi^*\otimes R_\pi xR_\pi^*]
=gh R_\pi^* [R_\pi y
R_\pi^*\otimes R_\pi xR_\pi^*]
=gh B(x\otimes y),\cr}$$
as required.
\vskip .1in
\noindent \bf Corollary~1. \it $X\boxtimes Y$ is isomorphic to
$Y\boxtimes X$ as an ${\cal L}G$--module.

\vskip .1in

\noindent \rm Let ${\cal R}$ be the representation ring of formal 
sums $\sum m_i H_i$ ($m_i\in {\Bbb Z}$) with multiplication extending
fusion. ${\cal R}$ is called the {\it fusion ring} (at level $\ell$). 
\vskip .1in
\noindent \bf Corollary~2. \it The fusion ring
${\cal R}$ is a commutative ring with an identity. 
\vskip .05in
\noindent \bf Proof. \rm ${\cal R}$ is commutative by the braiding
lemma and closed under multiplication by the previous lemmas.
Multiplication is associative because fusion is. 
\vskip .2in 
\noindent \bf 34. The general fusion rules (Verlinde
formulas). \rm In order to determine the general coefficients
$N_{ij}^k$ in the 
fusion rules  $H_i\boxtimes H_j = \bigoplus N_{ij}^k H_k$, we first have to
determine the structure of the fusion ring. Before doing so, we will
need some 
elementary facts on the affine Weyl group. The integer lattice
$\Lambda={\Bbb Z}^N$ acts by translation on ${\Bbb R}^n$. The symmetric
group $S_N$ acts on ${\Bbb R}^N$ by permuting the coordinates and
normalises $\Lambda$, so we get an action of the semidirect product
$\Lambda\rtimes S_N$. The subgroup $\Lambda_0=\{(N+\ell)(m_i):\sum
m_i=0\}\subset \Lambda$ is 
invariant under $S_N$, so we can consider the semidirect product
$W=\Lambda_0\rtimes S_N$. The sign of a permutation defines a homomorphism
$\det$ of $S_N$, and hence $W$, into $\{\pm 1\}$. 
\vskip .1in
\noindent \bf Lemma. \it (a) $\{(x_i):|x_i-x_j|\le N+\ell\}$ forms a
fundamental domain for the action of $\Lambda_0$ on ${\Bbb R}^N$.

\noindent (b) $D=\{(x_i):x_1\ge \cdots \ge x_N,\, x_1-x_N\le N+\ell\}$
forms a fundamental domain for the action of $\Lambda_0\rtimes S_N$ on
${\Bbb R}^N$.

\noindent (c) A point is in the orbit of the interior of $D$ consists
of points iff its stabiliser is trivial. For every other point $x$
there is an transposition $\sigma\in S_N$ such that $\sigma(x)-x$ lies in
$\Lambda_0$. 
\vskip .05in
\noindent \bf Proof. \rm (a) Take $(x_i)\in {\Bbb R}^N$. Write
$x_i=a_i +m_i$ with $0\le a_i <N+\ell$ and $m_i\in (N+\ell){\Bbb Z}$.
Without loss of generality, we may assume that $a_1\le \cdots \le
a_N$. Now $(m_i)$ can be written as the sum of a term 
$(b_i)=(N+\ell)(M,M,\dots,M,M-1,M-1,\dots,M-1)$ and an element $(c_i)$
of $\Lambda_0$. Thus $x=a+b +c$ with $c\in \Lambda_0$. It is easy to
see that $y=a+b$ satisfies $|y_i-y_j|\le N+\ell$. (b) follows
immediately from (a) since the domain there is 
invariant under $S_N$. Finally, since ${\rm int}(D)=
\{(x_i):x_1>\cdots > x_N,\, x_1-x_N< N+\ell\}$, it is easy to see that
any point in ${\rm int}(D)$ has trivial stabiliser. If
$x\in \partial D$, then either $x_i=x_{i+1}$ for some $i$, in which
case $\sigma=(i,i+1)$ fixes $x$; or $x_1-x_N=N+\ell$, in which case
$\sigma=(1,N)$ satisfies $\sigma(x)-x=(-N-\ell,0,\dots,0,N+\ell)$. 
Thus (c) follows for points in $D$ and therefore in general, since $D$
is a fundamental domain. 
\vskip .1in 
\bf \noindent Corollary. \it Let $\delta=(N-1,N-2,\dots,1,0)$. Then $m\in
{\Bbb Z}^N$ has trivial stabiliser in $W=\Lambda_0\rtimes S_N$ iff
$m=\sigma(f+\delta)$ 
for a unique $\sigma\in W$ and signature $f_1\ge f_2\ge \cdots f_N$ with
$f_1-f_N\le \ell$; $m$ has non--trivial stabiliser iff there is a
transposition $\sigma\in S_N$ such that $\sigma(m)-m$ lies in $\Lambda_0$.

\vskip .05in
\noindent \bf Proof. \rm In the first case $m=\sigma(x)$ for
$\sigma\in W$ and $x\in {\Bbb R}^N$ with 
$x_1>\dots >x_N$ and $x_1-x_N<N+\ell$. Since the $x_i$'s must be
integers, we can write $x=f+\delta$ with $f_1\ge \cdots \ge f_N$. Then
$f_1-f_N=x_1-x_N-(N-1)<\ell +1$, so that $f_1-f_N\le \ell$. 
\vskip .1in
Recall that the character
of $V_f$ is given by Weyl's character formula
$\chi_f(z)={\rm det}(z_i^{f_j+\delta_j})/{\rm
det}(z_i^{\delta_j})$. Let $S$ be the space of permitted
(normalised) signatures at 
level $\ell$, i.e.~$S=\{h:h_1\ge \cdots \ge h_N,\, h_1-h_N\le \ell,\,
h_N=0\}$. We now define a ring ${\cal S}$ as follows. For
$h\in S$, let $D(h)\in SU(N)$ be the diagonal matrix with
$D(h)_{kk}=\exp(2\pi i (h_k + N-k -H) /(N+\ell))$ where $H=(\sum h_k
+N-k)/N$ and set ${\cal T}=\{D(h):h\in S\}$. We denote the set of
functions on ${\cal T}$ by ${\Bbb C}^{\cal T}$. Let
$\theta:R(SU(N))\rightarrow {\Bbb 
C}^{\cal T}$ be the map of restriction of characters, i.e.~$\theta([V])=
\chi_V|_{{\cal T}}$. By definition $\theta$ is a ring *--homomorphism.
Set ${\cal S}=\theta(R(SU(N)))$ and let $\theta_f=\theta(V_f)$.

\vskip .1in
\noindent \bf Proposition. \it (1) $X_{\sigma
(f+\delta)-\delta}|_{\cal T}={\rm 
det}(\sigma) \, X_f|_{\cal T}$ for $\sigma\in S_N$ and $X_{f+m}|_{\cal
T} = X_f|_{\cal T}$ for $m\in \Lambda_0$. 
 
\noindent (2) The $\theta_f$'s with $f$ permissible
form a ${\Bbb Z}$--basis of ${\cal S}$. 

\noindent (3) ${\rm ker}(\theta)$
is the ideal in $R(SU(N))$ generated by $V_f$ with $f_1-f_N=\ell+1$.

\noindent (4) If $V_f\otimes V_{g} = \bigoplus N_{fg}^h V_h$, then 
$\theta_f \theta_{g}= \sum N_{fg}^h \det(\sigma_h)
\theta_{h^\prime}$ where $h$ ranges over those signatures in the
classical rule for which 
there is a $\sigma_h\in \Lambda_0\rtimes S_N$ (necessarily unique)
such that $h^\prime =\sigma_h(h+\delta)-\delta$ is permissible.

\noindent (5) If $f$, $h$ are permissible, then 
$|\{g_1: \hbox{$g_1$ permissible, $f<g_1<_k h$}\}|
=|\{g_2: \hbox{$g_2$ permissible, $f<_kg_2< h$}\}|$. 
\vskip .05in

\noindent \bf Proof. \rm The statements in (1) follow immediately from
the form of the $D(h)$'s. The $V_{[k]}$'s generate $R(SU(N))$ and,
if $f_1-f_N=\ell+1$, it is easy to see that $\chi_f(t)=0$ for all
$t\in {\cal T}$: for $f_1+N-1 -f_N=N+\ell$ and hence the numerator in
$\chi_f(t)$ must vanish. The $\theta_f$'s with $f$ permissible are
therefore closed under multiplication by $\theta_{[k]}$'s. Since 
the $\theta_{[k]}$'s  generate ${\cal S}$, the ${\Bbb Z}$--linear span of the 
$\theta_f$'s with $f$ permissible must equal ${\cal S}$. The characters
$\chi_{[k]}$ distinguish the points of ${\cal T}$ and $\chi_{[0]}=1$.
Hence ${\cal S}_{\Bbb C}$ is a unital subalgebra of ${\Bbb C}^{\cal
T}$ separating points. So given $x,y\in {\cal T}$, we can find $f\in
{\cal S}_{\Bbb C}$ such that $f(x)=1$ and $f(y)=0$. Multiplying these
together for all $y\ne x$, it follows that ${\cal S}_{\Bbb C}$
contains $\delta_x$ and hence coincides with ${\Bbb C}^{\cal T}$. So 
the $\theta_f$'s must be linearly independent over ${\Bbb C}$, so a
fortiori over ${\Bbb Z}$. This proves (2). 
Let $I\subset R(SU(N))$ be the ideal generated by the $[V_g]$'s with
$g_1-g_N=\ell+1$. Since $R(SU(N))$ is generated by the $V_{[k]}$'s and
we have the tensor
product rule $V_f\otimes V_{[k]}=\bigoplus_{g>_kf} V_g$, it follows
that $R(SU(N))/I$ is spanned by the image of the $[V_f]$'s as a ${\Bbb
Z}$--module. But $I\subseteq {\rm ker}(\theta)$ and the
$\theta([V_f])$'s are linearly independent over ${\Bbb Z}$. Hence the
images of the $[V_f]$'s give a ${\Bbb Z}$--basis of $R(SU(N))/I$ and
therefore $I={\rm ker}(\theta)$, so (3) follows.
The assertion in (4) follows from (1) by applying $\theta$ and using
the corollary to the lemma above. In fact, if $h+\delta$ has non--trivial
stabiliser, we can find a transposition $\sigma\in S_N$ such that
$\sigma( h+\delta)-h-\delta$ lies in $\Lambda_0$. Hence 
$X_h(t)=-X_{\sigma(h+delta)-\delta}(t) =-X_h(t)$, so that
$\chi_h(t)=X_h(t)=0$ for all $t\in {\cal T}$. When the stabiliser is
trivial, we clearly have $\theta_h=\det(\sigma_h) \theta_{h^\prime}$.
Finally (5) follows by comparing coefficients of $\theta_h$ in 
$\theta_f\theta_{[k]} \theta_\square= \sum_{g_1>_k f} \sum_{h>g_1} \theta_h
=\sum_{g_2>f} \sum_{h>_k g_2} \theta_h$.

\vskip .1in

\bf \noindent Theorem. \it (1) $H_{[k]}\boxtimes H_f=\bigoplus_{g>_k
f} H_g$, where the sum is over permissible $g$.

\noindent (2) The ${\Bbb Z}$--linear map ${\rm ch}:{\cal
R}\rightarrow {\cal S}$ defined by ${\rm ch}(H_f)=\chi_f|_{{\cal T}}$
is a ring isomorphism.

\noindent (3) The characters of ${\cal
R}$ are given by $[H_f]\mapsto {\rm ch}(H_f,h)=\chi_f(z)$
for $z\in {\cal T}$.

\noindent (4) The fusion coefficients $N_{ij}^k$'s can
be computed using the 
multiplication rules for the basis ${\rm ch}(H_f)$ of ${\cal S}$.

\noindent (5) Each representation $H_f$ has a unique conjugate representation
$\overline{H_f}$ such that $H_f\boxtimes \overline{H_f}$ contains $H_0$.
In fact $\overline{H_f}=H_{f^\prime}$, where $f^\prime_i=-f_{N-i+1}$, and $H_0$
appears exactly once in $H_f\boxtimes H_{f^\prime}$. The map
$H_f\mapsto 
\overline{H_f}$ makes ${\cal R}$ into an involutive ring and ${\rm
ch}$ becomes a *--isomorphism.
\vskip .05in
\noindent \bf Proof. \rm (1) We know that $H_f\boxtimes H_{[k]} \le
\bigoplus_{g>_k f} H_g$ with equality when $k=1$. We prove
by induction on $|f|=\sum f_j$ that 
$H_f\boxtimes H_{[k]} = \bigoplus_{g_1>_k f} H_{g_1}$. It suffices to show
that if this holds for $f$ then it holds for all 
$g$ with $g>f$. Tensoring by $H_\square$ and using part (5) of the
preceding proposition, 
we get 
$$\bigoplus_{g>f} H_g\boxtimes H_{[k]} =\bigoplus_{g_1>_k f}
\bigoplus_{h>g_1} H_h=\bigoplus_{g>f} \bigoplus_{h>_k g}
H_h.$$
Since $ H_g\boxtimes H_{[k]} \le \bigoplus_{h>_k g} H_h$, we must have
equality for all $g$, completing the induction.

\noindent (2) Let ${\rm ch}$ be the ${\Bbb Z}$--linear isomorphism ${\rm
ch}:{\cal R}\rightarrow {\cal S}$ extending ${\rm ch}(H_f)=\theta_f$.
Then by (1), ${\rm ch}(H_{[k]}\boxtimes H_f) =\theta_{[k]} \theta_f$.
This implies that ${\rm 
ch}$ restricts to a ring homomorphism on the subring of ${\cal R}$
generated by the $H_{[k]}$'s. On the other hand the $\theta_{[k]}$'s
generate ${\cal S}$, so the image of this subring is the whole of
${\cal S}$. Since ${\rm ch}$ is injective, the ring generated by the
$H_{[k]}$'s must be the whole of ${\cal R}$ and ${\rm ch}$ is thus
a ring homomorphism, as required. 

\noindent (3) and (4) follow immediately from the isomorphism ${\rm
ch}$ and the fact that ${\cal S}_{\Bbb C}={\Bbb C}^{\cal T}$. 

\noindent (5) We put an inner product on ${\cal S}_{\Bbb C}={\cal
R}_{\Bbb C}$ by taking $\theta_f$ as an orthonormal basis, so that
$(\theta_f,\theta_g)=\delta_{fg}$. We claim that
$(\theta_f\theta_g,\theta_h) =
(\theta_g,\overline{\theta_f}\theta_h)$ for all $\theta_f$. Note that
$\overline{\theta_f}=\theta_{f^\prime}$ where $f_i^\prime=-f_{N-i+1}$.
Let $\theta_f^*$ be the adjoint of multiplication by $\theta_f$. The
multiplication rules for $\theta_{[k]}$ imply that
$\theta_{[k]}^*=\overline{\theta_{[k]}}$ for $k=1,\dots, N$. Thus the
homomorphism $\theta_f \mapsto \overline{\theta_f}^*$ is the identity
on a set of generators of 
${\cal S}$ and therefore on the whole of ${\cal S}$, so the claim
follows. In particular 
$(\theta_f\theta_g,\theta_0)=(\theta_g,\overline{\theta_f})
=(\theta_g,\theta_{f^\prime})=\delta_{gf^\prime}$. Translating to
${\cal R}$, we retrieve all the assertions in (5). 
 
\vskip .1in
The following results are immediate consequences of the theorem and
preceding proposition.
\vskip .1in

\noindent \bf Corollary~1 (Verlinde formulas [38],[20]). \it  If the
``classical'' tensor 
product rules for $SU(N)$ are given by $V_f\otimes V_g=\bigoplus N_{fg}^h
V_h$, then the ``quantum'' fusion rules for $LSU(N)$ at level $\ell$
are given by $$H_f\boxtimes H_g= \bigoplus N_{fg}^h \det(\sigma_h)
H_{h^\prime},$$
where $h$ ranges over those signatures in the classical rule for which
there is a $\sigma_h \in \Lambda_0\rtimes S_N$ (necessarily unique)
such that $h^\prime =\sigma_h(h+\delta)-\delta$ is permissible. 
\vskip .1in
\noindent \bf Corollary~2 (Segal--Goodman--Wenzl rule [33],[13]). \it
The map 
$V_f\mapsto H_f$ extends to a *--homomorphism of $R(SU(N))$ (the
representation ring of $SU(N)$) onto the fusion ring ${\cal R}$. The
kernel of this homomorphism is the ideal generated by the
(non--permissible) representations $V_f$ with $f_1-f_N=\ell+1$.

\vskip .3in

\noindent \bf REFERENCES. \rm
\vskip .1in

\noindent 1. H.~Araki, {\it von Neumann algebras of local observables
for free scalar field}, J.~Math.~Phys. {\bf 5} (1964), 1--13.

\noindent 2. H.~Araki, {\it On quasi--free states of CAR and
Bogoliubov automorphisms}, Publ.~R.I.M.S. {\bf 6} (1970),
385--442.

\noindent 3. J.~Baez, I.~Segal and Z.~Zhou, ``Introduction to
algebraic and   constructive quantum field theory'', Princeton
University Press 1992.

\noindent 4. J.~Bisognano and E.~Wichmann, {\it On the duality
condition for a hermitian scalar field}, J.~Math.~Phys. {\bf 16}
(1975), 985--1007.

\noindent 5. R.~Borcherds, {\it Vertex algebras, Kac--Moody algebras and
the Monster}, Proc. Nat. Acad. Sci. U.S.A. {\bf 83} (1986), 3068--3071.

\noindent 6. H.~Borchers, {\it A remark on a theorem of Misra},
Commun.~Math.~Phys. {\bf 4} (1967), 315--323.

\noindent 7. H.~Borchers, {\it The CPT theorem in two--dimensional
theories of local observables}, Commun.~Math.~Phys. {\bf 143} (1992),
315--323.

\noindent 7. O.~Bratteli and D.~Robinson, ``Operator algebras and
quantum statistical mechanics'', Vol.~I, Springer 1979.

\noindent 8. A.~Connes, ``Non--commutative geometry'' (Chapter~V,
Appendix~B), Academic Press 1994.

\noindent 9. S.~Doplicher, R.~Haag and J.~Roberts, {\it Local observables
and particle statistics I, II}, Comm.~Math.~Phys. {\bf 23} (1971)
119--230 and {\bf 35} (1974) 49--85.

\noindent 10. J.~Glimm and A.~Jaffe, ``Quantum Physics'', 2nd edition,
Springer--Verlag 1987.

\noindent 11. P.~Goddard and D.~Olive, eds., ``Kac--Moody and Virasoro
algebras'', Advanced Series in Math.~Physics Vol.~3, World Scientific
1988.

\noindent 12. R.~Goodman and N.~Wallach, {\it Structure and unitary
cocycle representations of loop groups and the group of
diffeomorphisms of the circle}, J.~Reine Angew.~Math 347 (1984)
69--133. 

\noindent 13. F.~Goodman and H.~Wenzl, {\it Littlewood Richardson
coefficients for Hecke algebras at roots of unity}, Adv.~Math. {\bf
82} (1990), 244--265. 

\noindent 14. R.~Haag, ``Local quantum physics'', Springer--Verlag 1992.

\noindent 15. N.~Hugenholtz and J.~Wierenga, {\it On locally normal
states in quantum statistical mechanics}, Commun.~Math.~Phys. {\bf 11}
(1969), 183--197. 

\noindent 16. E.~Ince, ``Ordinary differential equations'', Dover 1956.

\noindent 17. V.~Jones, {\it Index for subfactors}, Invent.~math. {\bf 72}
(1983), 1--25.

\noindent 18. V.~Jones, {\it von Neumann algebras in mathematics and
physics}, Proc.~I.C.M. Kyoto 1990, 121--138.

\noindent 19. R.~Jost, ``The general theory of quantised fields'',
A.M.S.~1965.

\noindent 20. V.~Kac, ``Infinite dimensional Lie algebras'', 3rd edition,
C.U.P. 1990.

\noindent 21. D.~Kazhdan and G.~Lusztig, {\it Tensor structures arising
from affine Lie algebras IV}, Journal A.M.S. {\bf 7} (1994), 383--453.

\noindent 22. V.~Knizhnik and A.~Zamolodchikov, {\it Current algebra and
Wess--Zumino models in two dimensions}, Nuc.~Phys.~B {\bf 247} (1984),
83--103.

\noindent 23. P.~Leyland, J.~Roberts and D.~Testard, {\it Duality for
quantum free fields}, preprint, Luminy 1978. 

\noindent 24. T.~Loke, ``Operator algebras and conformal field theory
for the discrete series representations of ${\rm Diff}\, S^1$'',
thesis, Cambridge 1994. 

\noindent 25. F.~Murray and J.~von Neumann, {\it On rings of
operators}, Ann.~Math. {\bf 37} (1936), 116--229.

\noindent 26. T.~Nakanishi and A.~Tsuchiya, {\it Level--rank duality of
WZW models in conformal field theory}, Comm. Math. Phys. {\bf 144}
(1992), 351--372.

\noindent 27. E.~Nelson, ``Topics in dynamics I: Flows'', Princeton
University Press 1969.

\noindent 28. S.~Popa, ``Classification of subfactors and their
endomorphisms'', C.B.M.S. lectures, A.M.S.~1995.

\noindent 28. A.~Pressley and G. Segal, ``Loop groups'', O.U.P.~1986.

\noindent 29. M.~Reed and  B.~Simon, ``Methods of Mathematical Physics
I: Functional analysis'', Academic Press 1980.

\noindent 30. H.~Reeh and S.~Schlieder, {\it Bemerkungen zur
Unitarequivalenz von Lorentzinvarienten Feldern}, Nuovo Cimento {\bf
22} (1961), 1051--1068.

\noindent 31. M.~Rieffel and A.~van Daele, {\it A bounded operator
approach to Tomita--Takesaki theory}, Pacific J.~Math. {\bf 69}
(1977), 187--221.

\noindent 32. J.--L.~Sauvageot, {\it Sur le produit tensoriel relatif
d'espaces d'Hilbert}, J.~Op.~Theory {\bf 9} (1983), 237--252.

\noindent 33. G.~Segal, {\it Notes on conformal field theory}, unpublished
manuscript.

\noindent 34. R.~Streater and A.~Wightman, ``PCT, spin statistics and
all that'', Benjamin 1964.

\noindent 35. M.~Takesaki, {\it Conditional expectation in von Neumann
algebra}, J.~Funct.~Analysis {\bf 9} (1972), 306--321.

\noindent 36. A.~Thomae, ``Ueber die h\"oheren hypergeometrischen
Reihen, ...'' Math.~Ann. {\bf 2} (1870), 427--444. 

\noindent 37. A.~Tsuchiya and Y.~Kanie, {\it Vertex operators in conformal
field theory on ${\Bbb P}^1$ and monodromy representations of braid
group}, Adv.~Studies in Pure Math. {\bf 16} (1988), 297--372.

\noindent 38. E.~Verlinde, {\it Fusion rules and modular
transformations in 2D conformal field theory}, Nuclear Phys.~B {\bf
300} (1988), 360--376.

\noindent 39. A.~Wassermann, {\it Operator algebras and conformal
field theory}, Proc.~I.C.M. Zurich (1994), 966-979, Birkh\"auser.

\noindent 40. A.~Wassermann, with contributions by V.~Jones,
``Lectures on operator algebras and conformal field theory'', 
Proceedings of Borel Seminar, Bern 1994, to appear.

\noindent 41. A.~Wassermann, {\it Operator algebras and conformal
field theory~II: Fusion for von Neumann algebras and loop groups}, to
appear. 

\noindent 42. A.~Wassermann, {\it Operator algebras and conformal
field theory~IV: Loop groups, quantum invariant theory and
subfactors}, in preparation. 

\noindent 43. H.~Wenzl, {\it Hecke algebras of type $A_n$ and subfactors},
Invent.~math. {\bf 92} (1988), 249--383.

\noindent 44. H.~Weyl, ``The classical groups'', Princeton University
Press 1946.

\noindent 45. E.~Whittaker and G.~Watson, ``A course of modern analysis'',
4th edition, C.U.P.~1927.

\noindent 46. V.~Jones, {\it Fusion en alg\`ebres de von Neumann et
groupes de lacets (d'apr\`es A.~Wassermann), S\'eminaire Bourbaki
1994--95, No.~800. 

\vfill\eject
\end